\documentclass[times]{article}

\usepackage{amsmath,amssymb}
\usepackage{graphicx}
\usepackage{subfig}
\usepackage{stmaryrd}
\usepackage{bm}
\usepackage{epstopdf}

\date{June 24, 2016}

\title{CutFEM topology optimization of \\3D laminar incompressible flow problems}

\author{Hernan Villanueva and Kurt Maute \\ Department of Aerospace Engineering, \\ University of Colorado Boulder, 
 Boulder, CO 80309-0429, USA}

\begin{document}

\maketitle

\begin{abstract}
This paper studies the characteristics and applicability of the CutFEM approach \cite{BCH+:14} as the core of a robust topology optimization framework for 3D laminar incompressible flow and species transport problems at low Reynolds number ($Re < 200$). CutFEM is a methodology for discretizing partial differential equations on complex geometries by immersed boundary techniques. In this study, the geometry of the fluid domain is described by an explicit level set method, where the parameters of a level set function are defined as functions of the optimization variables. The fluid behavior is modeled by the incompressible Navier-Stokes equations. Species transport is modeled by an advection-diffusion equation. The governing equations are discretized in space by a generalized extended finite element method. Face-oriented ghost-penalty terms are added for stability reasons and to improve the conditioning of the system. The boundary conditions are enforced weakly via Nit\-sc\-he's method. The emergence of isolated volumes of fluid surrounded by solid during the optimization process leads to a singular analysis problem. An auxiliary indicator field is modeled to identify these volumes and to impose a constraint on the average pressure. Numerical results for 3D, steady-state and transient problems demonstrate that the CutFEM analyses are sufficiently accurate, and the optimized designs agree well with results from prior studies solved in 2D or by density approaches.
\end{abstract}


\section{Introduction}
\label{sec:intro}

The performance of a broad range of engineering applications is dependent on the characteristics of internal and external flows. Such applications include flows for momentum, energy and species transport. The performance of these systems can be improved by finding the optimal geometry of the fluid-solid interface. The increase of computational power and the development of improved numerical schemes allow for the systematic design of flow problems via mathematical optimization methods. In this paper, we study the influence of conceptual changes to the design geometry of fluid flow problems within the context of topology optimization. 

Initial ventures into fluid optimization focused on shape optimization \cite{MP:01,Gunzburger:03}. Topology optimization of fluid problems was pioneered by \cite{BP:03}, who adopted the concept of density methods to Stokes flows. Density methods were originally developed for problems in solid mechanics \cite{Bendsoe:89,ZR:91}, and describe the geometry of a body by the density distribution of a fictitious two-phase material. The density, or volume fraction, varies continuously to describe a smooth transition between two distinct phases, e.g.~fluid and solid. Fluid topology optimization was extended to a Darcy-Stokes flow model by \cite{GP:06}, and to the Navier-Stokes equations by \cite{GSH:05} and \cite{OOB:06}, among others.

Common to most density methods in fluid topology optimization is the use of a Brinkman model \cite{ABF:99} to weakly enforce no-slip conditions at the fluid-solid interface. This approach may suffer from errors if the optimization process converges to a ``0-1'' material distribution and the fluid-solid interface is not aligned with the mesh. In this case, the interface is approximated by a stair-step pattern that hinders the precise determination of the optimum geometry and may affect the accuracy of the flow solution. This issue is more pronounced for coarse meshes and high Reynolds numbers \cite{KM:11}, and can be mitigated by mesh refinement \cite{MR:95,MR:97,YNK+:11}. However, the required level of mesh refinement may lead to high computational costs. Adaptive mesh refinement may affect the convergence of the optimization process if a gradient-based optimization algorithm is applied \cite{SMR:00}.

The shortcomings of density methods have promoted the development of Level Set Methods (LSMs) for topology optimization, c.f.~\cite{PDR:12,DML+:13}. External and internal phase boundaries are described implicitly by the zero level set isosurfaces of a Level Set Function (LSF), $\phi\left(\bm{x}\right)$, where $\bm{x}$ is the position vector \cite{SW:00,WWG:03,AJT:04}. LSMs are well-suited for topology optimization because smooth differentiable changes in the LSFs lead to changes in the topology of the body, such as holes merging or splitting \cite{OS:88,AJT:02,WW:04c}; however, these changes may lead to discontinuities in the physical response \cite{JM:15}.

The LSF is typically discretized on a fixed background mesh and updated in the optimization process by solving the Hamilton-Jacobi equations \cite{WWG:03,AJT:04}. An alternate approach that is utilized in this work is to define the parameters of the discretized LSF as explicit functions of the optimization variables. The resulting parameter optimization problem is solved by standard nonlinear programming (NLP) methods, which provide flexibility with respect to the class of optimization problems one can model \cite{DML+:13}. In the context of flow topology optimization, NLP schemes allow, for example, including additional optimization variables in addition to the parameters of the discretized LSF, such as the position and shape of the inlets and outlets. For a detailed discussion of the LSM, the reader is referred to the comprehensive review by \cite{DML+:13} and \cite{GP:13}.

Several methods exist to represent the geometry and the material distribution described by the LSF in the mechanical model. The Ersatz material method \cite{WWG:03,AJ+:05} interpolates the physical properties of a fictitious material by using either an element-wise constant volume fraction or the LSF value at a point \cite{YYK+:15}. In a fluid-solid problem, the Ersatz material approach defines typically the local permeability of the Brinkman model. While the Ersatz material approach eases the computational complexity, the method faces the same issues as density methods in regards to enforcing boundary conditions.

In this work, we utilize the eXtended Finite Element Method (XFEM) to describe the material distribution in the mechanical model. The XFEM is an immersed boundary technique that does not require a mesh that conforms to the phase boundaries. The method was built upon the concept of partition of unity developed by \cite{NME:NME86}, and it was originally used to model crack propagation \cite{BB:99}. The XFEM augments the standard finite element interpolation space with additional enrichment functions to capture discontinuities in either the state variables or their spatial gradients within an element. The XFEM decomposes the elements cut by the zero level set isosurfaces into subdomains and interfaces that it uses to integrate the weak form of the governing equations. This approach avoids the need for material interpolation schemes used in density methods because each subdomain has a distinct phase. Boundary conditions on the interface are imposed weakly via penalty methods \cite{CSB:02}, stabilized Lagrange multipliers \cite{GW:08}, or via Nitsche's method \cite{BH:07}. For a general overview of the XFEM, the reader is referred to \cite{FB:10}. In the context of topology optimization of fluid flow problems, the enforcement of no-slip boundary conditions along the phase boundaries via the XFEM and a stabilized Lagrange multiplier method was adopted by \cite{KM:12} for a Navier-Stokes fluid model; a hydrodynamic Boltzmann model was employed in combination with a level set-based interface representation for generalized topology optimization of fluids by \cite{MM:14}. 

One challenge of the XFEM is that an ill-conditioned system of equations results when the ratio of the phases volumes in an intersected element is very small or very large. Such interface configurations are often unavoidable when using fixed meshes in topology optimization. In general, the ill-conditioning impedes the convergence of solvers for nonlinear problems and reduces the performance of iterative linear solvers. Several approaches have been proposed to avoid this ill-conditioning issue, such as the geometric preconditioner of \cite{LMD+:13}, the Jacobi preconditioner of \cite{SF:87}, the preconditioners of \cite{BMM+:05} and \cite{MB:11} based on a Cholesky decomposition, and face-oriented ghost-penalty methods \cite{BFH:06}. Face-oriented ghost-penalty methods have been studied in the context of fluid flow problems, where discontinuities in the spatial gradients of the velocities and the pressure across the common facets of intersected elements are penalized  \cite{BH:14,SRG+:14}. Here, we explore the characteristics and the performance of the ghost-penalty methods in the context of flow topology optimization.
	
\begin{figure}
	\centering
	\includegraphics[width=0.4\linewidth]{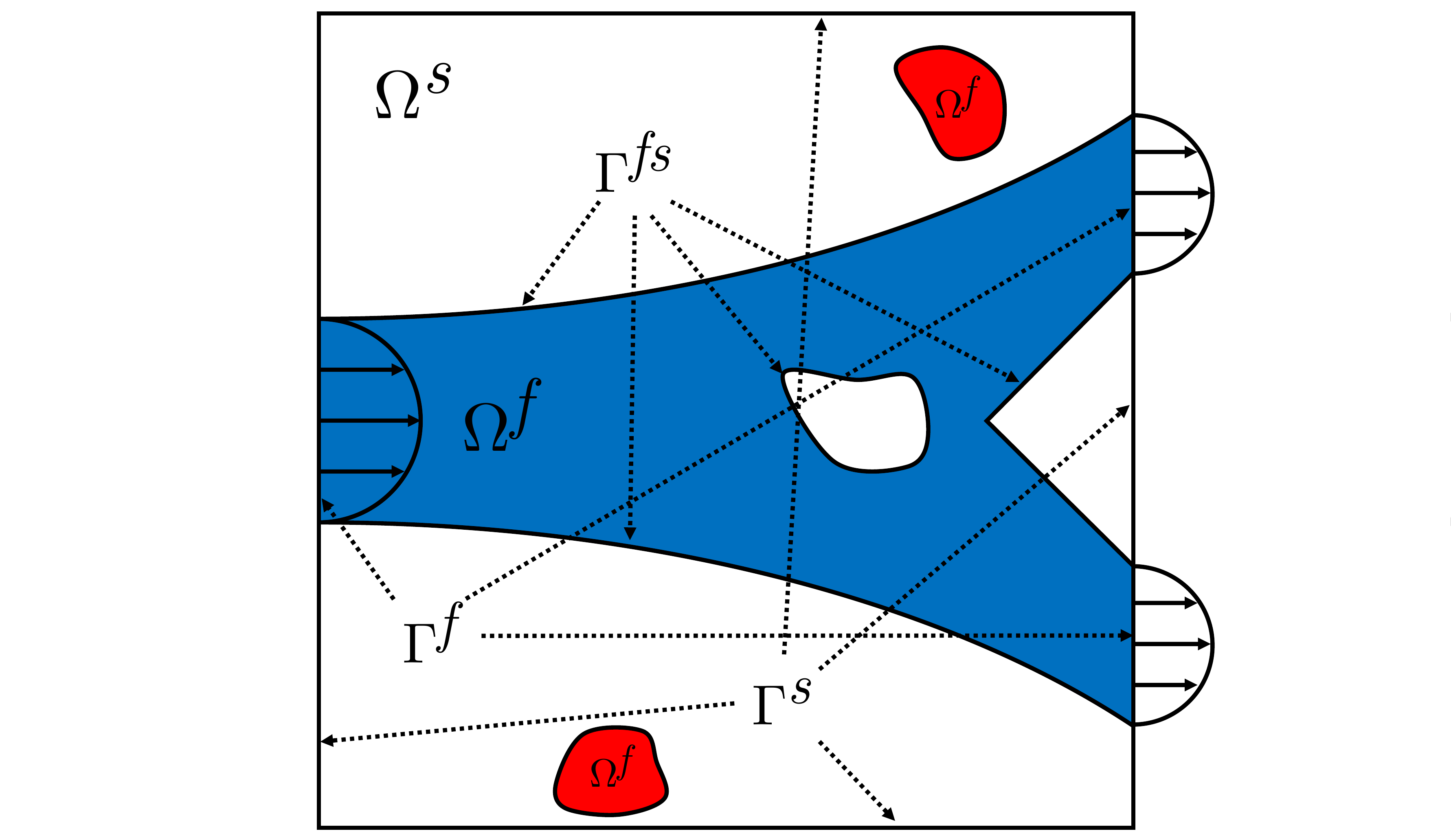}
	\caption{Geometry description of the fluid and solid domains with a LSF; blue region represents the fluid channels; isolated fluid regions surrounded by solid are denoted in red; white regions represent the solid domain.}
	\label{fig:level_set_method}
\end{figure}

A second challenge of the XFEM in the context of flow topology optimization is that isolated volumes of fluid surrounded by solid may emerge during the optimization process, as shown in Figure \ref{fig:level_set_method}. These ``puddles'' produce a singular analysis problem because the absolute value of the pressure field is not governed. A similar issue was encountered by \cite{VM:14} in the context of linear elasticity, where isolated regions of solid material could undergo rigid-body motion. The issue was avoided by placing the entire domain on a system of soft springs. A similar approach is adopted and studied here for incompressible flow, where we augment our fluid model with a penalty formulation to enforce a constraint on the average pressure. However, unlike \cite{VM:14}, where the penalty was applied to the entire domain, we model an auxiliary indicator field to detect these isolated volumes of fluid and only apply the penalty there.

A third challenge is the computational cost of performing topology optimization on 3D flows. Recent developments into large-scale studies using high-resolution density methods can be found in \cite{BP:01}, \cite{APG+:07}, and \cite{ASA:16}. As proof of concept, here we focus on laminar flows to study the characteristics of our proposed approach, avoiding thus the additional costs associated with studying large Reynolds number flows on highly refined meshes.

The combination of the LSM, the XFEM, and face-oriented ghost-penalty methods is referred to as CutFEM \cite{BCH+:14}. CutFEM can handle problems with a complex geometry description through the LSM. The boundary and interface conditions are part of the discrete formulation through the XFEM. Face-oriented ghost-penalty methods are added for stability. In this paper, we introduce and study the characteristics of the CutFEM approach as the core of a robust and broadly applicable topology optimization framework. We focus on 3D laminar flow problems modeled by the incompressible Navier-Stokes equations, as well as species transport problems. We consider two-phase problems, where one phase represents fluid, and the other represents solid. We expand the CutFEM by using a generalized Heaviside enrichment strategy for the XFEM, and adding terms to control the pressure in isolated volumes of fluid surrounded by solid. This work builds upon previous studies by \cite{KM:11} and \cite{MM:14}, but it enhances the flexibility, accuracy, and robustness of the XFEM flow analysis by weakly enforcing boundary conditions via a Nitsche formulation that accounts for viscous- and convection-dominated flow regimes, and applying face-oriented ghost-penalty methods to the velocity and pressure fields. We will show through numerical examples that these enhancements allow solving 3D steady-state and transient fluid topology optimization problems at low Reynolds number ($Re < 200$), such that relevant physical phenomena are sufficiently resolved and the boundary conditions are enforced well on the boundary geometries generated in the course of the optimization process.

The paper is structured as follows: Section \ref{sec:optimization_model} describes the formulation of our optimization problems. Section \ref{sec:geometry} details the description of the geometry using the design variables and the LSM. In Section \ref{sec:analysis}, we present the governing equations and the discretization by the XFEM. Numerical examples are studied in Section \ref{sec:numerical_examples}. The conclusions drawn from this study are presented in Section \ref{sec:conclusions}.


\section{Optimization Model}
\label{sec:optimization_model}

The optimization problems considered here are formulated with respect to an objective and one or more inequality constraints for some desired functionality. The objective and constraints are defined in terms of design criteria, such as drag, power dissipation, fluid volume, etc. These design criteria can depend explicitly on the state and optimization variables, e.g.~drag, or only on the optimization variables, e.g.~fluid volume. As the state variables may depend on time, the design criteria are either integrated over a specific period or at a given instance in time. The optimization problems of interest are formulated as follows:
\begin{equation}
	\label{eq:optimization_problem}
	\begin{split}
	\min_{\bm{s}} \quad & \mathcal{Z} = \int_{t_{1}}^{t_{2}} z \left(\bm{s}, \bm{u}\left(t,\bm{s}\right)\right) \, \mathrm{d}t\ , \\
	\mathrm{s.t.} \quad & g_{i}\left(\bm{s}, \bm{u}\left(t_{2},\bm{s}\right)\right) \le 0 \quad i = 1 \dots N_{g}\ , \\
	\bm{s} \in \bm{S} & = \left\{\mathbb{R}^{N_{\bm{s}}}\vert s_{i}^{L} \le s_{i} \le s_{i}^{U}, i = 1 \dots N_{s}\right\}\ , \\
	\bm{u}\left(t,\bm{s}\right) \in \bm{U} & = \left\{\mathbb{R}^{N_{\bm{u}}}\vert \bm{R}\left(\bm{s},\bm{u}\left(t\right)\right)= 0\right\}\ ,
	\end{split}
\end{equation}
%
where $\bm{s}$ is the vector of optimization variables, of size $N_{s}$, and $\bm{u}\left(t\right)$ is the vector of time-dependent discrete state variables, of size $N_{\bm{u}}$. The objective function $\mathcal{Z}$ is the integral of the time-dependent function $z$ over the interval $\lbrack t_{1}, t_{2} \rbrack$. The function $g_{i}$ is the $i$-th optimization constraint, and $N_{g}$ is the number of inequality constraints. In general, the constraints $g_{i}$ can be formulated analogously to the objective. However, in this study we consider state-dependent constraints that depend on the flow solution only at time $t_{2}$. The optimization variables $s_{i}$ are bounded by lower and upper limits, $s_{i}^{L}$ and $s_{i}^{U}$, respectively. The state variables satisfy the residual of the governing equations, $\bm{R}\left(\bm{s},\bm{u}\left(t\right)\right)= 0$. Note that we consider the state variables dependent on the optimization variables and, thus, solve the optimization problem in a reduced space over just $\bm{s}$.


\section{Geometry}
\label{sec:geometry}

In our CutFEM approach, the parameters of a discretized LSF are defined as explicit functions of the optimization variables. In the following subsections, we describe the concepts behind the LSM, and the parametrization of the LSF with respect to the optimization variables.


\subsection{Level Set Method}
\label{sec:level_set_method}

The LSM describes the geometry of a body immersed in a domain by the zero isosurfaces of the LSF, $\phi\left(\bm{x}\right)$. Considering a two-phase, fluid-solid problem, such as the one displayed in Figure \ref{fig:level_set_method}, the material distribution is defined as follows:
%
\begin{equation}
\label{eq:level_set}
	\begin{aligned}
		\phi(\bm{x}) &> 0    && \forall \ \bm{x} \in \mathrm{\Omega}^{s}\ , \\
		\phi(\bm{x}) &< 0    && \forall \ \bm{x} \in \mathrm{\Omega}^{f}\ , \\
		\phi(\bm{x}) &= 0    && \forall \ \bm{x} \in \mathrm{\Gamma}^{f\!s}\ ,
	\end{aligned}
\end{equation}
%
where $\mathrm{\Omega}^{s}$ is occupied by the solid phase, $\mathrm{\Omega}^{f}$ is occupied by the fluid phase, and $\mathrm{\Gamma}^{f\!s}$ is the fluid-solid interface. The external boundaries are defined by the outer surface of the design domain, $\Omega^d=\Omega^s\cup \Omega^f$, and are denoted by $\mathrm{\Gamma}^{m}$, where the superscript $m$ denotes the material phase; see Figure \ref{fig:level_set_method}.


\subsection{Parametrization of the Level Set Function}
\label{sec:parametrization_level_set}

We parametrize the LSF to describe a combination of geometric primitives and free-form shapes. In this study, geometric primitives, such as cylinders, are used to represent outflow ports; see Subsection \ref{sec:ports_geometry}. The geometry of the fluid-solid interface within the design domain is not restricted to a particular shape and topology. Free-form shapes are parametrized by a finite element approach; see Subsection \ref{sec:topology_optimization}. 

Using either parametrization scheme, the LSF is mapped onto the XFEM mesh. The level set values at the nodes are evaluated using the parametrized LSF. Within an XFEM element, the LSF is interpolated by standard trilinear finite element shape functions. Note that this approach restricts the number an element edge can be intersected by the fluid-solid interface to at most one. Therefore, the geometry resolution is limited by the size of an XFEM element, and convergence issues have been observed if sub-element-size features tend to emerge in the optimization process \cite{JM:15}. 

To prevent the formation of sub-element-size features, \cite{CM:15} proposes to add a constraint on the spatial gradient of the LSF in combination with restricting the nodal level set values to $\pm h / 2$, where $h$ is the element size. As we did not observe issues with the formation of sub-element-size features in the numerical results presented in Section \ref{sec:numerical_examples}, we did not include the feature size control of \cite{CM:15} in the formulations of our optimization problems, although our approach does allow adding additional constraints.


\subsubsection{Topology Optimization}
\label{sec:topology_optimization}

To allow for the emergence of free-form geometries within the design domain, the LSF is parametrized by local shape functions defined on a finite element mesh. In general, this mesh can differ from the XFEM mesh. To ease the computational complexity, here we use the XFEM mesh to discretize the LSF to describe a free-form geometry.

One optimization variable, $s_{i}$, for $i = 1 \dots N_{n}$, is assigned to each node of the XFEM mesh, where $N_{n}$ is the number of nodes. The LSF value of the $i$-th node, $\phi_{i}$, is defined via the following linear filter:
%
\begin{equation}
	\label{eq:smoothing_filter}
	\phi_{i} = \left(\sum^{N_{n}}_{j=1} w_{ij}\right)^{-1} \left(\sum^{N_{n}}_{j=1} w_{ij}s_{j}\right)
\qquad \textrm{with} 
\quad	w_{ij} = \max\left(0, r_{\phi} - \left\vert\bm{x}_{i} - \bm{x}_{j} \right\vert\right)\ ,
\end{equation}
%
where $r_{\phi}$ is the filter radius. The linear filter was used previously in the studies of \cite{KM:12,MM:14} to widen the zone of influence of the optimization variables, and to improve the convergence rate. Furthermore, the filter may promote (but does not guarantee) smooth shapes of the phase boundaries; however, in contrast to density or sensitivity filters used in density methods, the filter above does not guarantee control of the minimum feature size \cite{VM:14}.


\subsubsection{Ports Geometry}
\label{sec:ports_geometry}

In contrast to Hamilton-Jacobi-based updated schemes, formulating and solving the topology optimization problems via NLPs provide the flexibility to introduce additional optimization variables that do not stem from a finite element discretization of the LSF. Here we take advantage of this feature and define the geometry of cylindrical outflow ports in terms of optimization variables. The LSF of the $j$-th cylinder, $\phi_{c,j}\left(\tilde{\bm{x}}\right)$, is defined in its local coordinate system, $\tilde{\bm{x}}$, as follows:
%
\begin{equation}
	\label{eq:level_set_cylinder}
	\phi_{c,j}\left(\tilde{\bm{x}}\right) = r_{c,j} - \sqrt{\tilde{x}_1^{2} + \tilde{x}_2^{2}}
	\qquad \textrm{with} 
	\quad
	\tilde{\bm{x}}=\bm{T} \ \left({ \bm{x} - \bm{x}_c} \right)\ , 
\end{equation}
%
where $\tilde{x}_1$ and $\tilde{x}_2$ are the in-plane coordinates of the cylinder, $\tilde{x}_3$ points along the cylinder axis, and $r_{c,j}$ is the radius of the cylinder. The transformation from the global to the local coordinate system is defined by the rotation matrix $\bm{T}$ and the location of the cylinder center $\bm{x}_c$. The radius, the orientation of the local coordinate system, and the position of the cylinder center can be defined as functions of the optimization variables.

The nodal level set value, $\phi_{i}$, in the XFEM mesh is defined by the minimum value among all ports, ignoring the LSF associated with the free-form geometry, defined by \eqref{eq:smoothing_filter}. The minimum is approximated by the Kreisselmeier-Steinhauser function \cite{KS:83} to ensure differentiability of the formulation with respect to the cylinder parameters.


\section{Analysis}
\label{sec:analysis}

The main challenge in optimizing the topology of fluid flow and convective species transport problems is the modeling and numerical prediction of the flow and species fields. This section introduces the weak form of the governing equations and outlines the spatial discretization schemes. A schematic of the problem setup is shown in Figure \ref{fig:level_set_method}. The external surfaces and the fluid-solid interface, upon which Dirichlet boundary conditions are applied, are denoted by $\mathrm{\Gamma}^{m}_{\mathrm{D}}$ and $\mathrm{\Gamma}^{f\!s}_{\mathrm{D}}$, respectively. The surfaces $\mathrm{\Gamma}^{m}_{\mathrm{N}}$ and $\mathrm{\Gamma}^{f\!s}_{\mathrm{N}}$ are defined analogously for the Neumann boundary conditions.


\subsection{Governing Equations}
\label{sec:governing_equations}

In this study, we model the flow by the incompressible Navier-Stokes equations, which describe the transport of momentum and the conservation of mass. Species transport is described by coupling an advection-diffusion equation to our flow model. An indicator field is introduced to identify isolated volumes of fluid where a constraint on the average fluid pressure is enforced. The governing equations in the fluid phase are summarized subsequently. Note that the solid phase is considered ``void'' in this study; that is, we do not model any physical phenomena in the solid phase.


\subsubsection{Incompressible Navier-Stokes Equations}
\label{sec:icmp_navier_stokes}

The residual of the weak form of the incompressible Navier-Stokes equations, denoted as $r_{\bm{u},p}$, is decomposed into volumetric and surface contributions:
%
\begin{equation}
	\centering
	\label{eq:incompressible_navier_stokes_residual}
	r_{\bm{u}, p} = r_{\bm{u}, p}^{\mathrm{\Omega}} + r_{\bm{u}, p}^{\mathrm{\hat{\Omega}}} + r^{\mathrm{\Omega}}_{p,\psi} + r_{\bm{u}, p}^{\mathrm{D}} + r_{\bm{u}, p}^{f\!s} + r_{\bm{u}, p}^{\mathrm{N}} + r_{\bm{u}, p}^{\mathrm{GP}}\ ,
\end{equation}
%
where $r_{\bm{u}, p}^{\mathrm{\Omega}}$ and $r_{\bm{u}, p}^{\mathrm{\hat{\Omega}}}$ are the residuals of the volumetric contributions, non-stabilized and stabilized, respectively. The term $r^{\mathrm{\Omega}}_{p,\psi}$ is used to enforce a constraint on the average pressure in isolated volumes of fluid surrounded by solid. The terms $r_{\bm{u}, p}^{\mathrm{D}}$ and $r_{\bm{u}, p}^{f\!s}$ enforce the Dirichlet boundary conditions on the external surfaces, and at the fluid-solid interface, respectively. The residual of the Neumann conditions on the external boundaries is $r_{\bm{u}, p}^{\mathrm{N}}$. The stabilization term, $r_{\bm{u}, p}^{\mathrm{\hat{\Omega}}}$, depends on the discretization scheme, and the ghost-penalty term, $r_{\bm{u}, p}^{\mathrm{GP}}$, depends on the face-oriented ghost-penalty formulation; both are defined in Sections \ref{sec:subgrid_stabilization} and \ref{sec:ghost_penalty}, respectively.

The non-stabilized volumetric contribution is formulated as:
%
\begin{equation}
	\label{eq:weak_icmp_navier_stokes}
	r_{\bm{u}, p}^{\mathrm{\Omega}} = 
	  \int_{\mathrm{\Omega}^{f}} \Bigg(v_{i} \rho \left(\frac{\partial u_{i}}{\partial t} + u_{j} \frac{\partial u_{i}}{\partial x_{j}} \right) + \epsilon_{ij}\left(\bm{v}\right) \sigma_{ij}\left(\bm{u}, p\right)\Bigg) \, \mathrm{d}\mathrm{\Omega} 
	+ \int_{\mathrm{\Omega}^{f}} \Bigg(q \frac{\partial u_{i}}{\partial x_{i}}\Bigg) \, \mathrm{d}\mathrm{\Omega} \ .
\end{equation}
%
The first integral describes the momentum equations, with admissible test functions $v_{i}$; the second integral models the incompressibility condition, with admissible test function $q$. The fluid velocity is denoted by $u_{i}$, $p$ is the pressure, $\rho$ is the density, $\epsilon_{ij}\left(\bm{u}\right)$ is the strain rate tensor given by:
%
\begin{equation}
	\label{eq:strain_rate_tensor}
	\epsilon_{ij}\left(\bm{u}\right) = \frac{1}{2}\left(\frac{\partial u_{i}}{\partial x_{j}} + \frac{\partial u_{j}}{\partial x_{i}}\right) \ ,
\end{equation}
%
and $\sigma_{ij}\left(\bm{u}, p\right)$ is the Cauchy stress tensor for Newtonian fluids:
%
\begin{equation}
	\centering
	\label{eq:stress_tensor}
	\sigma_{ij}\left(\bm{u}, p\right) = -p \delta_{ij} + 2\mu\epsilon_{ij}\left(\bm{u}\right) \ ,
\end{equation}
%
where $\mu$ is the dynamic viscosity. 

The boundary conditions are defined as:
\begin{align}
	\label{eq:boundary_conditions}
	u_{i} &= \hat{u}^{f\!s}_{i}                 & \quad \forall \ \bm{x} & \in \mathrm{\Gamma}^{f\!s}_{\mathrm{D}}\ ,\\
	u_{i} &= \hat{u}_{i}                        & \quad \forall \ \bm{x} & \in \mathrm{\Gamma}^{f}_{\mathrm{D}}\ ,\\
	\sigma_{ij}\left(\bm{u}, p\right) n^{f}_{j} & = \hat{t}_{i}        & \quad \forall \ \bm{x} &\in \mathrm{\Gamma}^{f}_{\mathrm{N}}\ ,
\end{align}
%
where $\hat{u}^{f\!s}_{i}$ and $\hat{u}_{i}$ are the prescribed velocities, $\hat{t}_{i}$ is the traction, and $n^{f}_{j}$ is the normal on the surface pointing outwards. Dirichlet boundary conditions are enforced weakly on the fluid-solid interface, and on the external surfaces via Nitsche's method \cite{Nitsche:75}. The formulation adopted here is the one described in \cite{SRG+:14}. The surface residual of the external Dirichlet boundaries is defined as:
%
\begin{equation}
	\label{eq:weak_icmp_navier_stokes_surface}
	\begin{array}{ll}
	r_{\bm{u}, p}^{\mathrm{D}} & =
	  \int_{\mathrm{\Gamma}^{f}_{\mathrm{D}}} \bigg(v_{i} p \delta_{ij} n^{f}_{j}
	- v_{i} 2\mu \epsilon_{ij}\left(\bm{u}\right) n^{f}_{j}\bigg) \, \mathrm{d}\mathrm{\Gamma} \\
	& + \int_{\mathrm{\Gamma}^{f}_{\mathrm{D}}} \bigg(\beta_{p} q \delta_{ij} n^{f}_{j} \left(u_{i} - \hat{u}_{i}\right)
	- \beta_{\mu} 2\mu \epsilon_{ij}\left(\bm{v}\right) n^{f}_{j} \left(u_{i} - \hat{u}_{i}\right)\bigg) \, \mathrm{d}\mathrm{\Gamma} 
	+ \int_{\mathrm{\Gamma}^{f}_{\mathrm{D}}} \gamma_{\mathrm{N}, \bm{u}} v_{i} \left(u_{i} - \hat{u}_{i}\right) \, \mathrm{d}\mathrm{\Gamma} \ ,
	\end{array}
\end{equation}
%
where $\gamma_{\mathrm{N}, \bm{u}}$ is a penalty parameter. The first integral of \eqref{eq:weak_icmp_navier_stokes_surface} results from the integration by parts of the momentum equations; its two terms are denoted as the pressure and viscous \emph{standard} consistency terms, respectively. The second integral is the addition of pressure and viscous \emph{adjoint} consistency terms. The third integral introduces an additional penalty term that ensures coercivity of the viscous part of the formulation, and balances the lack of coercivity that is introduced by the viscous standard and adjoint consistency terms \cite{SRG+:14}. The terms $\beta_{p}$ and $\beta_{\mu}$ determine whether the adjoint consistency terms use a symmetric formulation ($\beta_{p} = +1$, $\beta_{\mu} = +1$), or a skew-symmetric formulation ($\beta_{p} = -1$, $\beta_{\mu} = -1$). In this work, we use the symmetric variant for the viscous adjoint consistency term because it leads to smaller errors compared to the skew-symmetric variant, as reported by \cite{Burman:12}. For the pressure adjoint consistency term, we use the skew-symmetric variation because it consistently controls the mass conservation, $u_{i} n^{f}_{i} = 0$ \cite{BMC+:10,SRG+:14}. A similar treatment is applied on the fluid-solid interface for the term $r_{\bm{u}, p}^{f\!s}$ by using $\hat{u}^{f\!s}_{i}$ instead of $\hat{u}_{i}$ in \eqref{eq:weak_icmp_navier_stokes_surface}.

The formulation of the penalty parameter $\gamma_{\mathrm{N}, \bm{u}}$ is taken from \cite{SRG+:14} and defined as:
%
\begin{equation}
	\centering
	\label{eq:nitsche_penalty_param}
	\gamma_{\mathrm{N}, \bm{u}} = \alpha_{\mathrm{N}, \bm{u}} \left(\frac{\mu}{h} + \frac{\rho {\left\Vert \bm{u} \right\Vert}_{\infty}}{6}\right)\ .
\end{equation}
%
The terms in \eqref{eq:nitsche_penalty_param} account for viscous-dominated and con\-vec\-tive-dominated flows, respectively. The term $\alpha_{\mathrm{N}, \bm{u}}$ is a constant problem-dependent penalty term, and the term ${\left\Vert \bm{u} \right\Vert}_{\infty}$ is the infinity norm evaluated at each integration point. The influence of this penalty term will be studied later in this paper.

The residual contribution from external Neumann boundary conditions is defined as:
%
\begin{equation}
	\centering
	\label{eq:weak_icmp_navier_stokes_external_neumann}
	r^{\mathrm{N}}_{\bm{u}, p} = \int_{\mathrm{\Gamma}^{f}_{\mathrm{N}}} v_{i} \hat{t}_{i} \, \mathrm{d}\mathrm{\Gamma} \ .
\end{equation}

\begin{figure}
	\centering
	\includegraphics[width=0.4\linewidth]{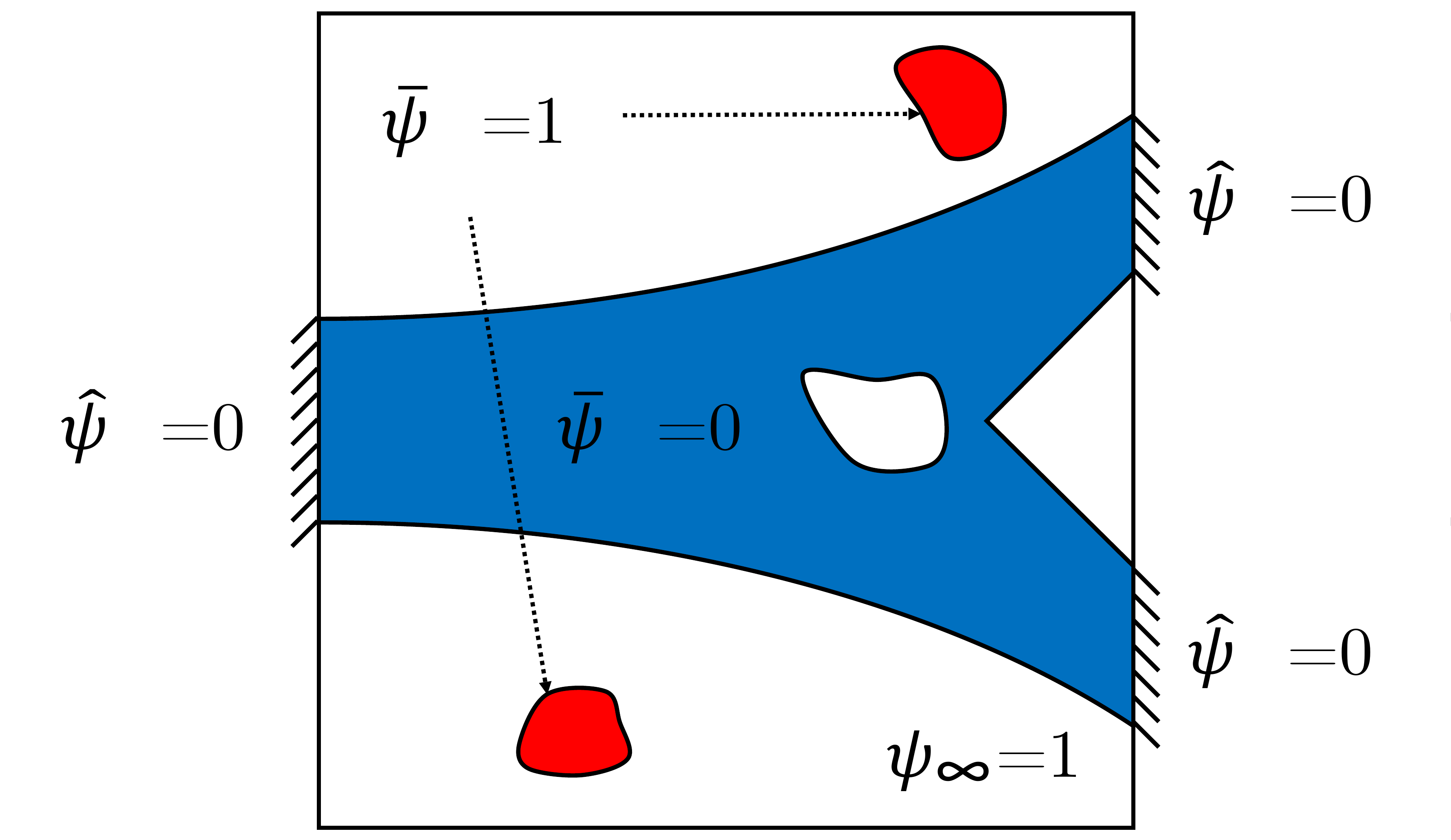}
	\caption{Modeling of the auxiliary indicator field. Isolated fluid regions surrounded by solid are denoted in red.}
	\label{fig:auxiliary_domains_fluid}
\end{figure}

An auxiliary indicator field is introduced to identify isolated volumes of fluid surrounded by the solid domain. These ``puddles'' lead to an ill-conditioned system of equations because the absolute value of the pressure is not governed. To stabilize the system, we add the following penalty formulation to enforce an average pressure:
%
\begin{equation}
	\label{eq:pressure_penalty_formulation}
	r^{\mathrm{\Omega}}_{p,\psi} = \int_{\mathrm{\Omega}^{f}} q k_{p} \bar{\psi} p \, \mathrm{d}\mathrm{\Omega} \ ,
\end{equation}
%
where $k_{p}$ is a scaling factor. The field $\bar{\psi}$ serves as a binary indicator with a value of $1$ in isolated fluid regions and zero everywhere else, as shown in Figure \ref{fig:auxiliary_domains_fluid}; the equations to model the field will be defined in Section \ref{sec:auxiliary_pressure_penalty}. The influence of the term \eqref{eq:pressure_penalty_formulation} on the conservation of mass will be studied in Section \ref{sec:pipebend}.


\subsubsection{Advection-Diffusion Equation}
\label{sec:advection_diffusion}

Species transport is modeled by an advection-dif\-fu\-sion equation. Similar to the Navier-Stokes equations \eqref{eq:weak_icmp_navier_stokes}, we denote the weak form as $r_{\bm{u}, c}$, and decompose it into volumetric and surface contributions:
%
\begin{equation}
	\centering
	\label{eq:weak_species_transport}
	r_{\bm{u}, c} = r^{\mathrm{\Omega}}_{\bm{u}, c} + r^{\mathrm{\hat{\Omega}}}_{\bm{u}, c} + r^{\mathrm{D}}_{\bm{u}, c} + r^{\mathrm{N}}_{\bm{u}, c} + r^{\mathrm{N},f\!s}_{\bm{u}, c} + r^{\mathrm{GP}}_{\bm{u}, c} \ ,
\end{equation}
%
where $r^{\mathrm{\Omega}}_{\bm{u}, c}$ and $r^{\mathrm{\hat{\Omega}}}_{\bm{u}, c}$ are the residuals of the volumetric contributions, non-stabilized and stabilized, respectively; the term $r^{\mathrm{D}}_{\bm{u}, c}$ represents the Dirichlet boundary conditions at the external boundaries; the terms $r^{\mathrm{N}}_{\bm{u}, c}$ and $r^{\mathrm{N},f\!s}_{\bm{u}, c}$ describe the Neumann boundary conditions at the external boundaries, and at the fluid-solid interface, respectively; and the term $r^{\mathrm{GP}}_{\bm{u}, c}$ models the ghost-penalty formulation. Similar to \eqref{eq:incompressible_navier_stokes_residual}, the residual $r^{\mathrm{\hat{\Omega}}}_{\bm{u}, c}$ is defined in Section \ref{sec:subgrid_stabilization}, and the residual $r^{\mathrm{GP}}_{\bm{u}, c}$ is defined in Section \ref{sec:ghost_penalty}.

The non-stabilized volumetric residual contribution, $r^{\mathrm{\Omega}}_{\bm{u}, c}$, is defined as:
%
\begin{equation}
	\label{eq:weak_species_transport_volume}
	r^{\mathrm{\Omega}}_{\bm{u}, c} = \int_{\mathrm{\Omega}^{f}} \Bigg(w \bigg(\frac{\partial c}{\partial t} + u_{i}\frac{\partial c}{\partial x_{i}}\bigg)
	+ \frac{\partial w}{\partial x_{i}} \ \bigg(J_{i}\left(c\right)\bigg)
	- w \hat{q}_{\mathrm{\Omega}}\Bigg) \, \mathrm{d\Omega} \ ,
\end{equation}
%
where $c$ is the species concentration, $w$ is an admissible test function, $u_{i}$ is the vector of fluid velocities, $\hat{q}_{\mathrm{\Omega}}$ is the volumetric source, and $J_{i}\left(c\right)$ is the diffusive flux defined as:
%
\begin{equation}
	\label{eq:diffusive_heat_flux}
	J_{i}\left(c\right) = k \ \delta_{ij} \ \frac{\partial c}{\partial x_{j}}  \ ,
\end{equation}
%
where $k$ is the isotropic diffusivity.

Dirichlet and Neumann boundary conditions are imposed on the external boundaries, $\mathrm{\Gamma}^{f}$, and Neumann conditions on the fluid-solid interface, $\mathrm{\Gamma}^{f\!s}$, as:
%
\begin{align}
	\label{eq:advection_diffusion_bcs}
	c                                 &= \hat{c}                          &\quad \forall \ \bm{x}  &\in \mathrm{\Gamma}^{f}_{\mathrm{D}}\ ,\\
	J_{i}\left(c\right) n^{f}_{i}     &= \hat{q}_{\mathrm{\Gamma}}        &\quad \forall \ \bm{x}   &\in \mathrm{\Gamma}^{f}_{\mathrm{N}}\ ,\\
	J_{i}\left(c\right) n^{f\!s}_{i}  &= \hat{q}^{f\!s}_{\mathrm{\Gamma}} &\quad \forall \ \bm{x}   &\in \mathrm{\Gamma}^{f\!s}_{\mathrm{N}}\ ,
\end{align}
%
where $n^{f\!s}_{i}$ is the normal vector on the fluid-solid interface pointing towards the solid phase, $\hat{c}$ is a prescribed concentration value, and $\hat{q}_{\mathrm{\Gamma}}$ and $\hat{q}^{f\!s}_{\mathrm{\Gamma}}$ are prescribed species flux values.

The weak enforcement of Dirichlet boundary conditions is modeled using Nitsche's method \cite{Nitsche:75}. The residual contributions of the Dirichlet conditions at the external surfaces are defined as:
%
\begin{equation}
		\label{eq:weak_species_transport_external_dirichlet}
	r^{\mathrm{D}}_{\bm{u}, c} = \int_{\mathrm{\Gamma}^{f}_{\mathrm{D}}} \Bigg(-w J_{i} \left(c\right) n^{f}_{i}
	+ J_{i} \left(w\right) n^{f}_{i} c	+ \alpha_{\mathrm{N},c} h^{-1} w \left(c - \hat{c}\right)\Bigg) \, \mathrm{d}\mathrm{\Gamma}  \ ,
\end{equation}
%
where the penalty parameter $\alpha_{\mathrm{N},c}$ is a problem-dependent constant. The first, second, and third terms correspond to the standard consistency term, the adjoint consistency term, and the penalty term of the Nitsche formulation, respectively.

The Neumann contributions at the fluid-solid interface and  the external boundaries are defined as:
%
\begin{equation}
	\centering
	\label{eq:weak_species_transport_internal_neumann}
	r^{\mathrm{N}, f\!s}_{\bm{u}, c} = \int_{\mathrm{\Gamma}^{f\!s}_{\mathrm{N}}} w \hat{q}^{f\!s}_{\mathrm{\Gamma}} \, \mathrm{d}\mathrm{\Gamma}
	\qquad
	\textrm{and}
	\qquad
	r^{\mathrm{N}}_{\bm{u}, c} = \int_{\mathrm{\Gamma}^{f}_{\mathrm{N}}} w \hat{q}^{f}_{\mathrm{\Gamma}} \, \mathrm{d}\mathrm{\Gamma}\ .
\end{equation}


\subsection{Auxiliary Indicator Field}
\label{sec:auxiliary_pressure_penalty}

The auxiliary field is modeled as a linear diffusion problem in the fluid domain. The residual of the weak form is defined as:
%
\begin{equation}
	\label{eq:auxiliary_species_field_residual}
	r_{\psi} = r^{\mathrm{\Omega}}_{\psi} + r^{\mathrm{D}}_{\psi} + r^{\mathrm{GP}}_{\psi}\ ,
\end{equation}
%
where the residual of the Dirichlet boundary conditions, $r^{\mathrm{D}}_{\psi}$, is formulated in the same way as in  \eqref{eq:weak_species_transport_external_dirichlet}. Dirichlet boundary conditions are imposed on all inlets and outlets, by setting $\hat{\psi}$ to $0$. Adiabatic boundary conditions are imposed on the fluid-solid interface. 

The volumetric residual contribution, $r^{\mathrm{\Omega}}_{\psi}$ is defined as:
%
\begin{equation}
	\label{eq:weak_auxiliary_species_field_volume}
	r^{\mathrm{\Omega}}_{\psi} = \int_{\mathrm{\Omega}^{f}} \Bigg(\frac{\partial \xi}{\partial x_{i}} J_{i}\left(\psi\right) - \xi h_{\psi}\left(\psi - \psi_{\infty}\right)\Bigg) \, \mathrm{d}\mathrm{\Omega}\ ,
\end{equation}
%
where $\xi$ is an admissible test function, $h_{\psi}$ is the convection coefficient, and $\psi_{\infty}$ is the reference indicator value. The diffusion coefficient in \eqref{eq:diffusive_heat_flux} is set to $1$. The parameters $h_{\psi}$ and $\psi_{\infty}$ are set to $0.01$ and $1$, respectively, so that fluid channels that are connected to the inlets and outlets will have an indicator field value close to $0$, while the fluid ``puddles'' will have a value close to $1$.

A smooth-Heaviside projection scheme is applied to the indicator field to project the values of the solution either to $0$ or to $1$, and is defined as:
%
\begin{equation}
	\label{eq:pressure_penalty_projection}
	\bar{\psi} = \frac{1}{2} + \frac{1}{2}\tanh\left(k_{w}\left(\psi - k_{t} \psi_{\infty}\right)\right)\ ,
\end{equation}
%
where the sharpness of the projection increases with the parameter $k_{w}$, and $k_{t} \psi_{\infty}$ is the threshold. We adopt the values of $k_{w} = 1000$ and $k_{t} = 0.99$ to effectively turn the $\bar{\psi}$ term into a binary switch, where a value of $0$ corresponds to a point connected to the inlet and/or outlet ports, and a value of $1$ corresponds to an isolated volume of fluid.


\subsection{Spatial Discretization}
\label{sec:spatial_discretization}

The governing equations in the fluid phase are discretized in space by the XFEM. This study adopts a generalized enrichment strategy based on the Heaviside-step enrichment of \cite{HH:04}, which interpolates consistently the solution fields in the presence of small features, and does not suffer from artificial coupling of disconnected phases. This particular approach has been used by \cite{MM:13,VM:14}, \cite{KM:12}, \cite{LMD+:13}, and \cite{MM:14}, who considered linear elasticity, incompressible Navier-Stokes, linear diffusion, and advection-diffusion problems, respectively.

Here we discretize the fluid state variables, the species concentration, and the indicator field in the fluid domain by the XFEM. Using the symbol $\omega$ to represent any of these state variables, the approximation of $\omega$ within an element, $\tilde{\omega}$, can be written as follows:
%
\begin{equation}
	\centering
	\label{eq:heaviside_enrichment}
	\omega\left(\bm{x}\right) \approx \tilde{\omega}\left(\bm{x}\right) = \sum^{N_{l}}_{l=1} \
	\sum^{N_{n}}_{i=1} \mathcal{N}_{i}(\bm{x}) \ \delta^{i}_{lk} \ \omega^{i}_{l}  \qquad \forall \ \bm{x} \in \Omega^{f}\ ,
\end{equation}
%
where $l$ is the enrichment level, $N_{l}$ is the maximum number of enrichment levels, $\mathcal{N}_{i}(\bm{x})$ are the nodal basis functions, $\omega^{i}_{l}$  is the de\-grees-of-free\-dom of enrichment level $l$ at node $i$ in the fluid phase, respectively. The Kronecker delta, $\delta^{i}_{l k}$, selects the de\-gree-of-free\-dom that is active at node $i$. At any given point, only one de\-gree-of-free\-dom per node is used to interpolate the solution, ensuring that the partition of unity is satisfied.

Multiple enrichment levels, i.e. sets of shape functions, may be necessary to interpolate the state variables in multiple, physically disconnected regions of the same phase, c.f. \cite{TAY:03}, \cite{TLL:11}, and \cite{MM:13}. When interpolating the level set field by element-wise linear functions in a structured grid, a maximum of $14$ enrichment levels is needed in 3D problems \cite{VM:14}. In order to accurately integrate the weak form of the governing equations by Gauss quadrature, intersected elements are decomposed into tetrahedrons using Delaunay triangulation. The reader is referred to \cite{MM:13} for more details on the particular XFEM implementation used in this study.

The Heaviside-step enrichment formulation \eqref{eq:heaviside_enrichment} is ill-defined for cases in which the fluid-solid interface lies exactly on a node, i.e. the level set value $\phi_{i}$ at node $i$ equals $0$. To avoid this issue, we adopt the level set perturbation approach outlined in \cite{CHM:12} and \cite{LMD+:13}. If the magnitude of the level set value at a node is smaller than some critical value, $\phi_{c}$, the level set value is modified to a shift value, $\phi_{s}$. This perturbation results in the fluid-solid interface moving away from the node, solving the singularity issue. In this study, we adopt the values of $\phi_{c} = \phi_{s} = 10^{-6} \times h$. Numerical studies have shown that the influence of this perturbation is negligible for the problems considered here \cite{CM:15b}.


\subsection{Subgrid Stabilization}
\label{sec:subgrid_stabilization}

The convective terms in the incompressible Navier-Stokes and advection-diffusion equations may cause spurious node-to-node velocity oscillations. Furthermore, the equal-order approximations used for $u_{i}$ and $p$ may give rise to spurious pressure oscillations. To prevent these numerical instabilities, we augment the incompressible Navier-Stokes equations with the Streamline Upwind Petrov-Galerkin (SUPG) and the Pressure Stabilized Petrov-Galerkin (PSPG) stabilization formulations introduced by \cite{TMR+:92}. The stabilized volumetric residual contribution of \eqref{eq:incompressible_navier_stokes_residual}, $r_{\bm{u}, p}^{\hat{\mathrm{\Omega}}}$, is defined as:
%
\begin{multline}
	\label{eq:weak_icmp_navier_stokes_subgrid}
	r_{\bm{u}, p}^{\hat{\mathrm{\Omega}}} = \sum_{\mathrm{\Omega}_{e} \in \mathrm{\Omega}} \int_{\mathrm{\Omega}_{e} \cap \mathrm{\Omega}} 
	\Bigg( \tau_{\mathrm{SUPG}, \bm{u}} \left(u_{j} \frac{\partial v_{i}}{\partial x_{j}}\right) + \tau_{\mathrm{PSPG}} \left(\frac{1}{\rho}\frac{\partial q}{\partial x_{i}}\right) \Bigg) \ \cdot \\
	\Bigg(\rho\left(\frac{\partial u_{i}}{\partial t} + u_{j} \frac{\partial u_{i}}{\partial x_{j}}\right) + \frac{\partial p}{\partial x_{j}}\delta_{ij} - 2 \mu \frac{\partial}{\partial x_{j}}\bigg(\epsilon_{ij}\left(\bm{u}\right)\bigg) \Bigg) \ \mathrm{d}\mathrm{\Omega} \ ,
\end{multline}
%
where $\mathrm{\Omega}_{e}$ denotes the set of all elements in the domain $\mathrm{\Omega}$, and the stabilization terms $\tau_{\mathrm{SUPG}, \bm{u}}$ and $\tau_{\mathrm{PSPG}}$ are defined in \cite{TMR+:92}.

The stabilized volumetric residual contribution of the advection-diffusion model \eqref{eq:weak_species_transport}, $r_{\bm{u}, c}^{\hat{\mathrm{\Omega}}}$, uses the SUPG method, and is defined as:
\begin{equation}
	\label{eq:weak_species_transport_subgrid}
	r_{\bm{u}, c}^{\hat{\mathrm{\Omega}}} = \sum_{\mathrm{\Omega}_{e} \in \mathrm{\Omega}} \int_{\mathrm{\Omega}_{e} \cap \mathrm{\Omega}} 
	\tau_{\mathrm{SUPG}, c} \ \left( \frac{\partial w}{\partial x_{i}} \right) \ \Bigg( \frac{\partial c}{\partial t} + u_{i} \ \frac{\partial c}{\partial x_{i}} - \frac{\partial}{\partial x_{i}} \bigg(J_{i}\left(c\right)\bigg) \Bigg) \, \mathrm{d}\mathrm{\Omega} \ ,
\end{equation}
%
\noindent where the stabilization terms $\tau_{\mathrm{SUPG}, c}$ is defined in \cite{FFH:92}.


\subsection{Face-oriented Ghost-penalty Methods}
\label{sec:ghost_penalty}

As the geometry of the design evolves during the optimization process, the fluid-solid interface geometry may lead to intersection configurations where certain degrees-of-freedom interpolate in very small subdomains. This produces an ill-conditioning of the system, which manifests itself through an increase in the condition number of the linearized system, and may slow down or prevent the convergence of the nonlinear problem. To guarantee stability, as well as to improve the conditioning of the system, face-oriented ghost-penalty stabilization terms are used in the vicinity of the fluid-solid interface, c.f. \cite{BFH:06}. The ghost-penalty terms for the residual contribution of the incompressible Navier-Stokes equations are defined as:
%
\begin{equation}
	\label{eq:ghost_penalty_incompressible_residual}
	r_{\bm{u}, p}^{\mathrm{GP}} = r_{\bm{u}, p}^{\mathrm{GP}, \mu} + r_{\bm{u}, p}^{\mathrm{GP}, p} + r_{\bm{u}, p}^{\mathrm{GP}, \bm{u}}  \ ,
\end{equation}
%
where $r_{\bm{u}, p}^{\mathrm{GP}, \mu}$, $r_{\bm{u}, p}^{\mathrm{GP}, p}$, and $r_{\bm{u}, p}^{\mathrm{GP}, \bm{u}}$ are the viscous, pressure, and convective ghost-penalty formulations, respectively.

\begin{figure}
	\centering
		\includegraphics[width=0.25\linewidth]{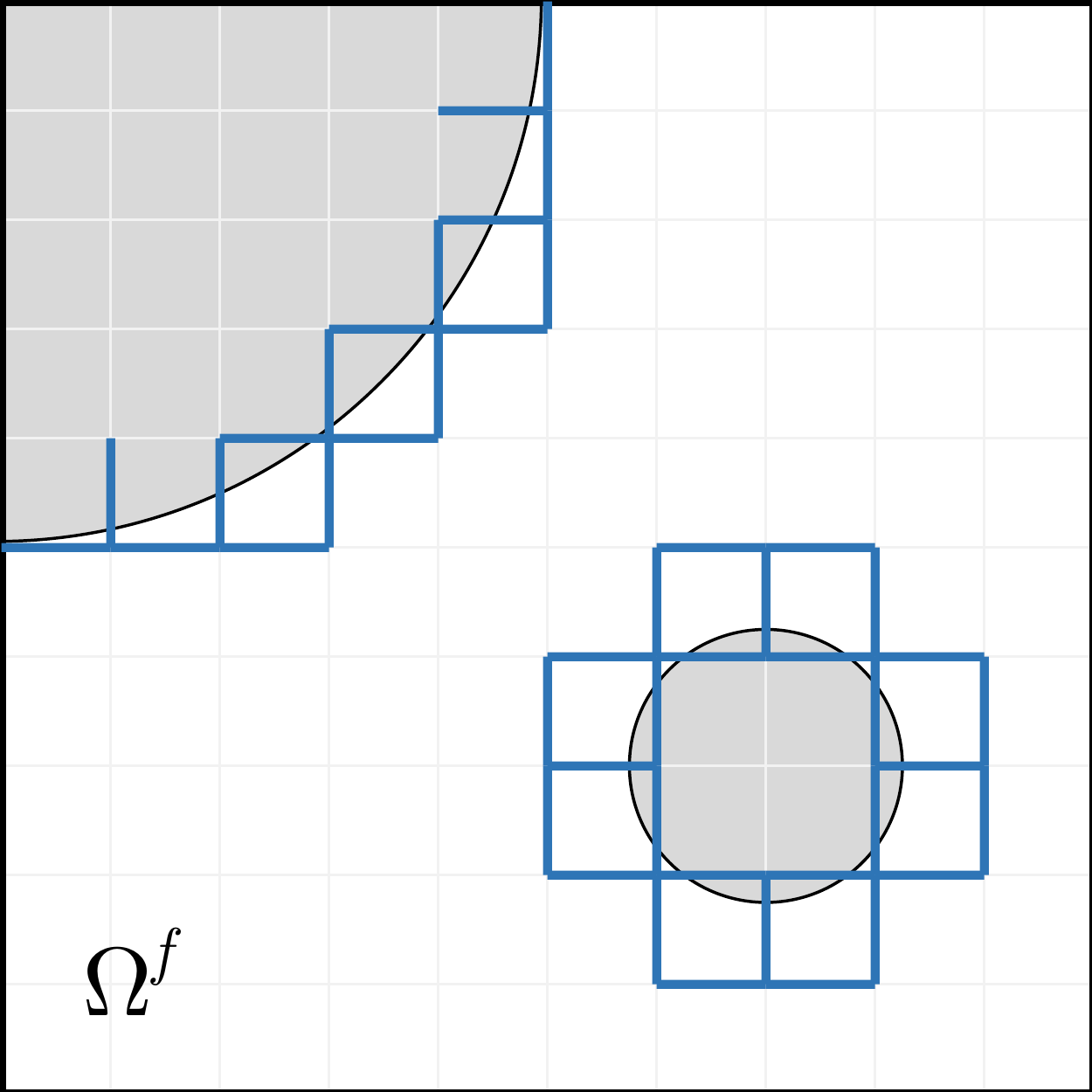}
	\caption{Integration domains for the face-oriented ghost-penalty method; facets $\mathrm{\Xi}$ in domain $\mathrm{\Omega}^{f}$.}
	\label{fig:geometry_ghost_penalty}
\end{figure}

We adopt the viscous face-oriented ghost-penalty formulation as proposed by \cite{BH:14}:
%
\begin{equation}
	\label{eq:viscous_ghost_penalty}
	r_{\bm{u}, p}^{\mathrm{GP}, \mu} = \sum_{F \in \mathrm{\Xi}} \int_{F} \Bigg(\gamma_{\mathrm{GP}, \mu} \left\llbracket \frac{\partial v_{i}}{\partial x_{j}} \right\rrbracket n^{f}_{j} \left\llbracket \frac{\partial u_{i}}{\partial x_{k}} \right\rrbracket  n^{f}_{k}\Bigg) \, \mathrm{d} \mathrm{\Gamma} \ ,
\end{equation}
%
where $\gamma_{\mathrm{GP}, \mu}$ is a penalty parameter defined as:
%
\begin{equation}
	\label{eq:viscous_ghost_penalty_gamma}
	\gamma_{\mathrm{GP}, \mu} = \alpha_{\mathrm{GP}, \mu} \mu h \ ,
	\end{equation}
%
and $\alpha_{\mathrm{GP}, \mu}$ is a constant scaling factor. The jump operator is defined as:
%
\begin{equation}
	\label{eq:ghost_penalty_jump_op}
	\llbracket \zeta \rrbracket = {\zeta\vert}_{\mathrm{\Omega}^{1}_{e}} - {\zeta\vert}_{\mathrm{\Omega}^{2}_{e}}\ ,
	\end{equation}
%
and is evaluated at the facet between two adjacent elements, $\mathrm{\Omega}^{1}_{e}$ and $\mathrm{\Omega}^{2}_{e}$. This formulation overcomes the issue of having a small ratio of volumes on elements bisected by the interface because the domain of integration is the entire edge, regardless of the intersection configuration. As illustrated in Figure \ref{fig:geometry_ghost_penalty}, the set $\mathrm{\Xi}$ belonging to the domain $\mathrm{\Omega}$ contains all facets $F$ in the immediate vicinity of the fluid-solid interface, for which at least one of the two adjacent elements is cut by the interface.

To control pressure instabilities due to a violated inf-sup condition for equal-order approximations used for $u_{i}$ and $p$ \cite{SRG+:14}, a pressure ghost-penalty stabilization term is applied:
%
\begin{equation}
	\label{eq:pressure_ghost_penalty}
	r_{\bm{u}, p}^{\mathrm{GP}, p} = \sum_{F \in \mathrm{\Xi}} \int_{F} \Bigg(\gamma_{\mathrm{GP}, p} \left\llbracket \frac{\partial q}{\partial x_{j}} \right\rrbracket n^{f}_{j} \left\llbracket \frac{\partial p}{\partial x_{k}} \right\rrbracket  n^{f}_{k}\Bigg) \, \mathrm{d} \mathrm{\Gamma} \ ,
\end{equation}
%
where $\gamma_{\mathrm{GP}, p}$ is a penalty parameter defined as:
%
\begin{equation}
	\label{eq:pressure_ghost_penalty_gamma}
	\gamma_{\mathrm{GP}, p} = \alpha_{\mathrm{GP}, p} {\left(\frac{\mu}{h} + \frac{\rho{\left\Vert\bm{u}\right\Vert}_{\infty}}{6}\right)}^{-1} h^{2}\ ,
	\end{equation}
%
and accounts for the viscous and convective flow regimes, c.f. \cite{BFH:06}. The term $\alpha_{\mathrm{GP}, p}$ is a constant scaling parameter.

For high Reynolds number flows, a convective ghost-penalty formulation was proposed by \cite{SW:14} to have sufficient control over the convective derivative, $u_{i} \nabla u_{i}$, of the incompressible Navier-Stokes equations. This formulation is defined as:
%
\begin{equation}
	\label{eq:convective_ghost_penalty}
	r_{\bm{u}, p}^{\mathrm{GP}, \bm{u}} = \sum_{F \in \mathrm{\Xi}} \int_{F} \Bigg(\gamma_{\mathrm{GP}, \bm{u}} \left\llbracket \frac{\partial v_{i}}{\partial x_{j}} \right\rrbracket n^{f}_{j} \left\llbracket \frac{\partial u_{i}}{\partial x_{k}} \right\rrbracket  n^{f}_{k}\Bigg) \, \mathrm{d} \mathrm{\Gamma}\ ,
\end{equation}
%
where the parameter $\gamma_{\mathrm{GP}, \bm{u}}$ is a penalty factor defined as:
%
\begin{equation}
	\label{eq:convective_ghost_penalty_gamma}
	\gamma_{\mathrm{GP}, \bm{u}} = \alpha_{\mathrm{GP}, \bm{u}} \rho {\left\Vert u_{i} n^{f}_{i}\right\Vert} h^{2}\ ,
\end{equation}
%
and $\alpha_{\mathrm{GP}, \bm{u}}$ is a constant scaling parameter. Additional ghost-penalty measures have been proposed in the literature, for example, to control instabilities arising from the incompressibility constraint. However, these additional formulations are not considered in this paper because previous studies have not revealed any further improvement for the laminar flow situations analyzed here, as stated by \cite{SRG+:14}.

To stabilize the concentration field \eqref{eq:weak_species_transport_volume}, we use the formulation from \cite{BH:12}:
%
\begin{equation}
	\label{eq:species_ghost_penalty}
	r_{\bm{u}, c}^{\mathrm{GP}} = \sum_{F \in \mathrm{\Xi}} \int_{F} \Bigg(\gamma_{\mathrm{GP}, c} \left\llbracket \frac{\partial w}{\partial x_{i}} \right\rrbracket n^{f}_{i} \left\llbracket J_{j}\left(c\right) \right\rrbracket n^{f}_{j}\Bigg) \, \mathrm{d} \mathrm{\Gamma}\ ,
\end{equation}
%
where $\gamma_{\mathrm{GP}, c}$ is a penalty parameter defined as:
%
\begin{equation}
	\label{eq:species_ghost_penalty_gamma}
	\gamma_{\mathrm{GP}, c} =  \alpha_{\mathrm{GP}, c} \ h\ ,
\end{equation}
%
and $\alpha_{\mathrm{GP}, c}$ is a scaling constant.

The ghost-penalty formulation for the auxiliary indicator field is identical to \eqref{eq:species_ghost_penalty}, except that it operates on a different admissible test function, $\xi$, a different set of degrees-of-freedom, $\psi$, and a different scaling factor, $\alpha_{\mathrm{GP}, \psi}$.

The values for $\alpha_{\mathrm{GP}, \mu}$, $\alpha_{\mathrm{GP}, p}$, $\alpha_{\mathrm{GP}, \bm{u}}$, and $\alpha_{\mathrm{GP}, c}$ are set on a per-problem basis. The value of $\alpha_{\mathrm{GP}, \psi}$ is set to $0.05$ for all numerical examples, in accordance to the parameter used by \cite{BH:12} for a linear diffusion field.


\section{Numerical Examples}
\label{sec:numerical_examples}

In the following, we study the characteristics of the proposed CutFEM topology optimization framework for steady-state and transient laminar flow and species transport problems in 3D. We first compare the accuracy and convergence of the analysis of the proposed CutFEM framework against a body-fitted problem taken from the literature. Then, we study the effects on the mass conservation caused by the penalty formulation introduced in \eqref{eq:pressure_penalty_formulation}. Finally, we apply topology optimization to several problems in order to study the characteristics of the framework with respect to different physical phenomena. Unless otherwise stated, geometric and material parameters are given in non-dimensional and self-consistent units.

In this study, time integration is performed by a two-step backward differentiation scheme. A sufficiently large time step is chosen to simulate steady-state conditions. For all time steps, equilibrium is satisfied by solving the associated nonlinear system of equations via Newton's method. Linear problems are solved by the Generalized Minimal RESidual (GMRES) iterative method \cite{SS:86}, with an Incomplete LU factorization with dual Threshold (ILUT) preconditioner \cite{Saad:94}. For all examples studied in the following, the nonlinear and linear problems are considered converged if the relative residuals are less then $10^{-6}$.

The optimization problems (\ref{eq:optimization_problem}) are solved by a gradient-based algorithm in the reduced space. The gradients of the objective and constraint functions with respect to the optimization variables, $s_{i}$, are computed via the adjoint method. In this work, we adopt the discrete adjoint formulation for nonlinear fluid and coupled systems of \cite{KM:11} and \cite{GMD:12}. The problems are solved via the Globally Convergent Method of Moving Asymptotes (GCMMA) of \cite{Svanberg:95}. The GCMMA parameters are given in Table \ref{tab:gcmma_params}. The optimization problem is considered converged if the change of the objective function relative to the previous objective value is less than $10^{-6}$, and all constraints are satisfied.
\begin{table}
	\centering
	\begin{tabular*}{0.625\linewidth}{@{\extracolsep{\fill}} l l}
	\hline
	                             & Value \\
    \hline
	Relative step size           & $0.04$ \\
	Minimum asymptote adaptivity & $0.5$ \\
	Initial asymptote adaptivity & $0.7$ \\
	Maximum adaptivity           & $1.43$ \\
	Constraint penalty           & $100$ \\
    \hline
	\end{tabular*}
	\caption{GCMMA parameters for the topology optimization problems.}
	\label{tab:gcmma_params}
\end{table}


\subsection{Optimization criteria}
\label{sec:optimization_criteria}

Here we summarize the design criteria used to formulate the objectives and the constraints in the optimization problems studied in this section. These include:

{\it Drag Coefficient:} The drag coefficient is used to qualify the forces of laminar flow on a surface, and is defined as:
%
\begin{equation}
	\label{eq:drag_coeff_criteria}
	\begin{split}
	c_{D} = -2e_{i}{\Bigg(\rho {\vert\bm{u}_{c}\vert}^{2}L_{c}\Bigg)}^{-1} \int_{\mathrm{\Gamma}^{f}} \sigma_{ij}\left(\bm{u}, p\right) n^{f}_{j} \, \mathrm{d}\mathrm{\Gamma} \ ,
	\end{split}
\end{equation}
%
where $e_{i}$ is a unit vector pointing in the direction of the flow velocity, $\bm{u}_{c}$ is the characteristic velocity, and $L_{c}$ is the characteristic length.

{\it Mass Flow Rate:} The mass flow rate criterion computes the mass of the fluid that passes through a surface per unit of time, and is defined as:
%
\begin{equation}
	\centering
	\label{eq:mass_flow_rate_criteria}
	\dot{m} = \int_{\mathrm{\Gamma}^{f}} \Bigg(\rho u_{i} \, n^{f}_{i} \Bigg)\,\mathrm{d}\mathrm{\Gamma} \ .
\end{equation}

{\it Pressure difference:} The total pressure criterion measures the sum of the static and dynamic pressures over a surface:
%
\begin{equation}
	\label{eq:pressure_drop_criteria}
	\mathcal{T} = \int_{\mathrm{\Gamma}^{f}} \Bigg(p + \frac{\rho \vert \bm{u} \vert^2}{2}\Bigg) \, \mathrm{d}\mathrm{\Gamma}  \ .
\end{equation}

{\it Volume:} The volumes of the fluid and solid domains are computed as:
%
\begin{equation}
	\label{eq:volume_criteria}
	\mathcal{V}^{f} = \int_{\mathrm{\Omega}^{f}} \, \mathrm{d}\mathrm{\Omega} \qquad \textrm{and} \qquad \mathcal{V}^{s} = \int_{\mathrm{\Omega}^{s}} \, \mathrm{d}\mathrm{\Omega}\ .
\end{equation}

{\it Surface Area:} The surface area criterion is computed at the fluid-solid interface, and is defined as:
%
\begin{equation}
	\centering
	\label{eq:surface_criteria}
	\mathcal{S} = \int_{\mathrm{\Gamma}^{f\!s}} \, \mathrm{d}\mathrm{\Gamma} \ .
\end{equation}
%
Reducing the value of this measure, either with a penalty in the objective functional or with a constraint, discourages the emergence of small geometric features and oscillatory shapes in the optimization problem. While a surface area penalty does not allow explicitly controlling the local shape and the feature size \cite{VM:14}, it has been reported effective in regularizing flow optimization problems \cite{DML+:13,MM:14}.

{\it Target Scalar Value:} To measure the maximum difference between a current species concentration and a target species concentration, $c_{\mathrm{ref}}$, over the fluid phase, we use the Kreisselmeier-Steinhauser function \cite{KS:83}:
%
\begin{equation}
	\centering
	\label{eq:target_scalar_value_criteria}
	\mathcal{K} = \frac{1}{\beta_{KS}} \ln \int_{\mathrm{\Gamma}^{f}} \Bigg(e^{\beta_{KS}\left(c - c_{\mathrm{ref}} \right)^2}\Bigg) \,\mathrm{d}\mathrm{\Gamma} \ .
\end{equation}
%
A large value of $\beta_{KS}$ increases the accuracy of approximating the maximum, but may result in large design sensitivities that affect the convergence of the optimization problem.


\subsection{Verification of the CutFEM analysis}
\label{sec:dfg_cylinder}

\begin{figure*}
	\centering
	\includegraphics[width=0.67\linewidth]{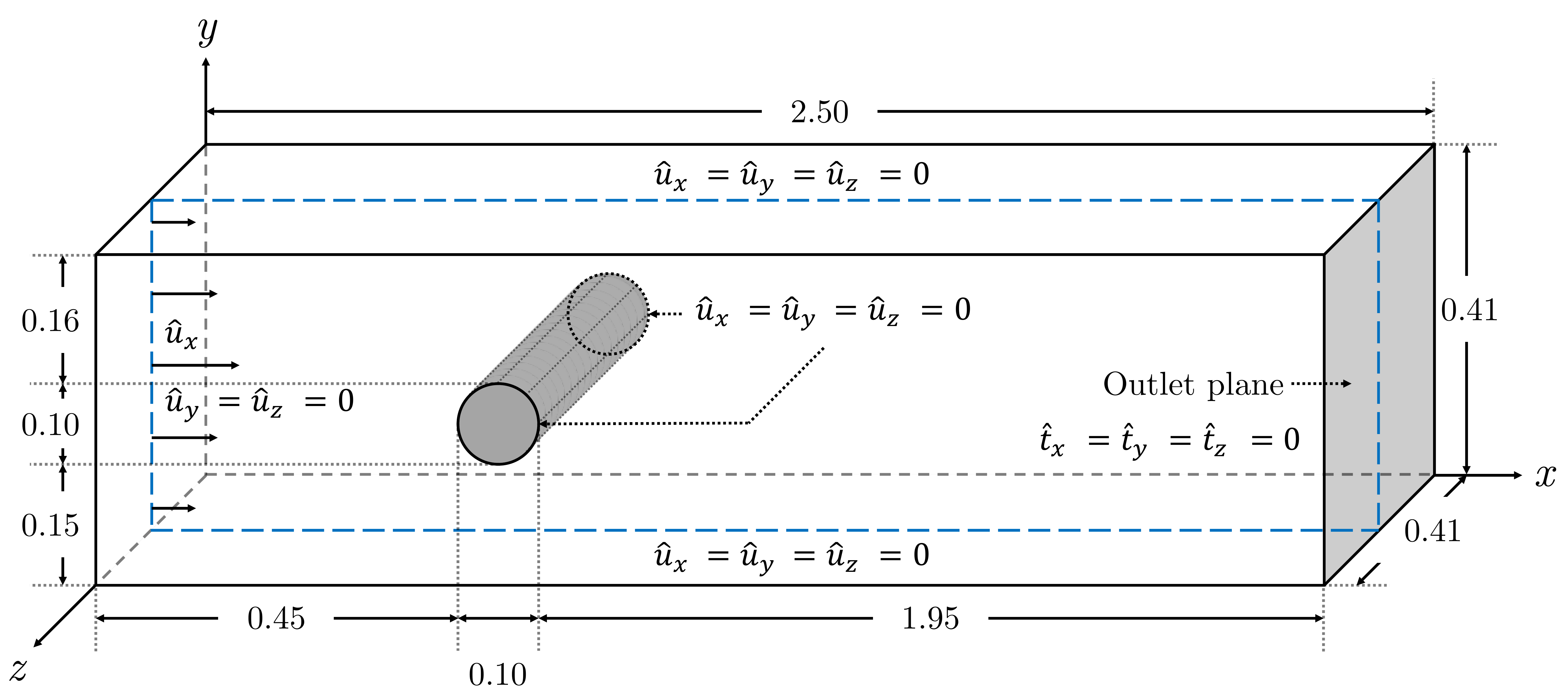}
	\caption{Problem setup for the flow around a cylinder example. Blue dashes lines denote the symmetry plane.}
	\label{fig:setup_dfg_cylinder}
\end{figure*}

In this first example, we seek to verify the analysis capabilities of our CutFEM optimization framework. The validity of our approach is measured in terms of the accuracy of the flow solutions with respect to solutions from the literature, and in terms of the convergence of the flow solutions with respect to different levels of mesh refinement and different values for the Nitsche penalty parameter, $\alpha_{\mathrm{N},\bm{u}}$. We study the 3D-1Z problem from \cite{STD+:96}, which considers a 3D laminar steady-state flow around a cylinder. No indicator fields are modeled in this example because this configuration does not include isolated volumes of fluid and the geometry does not change during the analysis. 

The problem setup is shown in Figure \ref{fig:setup_dfg_cylinder}. The inflow condition is:
%
\begin{equation}
	\label{eq:dfg_cylinder_inflow}
	\hat{u}_{x}\left(0, y, z\right) = 16 \ u_{m} \ y \ z \ \left(\frac{\left(0.41 - y\right)\left(0.41 - z\right)}{0.41^4}\right) \qquad \hat{u}_{y} = \hat{u}_{z} = 0 \ ,
\end{equation}
%
where the maximum inflow velocity is set to $u_{m}=0.45$. A traction-free boundary condition is imposed on the outlet. No-slip boundary conditions are imposed on the surface of the cylinder and on all other planes. The characteristic velocity is defined by the average inflow velocity, that is $u_{c}=4/9 \ u_m=0.2$. The characteristic length, $L_{c}$, is defined as the diameter of the cylinder, $0.1$, which yields a Reynolds number of $20$. Note that in our CutFEM approach the cylinder is represented by the level set function \eqref{eq:level_set_cylinder}. 

The following quantities are used to monitor mesh convergence: the drag coefficient \eqref{eq:drag_coeff_criteria} around the cylinder, $c_{D}$, and the total pressure drop \eqref{eq:pressure_drop_criteria} between the inlet and the outlet planes, $\mathcal{T}_{\mathrm{in}} - \mathcal{T}_{\mathrm{out}}$. The numerical solutions provided in the study by \cite{STD+:96} give the values of $6.05$ and $6.25$ as the lower and upper bounds for the drag, respectively. These bounds were computed from the numerical results provided by several research groups through different numerical schemes, such as Finite Difference, Finite Volume, and Finite Element Methods, among others, and mesh convergence studies. No reference solution is provided for the total pressure difference; however, given that this measure is used as the objective of the optimization problems in several numerical examples ahead, we include it in the convergence study.

The mesh for the body-fitted problem utilizes a boundary layer around the cylinder, with $128$ elements on the surface and $64$ layers. The width of the first layer is $1.0832 \times 10^{-9}$, and the exponential growth factor of each subsequent layer is $1.2$. The number of elements on the inlet and outlet surfaces is $32 \times 32$ and $16 \times 16$, respectively. We do not consider symmetry, and solve the problem on the entire flow domain. For the body-fitted mesh, all boundary conditions are imposed in the strong form, and no face-oriented ghost-penalty formulation is applied. The total number of degrees-of-freedom is $1,556,066$. The remaining parameters used for this problem are shown in Table \ref{tab:dfg_cylinder_params_bf}.
\begin{table}
	\centering
	\begin{tabular*}{0.625\linewidth}{@{\extracolsep{\fill}} l l}
	\hline
	                                   & Value \\
  \hline
	Characteristic velocity            & $u_{c} = 0.2$ \\
	Characteristic length              & $L_{c} = 0.1$ \\
	Dynamic viscosity                  & $\mu = 0.001$ \\
	Density                            & $\rho = 1$ \\
  Nitsche velocity penalty           & $\alpha_{\mathrm{N}, \bm{u}} = 0$ \\
  Viscous ghost-penalty              & $\alpha_{\mathrm{GP}, \mu} = 0$ \\
  Pressure ghost-penalty             & $\alpha_{\mathrm{GP}, p} = 0$ \\
  Convective ghost-penalty           & $\alpha_{\mathrm{GP}, \bm{u}} = 0$ \\
	Pressure constraint parameter      & $k_{p} = 0$ \\
  \hline
	\end{tabular*}
	\caption{Problem parameters for the flow around a cylinder example (body-fitted).}
	\label{tab:dfg_cylinder_params_bf}
\end{table}

\begin{figure}
	\centering
	\subfloat[Cross-section of the original structured mesh, along the $z=0.205$ symmetry plane. The red lines denote the zooming area for the figures below.]{
		\includegraphics[width=0.75\linewidth]{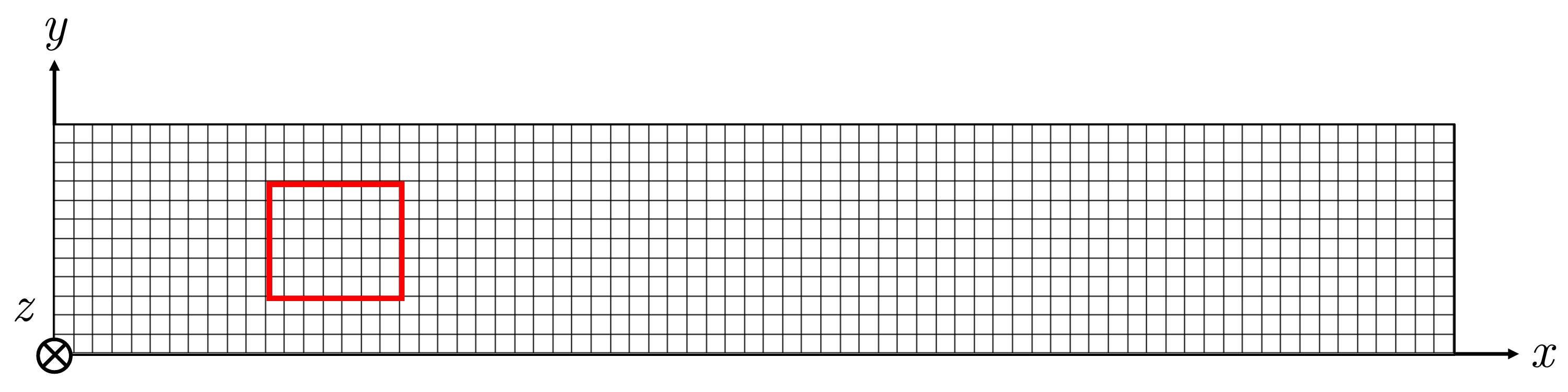}
		\label{fig:dfg_cylinder_mesh_xfem}
	} \\
	\subfloat[$h = 1.14 \times 10^{-2}$] {
		\includegraphics[width=0.175\linewidth]{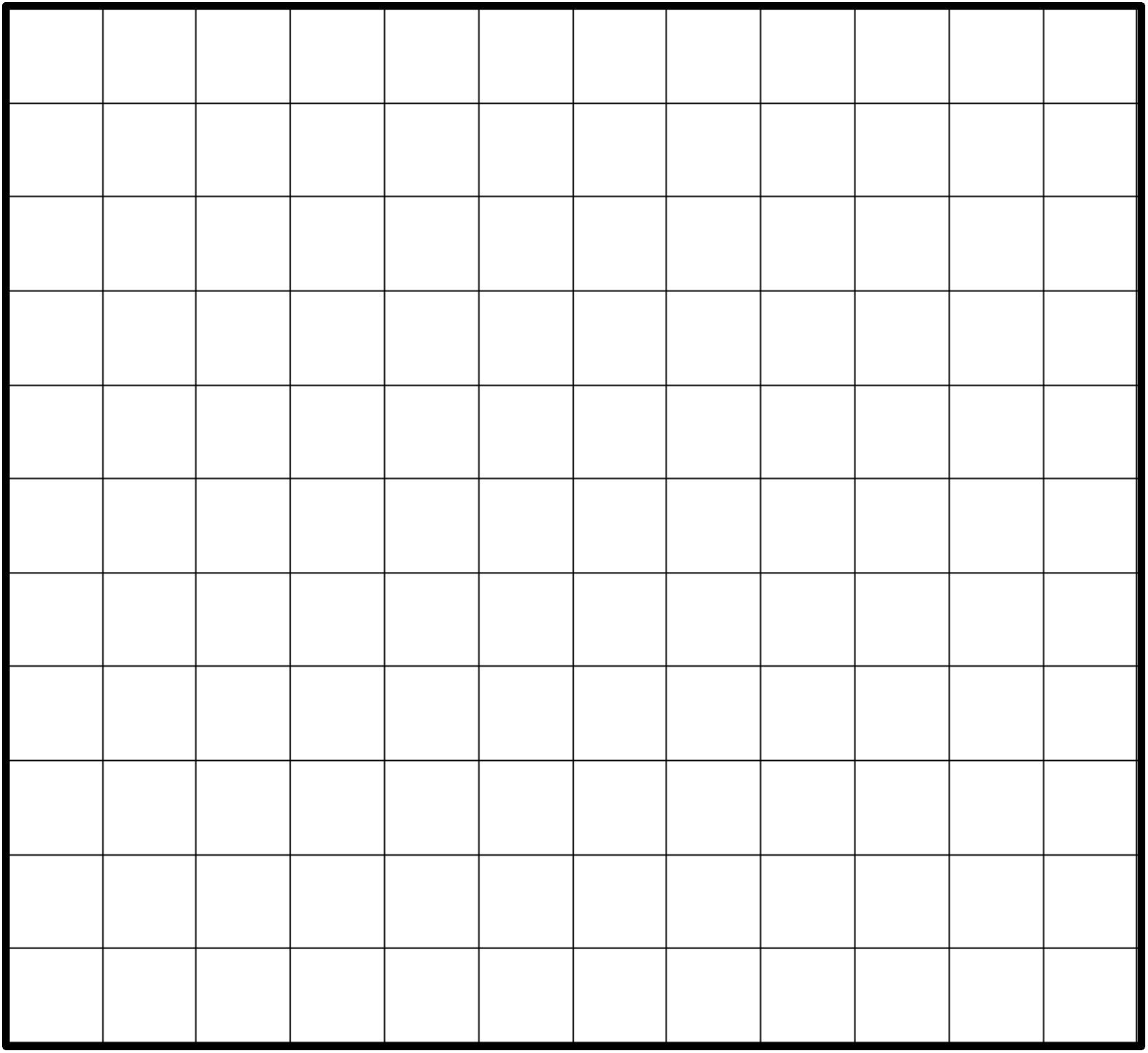}
		\label{fig:dfg_cylinder_mesh_xfem_cylbox}
	} \qquad
	\subfloat[$h = 3.79 \times 10^{-3}$] {
		\includegraphics[width=0.175\linewidth]{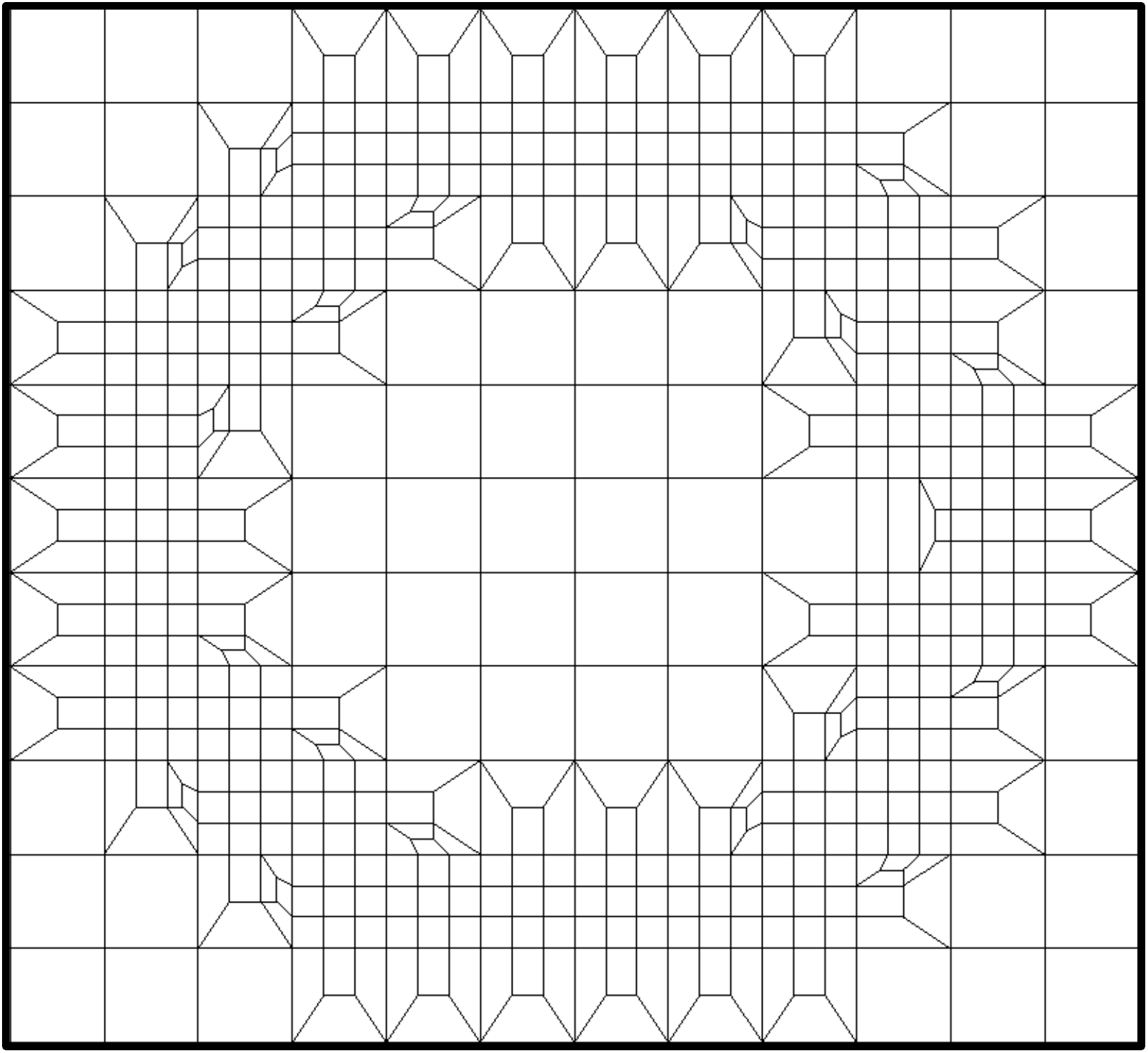}
		\label{fig:dfg_cylinder_mesh_xfem_fine}
	} \qquad
	\subfloat[$h = 1.27 \times 10^{-3}$] {
		\includegraphics[width=0.175\linewidth]{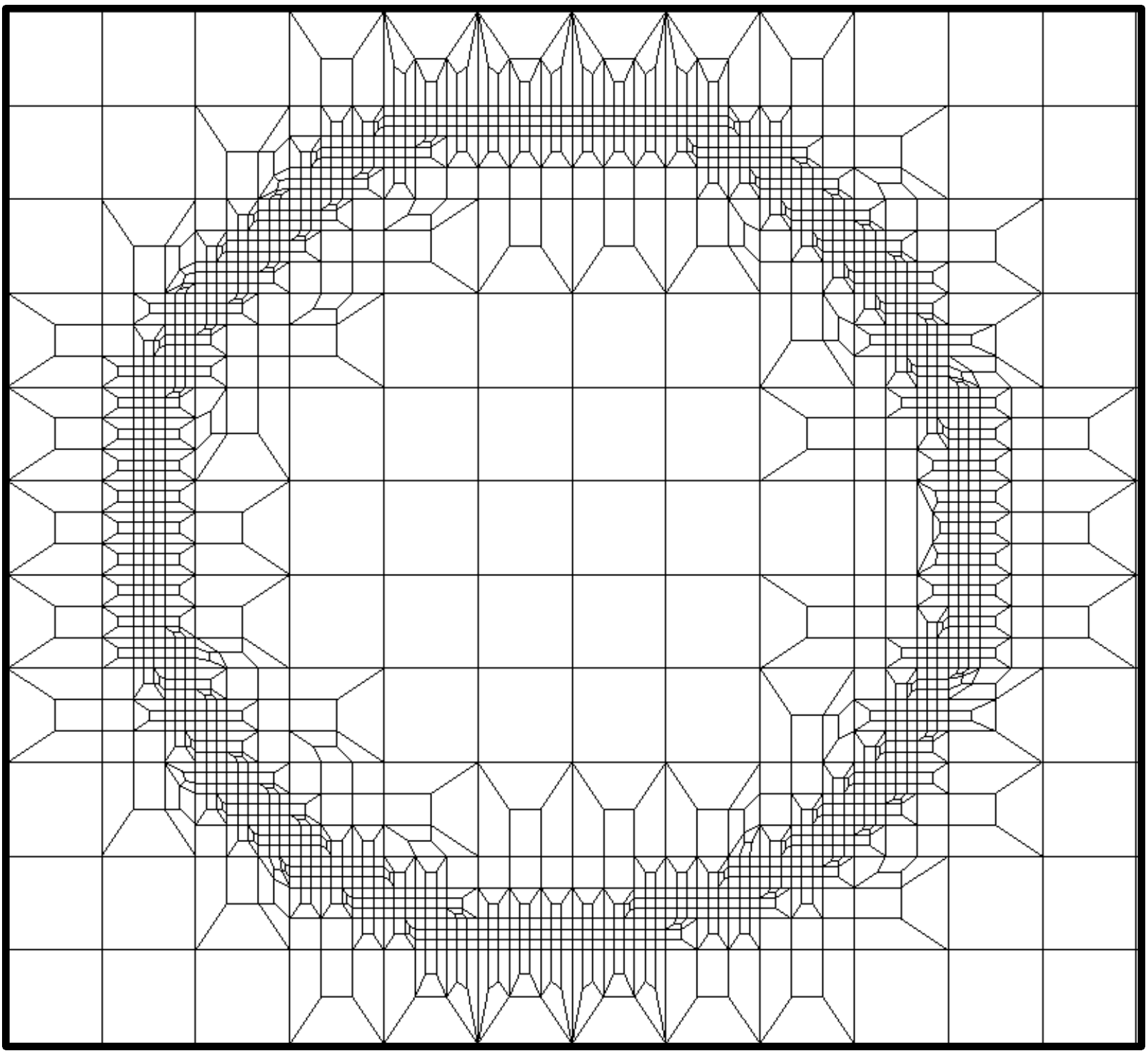}
		\label{fig:dfg_cylinder_mesh_xfem_finer}
	} \qquad
	\subfloat[$h = 4.22 \times 10^{-4}$] {
		\includegraphics[width=0.175\linewidth]{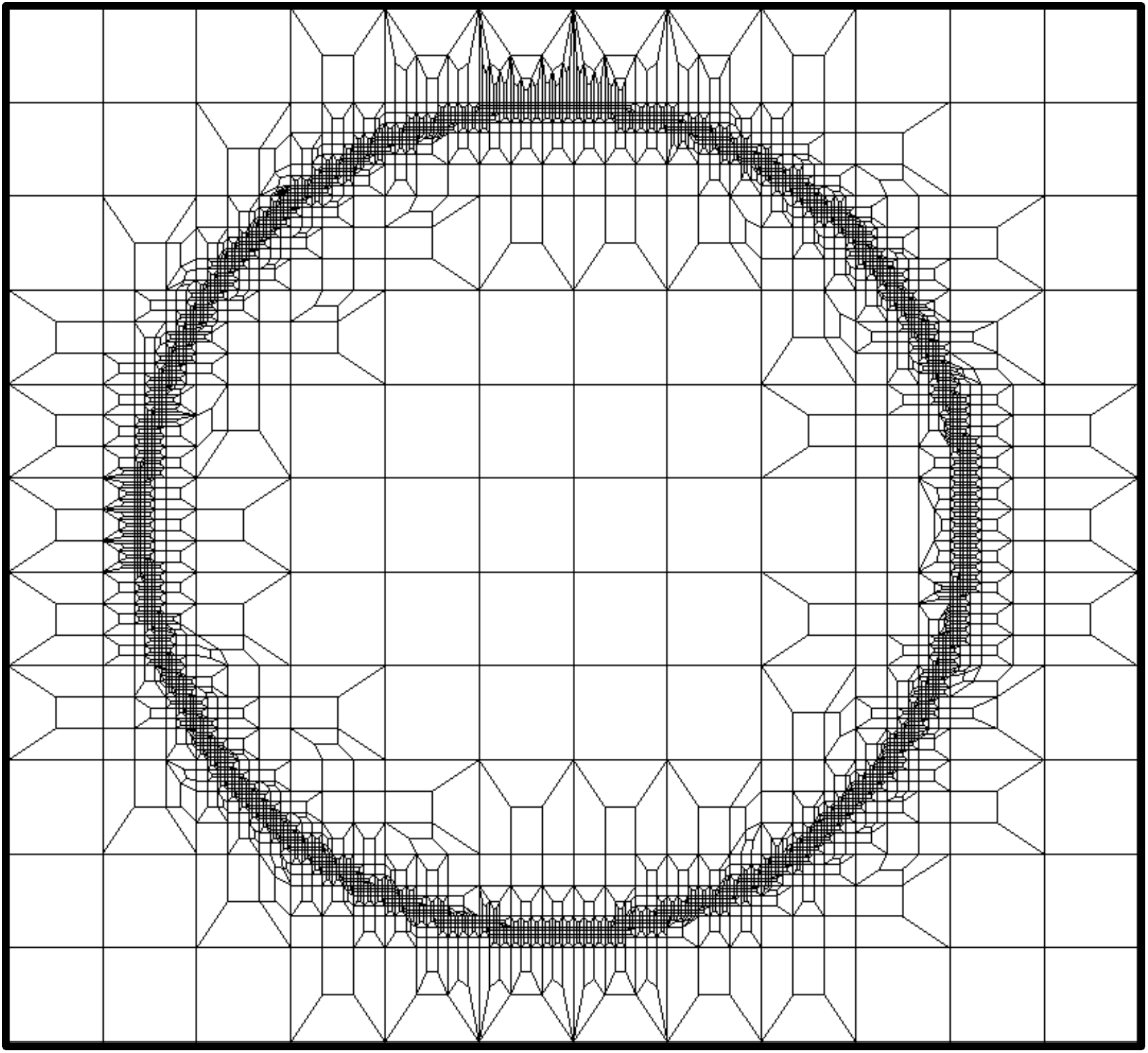}
		\label{fig:dfg_cylinder_mesh_xfem_finest}
	}
	\caption{Mesh refinement levels using a local hierarchical mesh refinement for validation of the CutFEM framework. The $h$ values represent the minimum element sizes in the mesh.}
	\label{fig:dfg_cylinder_mesh}
\end{figure}

The meshes for the CutFEM convergence study are constructed using a local hierarchical refinement strategy along the fluid-solid interface. We deliberately do not perform adaptive mesh refinement nor generate boundary layer meshes in the CutFEM approach, as we do not consider these techniques in the topology optimization problems ahead. The first level of mesh refinement has an element size of $h=1.14 \times 10^{-2}$ for elements located at $x<1$ and $h=3.41667 \times 10^{-2}$ for all other elements. Subsequently, we only refine elements that are intersected by the surface of the cylinder; we perform this process thrice. In consequence, we generate the four refinement levels shown in Figure \ref{fig:dfg_cylinder_mesh}. Note that the hierarchical mesh refinement utilized here does not provide the same resolution as the boundary layer meshing scheme above. 

The hierarchical mesh refinement scheme leads to a larger number of elements compared to the boundary layer approach. To reduce the computational cost, we model half of the flow domain and impose symmetry boundary conditions by setting $\hat{u}_{z} = 0$ along the plane $z=0.205$. The inflow and outflow conditions are the same as in the body-fitted problem setup. All boundary conditions are enforced weakly. The total number of degrees-of-freedom for the finest mesh is $8,019,736$. The face-oriented ghost-penalty methods increase the bandwidth of the sparse matrix of the system, which may increase the linear solve time; see also discussion on the computational costs below. The remaining parameters used for this problem are shown in Table \ref{tab:dfg_cylinder_params}.

\begin{table}
	\centering
	\begin{tabular*}{\linewidth}{@{\extracolsep{\fill} } l l}
	\hline
	                                   & Value \\
    \hline
	Element size                       & $h = \left\{1.14 \times 10^{-2}, 3.79 \times 10^{-3}, 1.27 \times 10^{-3}, 4.22 \times 10^{-4}\right\}$ \\
	Characteristic velocity            & $u_{c} = 0.2$ \\
	Characteristic length              & $L_{c} = 0.1$ \\
	Dynamic viscosity                  & $\mu = 0.001$ \\
	Density                            & $\rho = 1$ \\
    Nitsche velocity penalty           & $\alpha_{\mathrm{N}, \bm{u}} = \left\{10, 10^{2}, 10^{3}, 10^{4}\right\}$ \\
    Viscous ghost-penalty              & $\alpha_{\mathrm{GP}, \mu} = 0.05$ \\
    Pressure ghost-penalty             & $\alpha_{\mathrm{GP}, p} = 0.005$ \\
    Convective ghost-penalty           & $\alpha_{\mathrm{GP}, \bm{u}} = 0.05$ \\
    Pressure constraint parameter      & $k_{p} = 0$ \\
    \hline
	\end{tabular*}
	\caption{Problem parameters for the flow around a cylinder example (CutFEM).}
	\label{tab:dfg_cylinder_params}
\end{table}

\begin{figure}
	\centering
	\subfloat[Drag coefficient.]{
		\includegraphics[width=0.45\linewidth]{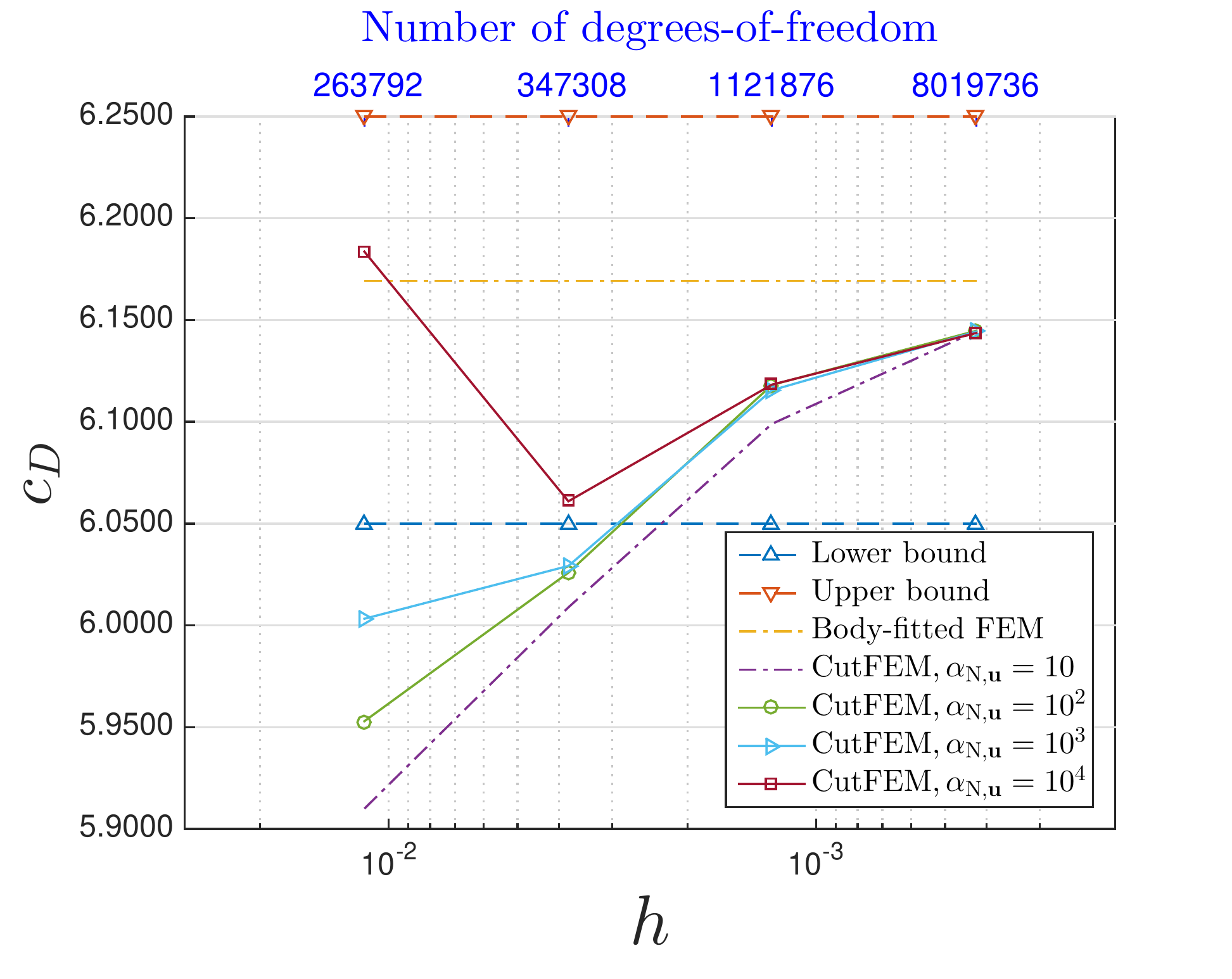}
		\label{fig:ghost_validation_cd}
	}
	\subfloat[Total pressure difference.]{
		\includegraphics[width=0.45\linewidth]{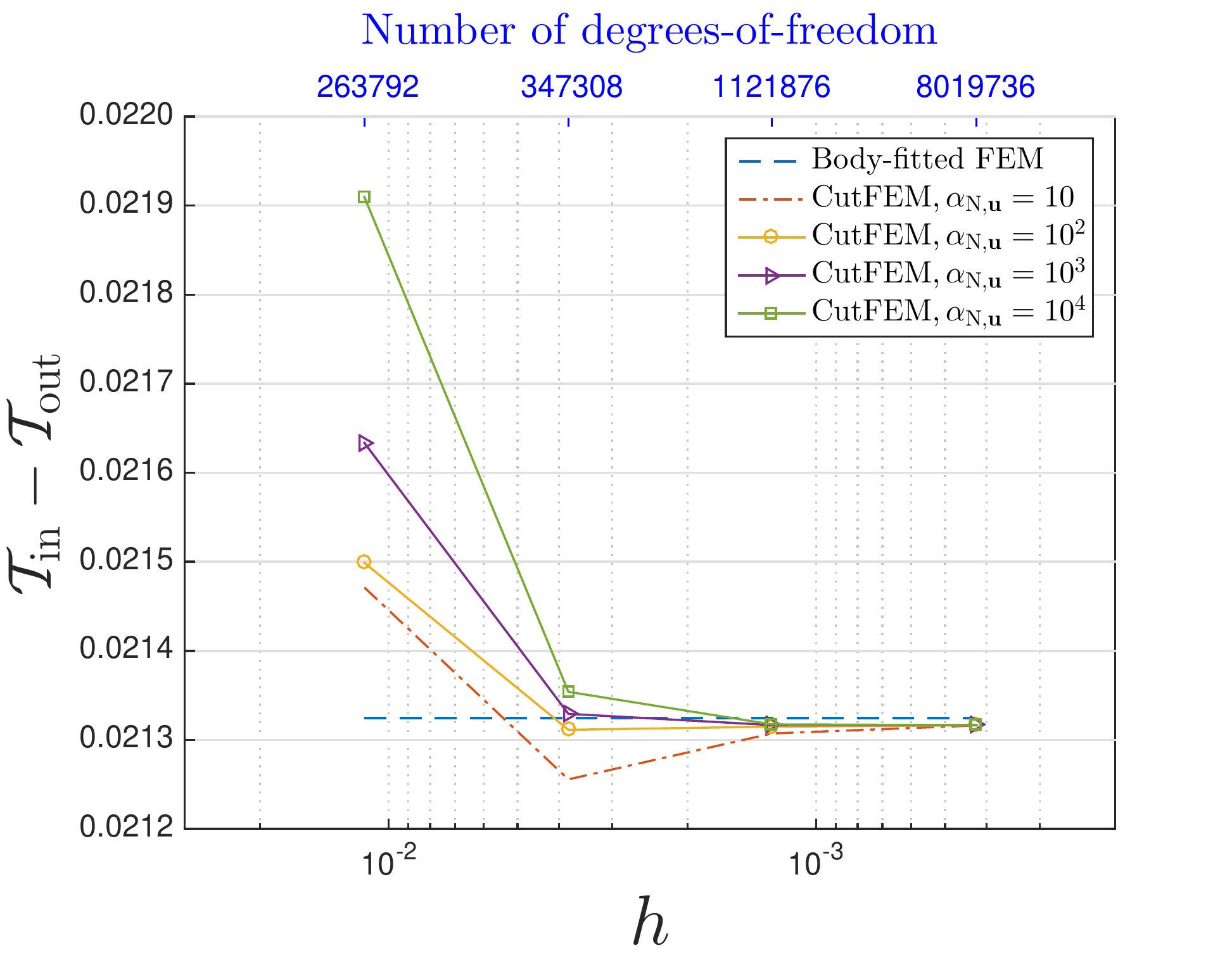}
		\label{fig:ghost_validation_tp}
	} \\
	\subfloat[Mass flow rate through surface of cylinder.]{
		\includegraphics[width=0.45\linewidth]{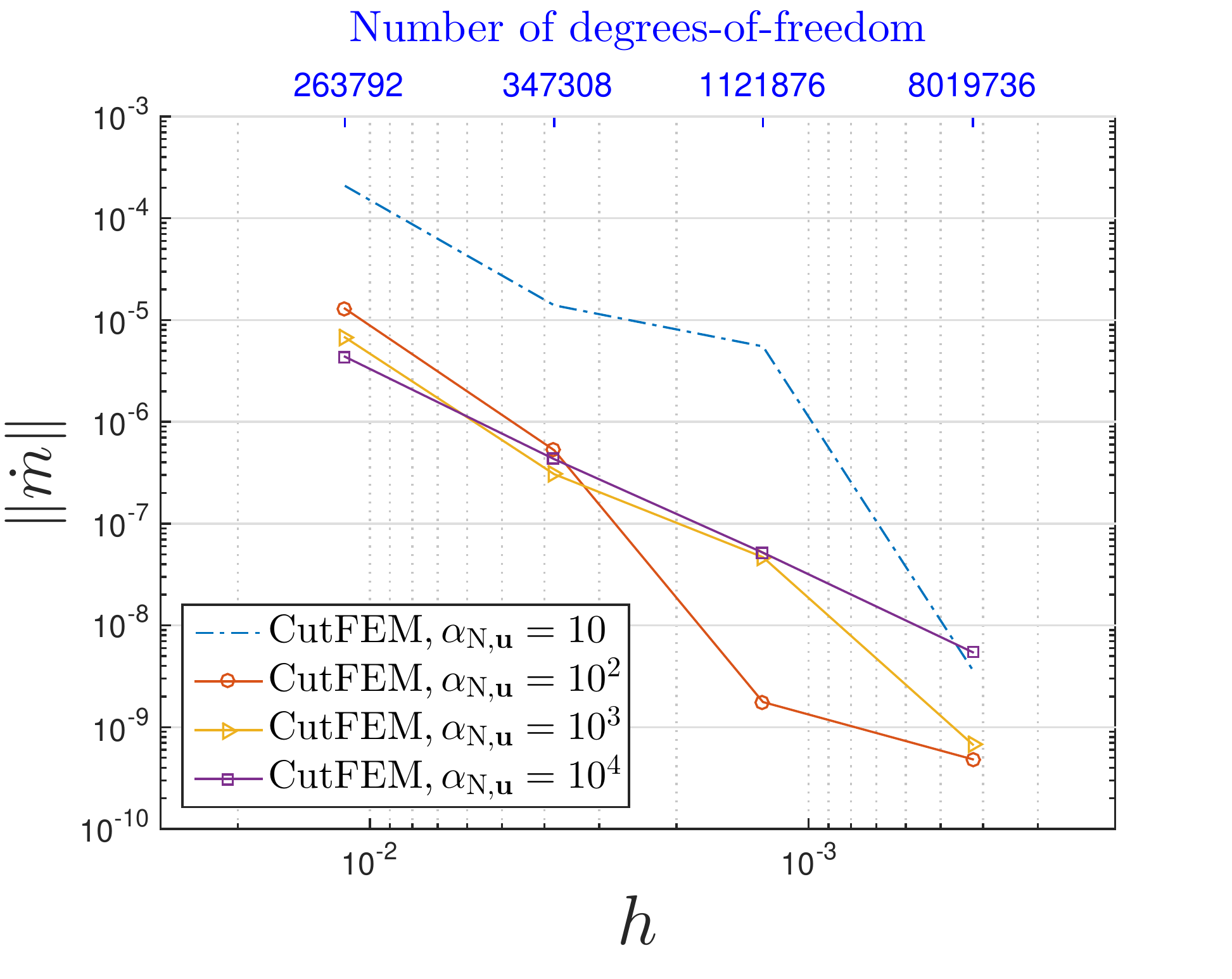}
		\label{fig:ghost_validation_mass_flow}
	}
	\caption{Comparison of the body-fitted and CutFEM solutions for the flow around a cylinder example.}
	\label{fig:ghost_validation}
\end{figure}

The results for the body-fitted and the CutFEM problems are shown in Figure \ref{fig:ghost_validation}. The drag coefficient for the body-fitted problem is $6.169$, well within the lower and upper bounds established in the study by \cite{STD+:96}. The total pressure difference is $0.0213$. The results of the CutFEM analysis are shown in Figure \ref{fig:ghost_validation_cd}. Although the drag coefficient is not fully converged, due to the lack of sufficiently resolving the boundary layer, the values for the finest mesh are well within the lower and upper bounds. Due to limited computational resources, finer CutFEM meshes could not be considered. The results of the total pressure drop are shown in Figure \ref{fig:ghost_validation_tp} and display a higher convergence rate than the drag coefficient. The total pressure drop characterizes the global flow solution, while the drag is a local measure along the cylinder surface that depends on the spatial gradients of the flow field. The total pressure difference reaches the same solution as its body-fitted counterpart. Spurious mass flow rates through the surface of the cylinder are shown in Figure \ref{fig:ghost_validation_mass_flow}. We observe that as the mesh is refined, and the value of the Nitsche velocity penalty is increased, the spurious mass flow converges to zero within numerical precision.

As the mesh is refined, the influence of the Nitsche velocity parameter in \eqref{eq:nitsche_penalty_param} vanishes. The relative maximum difference between the body-fitted problem and the CutFEM solutions at the finest mesh is $0.4\%$ for the drag coefficient and $0.03\%$ for the total pressure drop. These results suggest that our CutFEM framework is sufficiently accurate for topology optimization purposes. However, our discretization scheme may suffer from inaccuracies in local quantities that strongly depend on the resolution of boundary layer phenomena. If such quantities are used in the formulation of the optimization problem, we recommend that the performance of the optimized design is verified using a body-fitted mesh with a resolved boundary layer.

To provide insight into the computational costs, wall clock timing results for individual steps in the forward analysis with the body-fitted mesh and the CutFEM mesh with $h = 1.27 \times 10^{-3}$ are given in Table \ref{tab:dfg_cylinder_computational_time}, using $60$ and $120$ cores. These particular meshes are chosen as they have comparable numbers of elements and unconstrained DOFs. For the CutFEM analysis, the case with $\alpha_{\mathrm{N}, \bm{u}} = 10$ is considered; other values for $\alpha_{\mathrm{N}, \bm{u}}$ lead to similar results. The wall clock times for assembling and solving the linearized systems are averaged over all Newton iterations. For both, the body-fitted and the CutFEM analysis the same convergence criteria for the nonlinear and linear problems are applied, as defined above; the ILUT fill is set to $7.0$. The computations were performed on a cluster equipped with Quad-Core AMD Opteron 2378 Processors. These results show that the computational effort for building and updating the CutFEM model is insignificant when compared to the time spent for solving the nonlinear system. Assembly of the CutFEM model is more costly than the body-fitted mesh due to the integration of the ghost-penalty terms and the increased number of integration points of intersected elements. However, the assembly operation scales well with the number of cores. The timing results presented in Table \ref{tab:dfg_cylinder_computational_time} are representative of the computational cost of all the other examples studied in this paper, including cost per time step of transient problems and the cost for solving adjoint problems.

\begin{table}
	\centering
	\def\arraystretch{1.5}
	\begin{tabular*}{1.0\linewidth}{@{\extracolsep{\fill}} l c c c c}
	\hline
	                               & \multicolumn{2}{c}{Body-fitted}  & \multicolumn{2}{c} {CutFEM}      \\
	                               & \multicolumn{2}{c}{}             & \multicolumn{2}{c} {$\left({h = 1.27 \times 10^{-3}}\right)$}      \\
	\hline 
	Number of nodes                & \multicolumn{2}{c}{416,576}      & \multicolumn{2}{c}{380,492}   \\
	Number of elements             & \multicolumn{2}{c}{397,792}      & \multicolumn{2}{c}{372,492}   \\
	Number of DOFs                 & \multicolumn{2}{c}{1,556,066}    & \multicolumn{2}{c}{1,121,876} \\
	Number of intersected elements & \multicolumn{2}{c}{0}            & \multicolumn{2}{c}{51,192}    \\
  \cline{2-3}\cline{4-5}
	Number of cores used                   &        60 &      120             &      60 &         120         \\
	Number of Newton steps                 &         5 &       5              &       6 &           5         \\
	Time for updating CutFEM model $[s]$   &         0 &       0              &    24.1 &        17.6         \\
	Average assembly time          $[s]$   &       8.5 &     4.3              &   22.56 &        12.3         \\
	Average GMRES iterations               &     190.2 &   226.6              &    71.8 &       123.4         \\
	Average ILUT/GMRES time         $[s]$  &     106.4 &    68.1              &   107.3 &        67.7         \\
	\hline
	\end{tabular*}
	\caption{Computational cost of forward analysis for the body-fitted mesh and the CutFEM mesh with $h = 1.27 \times 10^{-3}$ for the flow around a cylinder example.}
	\label{tab:dfg_cylinder_computational_time}
\end{table}


\subsection{Verification of the Average Pressure Constraint}
\label{sec:pipebend}

\begin{figure}
	\centering
	\includegraphics[width=0.4\linewidth]{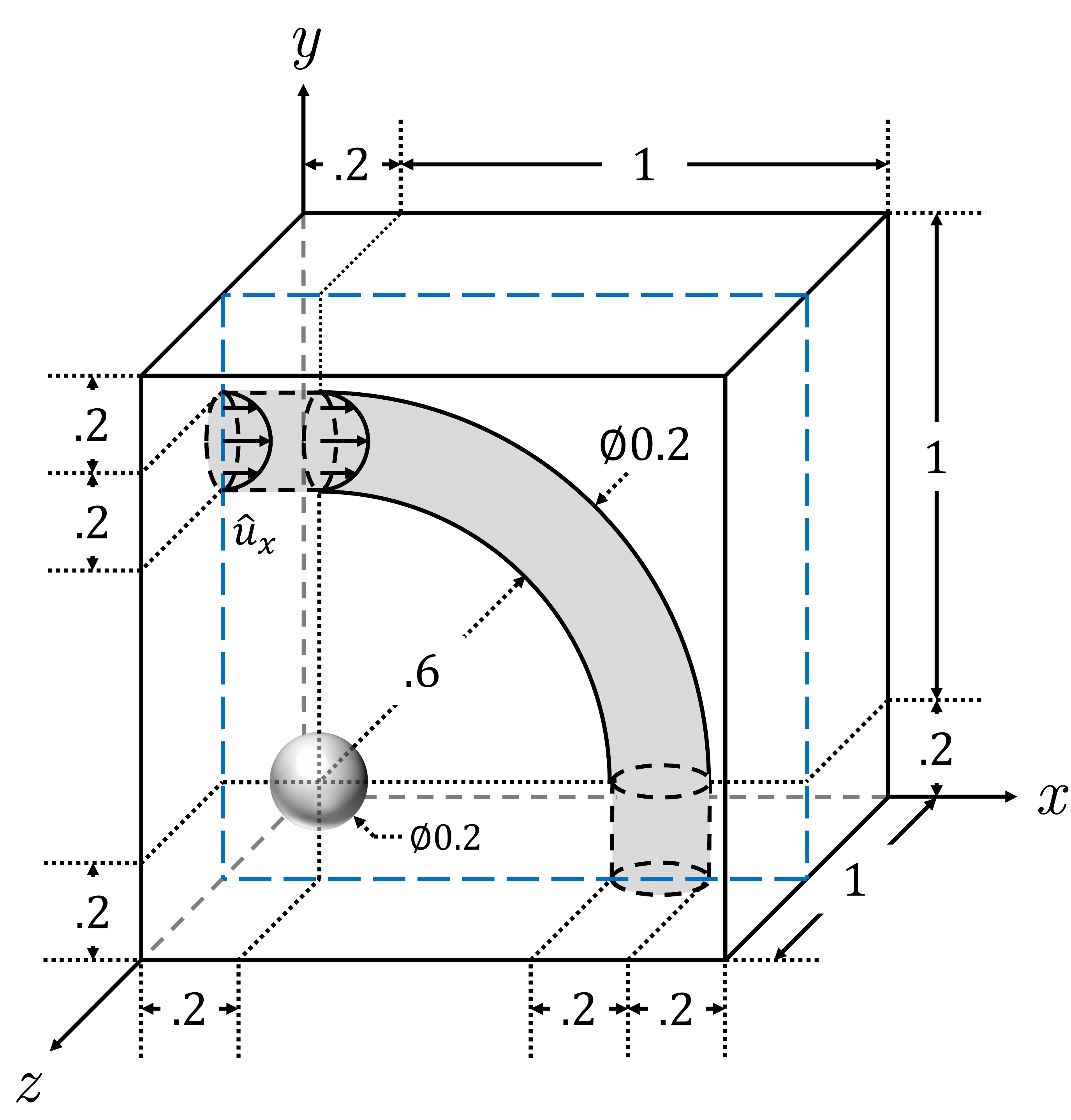}
	\caption{Problem setup for the bent pipe example. Blue dashes lines denote the symmetry plane.}
	\label{fig:setup_pipebend}
\end{figure}

In this second example, we seek to verify the penalty formulation in \eqref{eq:pressure_penalty_formulation} with respect to the accuracy to which mass conservation is satisfied. We model a steady-state flow through a bent pipe and measure the relative mass flow rate difference between the inlet and the outlet. The problem setup is shown in Figure \ref{fig:setup_pipebend}. The presence of an isolated spherical inclusion of fluid causes a singular analysis problem because the absolute value of the pressure is not governed within it. We compare flow solutions with and without the spherical inclusion and study the influence of the pressure penalty parameter, $k_{p}$.

\begin{table}
	\centering
	\begin{tabular*}{0.75\linewidth}{@{\extracolsep{\fill}} l l}
	\hline
	                                      & Value \\
    \hline
    Mesh size                           & $120 \times 120 \times 50$ (half domain) \\
	  Element size                        & $h = 0.01$ \\
	  Characteristic velocity             & $u_{c} = 200$ \\
	  Characteristic length               & $L_{c} = 0.2$ \\
	  Dynamic viscosity                   & $\mu = 1$ \\
	  Density                             & $\rho = 1$ \\
    Nitsche velocity penalty            & $\alpha_{\mathrm{N}, \bm{u}} = 100$ \\
    Nitsche indicator field penalty     & $\alpha_{\mathrm{N}, \psi} = 1$ \\
    Viscous ghost-penalty               & $\alpha_{\mathrm{GP}, \mu} = 0.5$ \\
    Pressure ghost-penalty              & $\alpha_{\mathrm{GP}, p} = 0.05$ \\
    Convective ghost-penalty            & $\alpha_{\mathrm{GP}, \bm{u}} = 0.5$ \\
	  Pressure constraint parameter       & $k_{p} = \left\{10^{-8}, 10^{-6}, 10^{-4}, 10^{-2}, 1\right\}$ \\
    \hline
	\end{tabular*}
	\caption{Problem parameters for the bent pipe example.}
	\label{tab:pipebend_params}
\end{table}

The design domain is discretized by a uniform structured mesh, and the geometries of the pipe and the spherical inclusion are described by the LSFs presented in \cite{BCH+:14}. We only model half of the domain, and apply symmetry boundary conditions along the plane $z=0.5$. We study the pressure penalty formulation in \eqref{eq:pressure_penalty_formulation} and, for illustration purposes, a formulation in which the pressure penalty is applied over the entire fluid domain, rather than on the isolated ``puddles'' exclusively. 

The problem parameters are given in Table \ref{tab:pipebend_params}. The inflow condition is:
%
\begin{equation}
	\label{eq:pipebend_inflow}
	\hat{u}_{x}\left(0, y, z\right) = u_{c} \cdot \left(
	\left(-\frac{4}{L_{c}^{2}}\right) \cdot
	\left({\left(y - y_{c}\right)}^{2} + {\left(z - z_{c}\right)}^{2}\right)
	 + 1\right)\ , \qquad \hat{u}_{y} = \hat{u}_{z} = 0\ ,
\end{equation}
%
where $y_{c}$ and $z_{c}$ are the coordinates at the center of the inflow, as defined in Figure \ref{fig:setup_pipebend}. The characteristic velocity is $u_{c} = 200$, and the characteristic length is defined as the diameter of the pipe, $L_{c} = 0.2$, for a Reynolds number of $40$. A traction-free boundary condition is imposed on the outlet. No-slip boundary conditions are imposed on the surfaces of the pipe and of the sphere.

\begin{figure}
	\centering
	\includegraphics[width=0.75\linewidth]{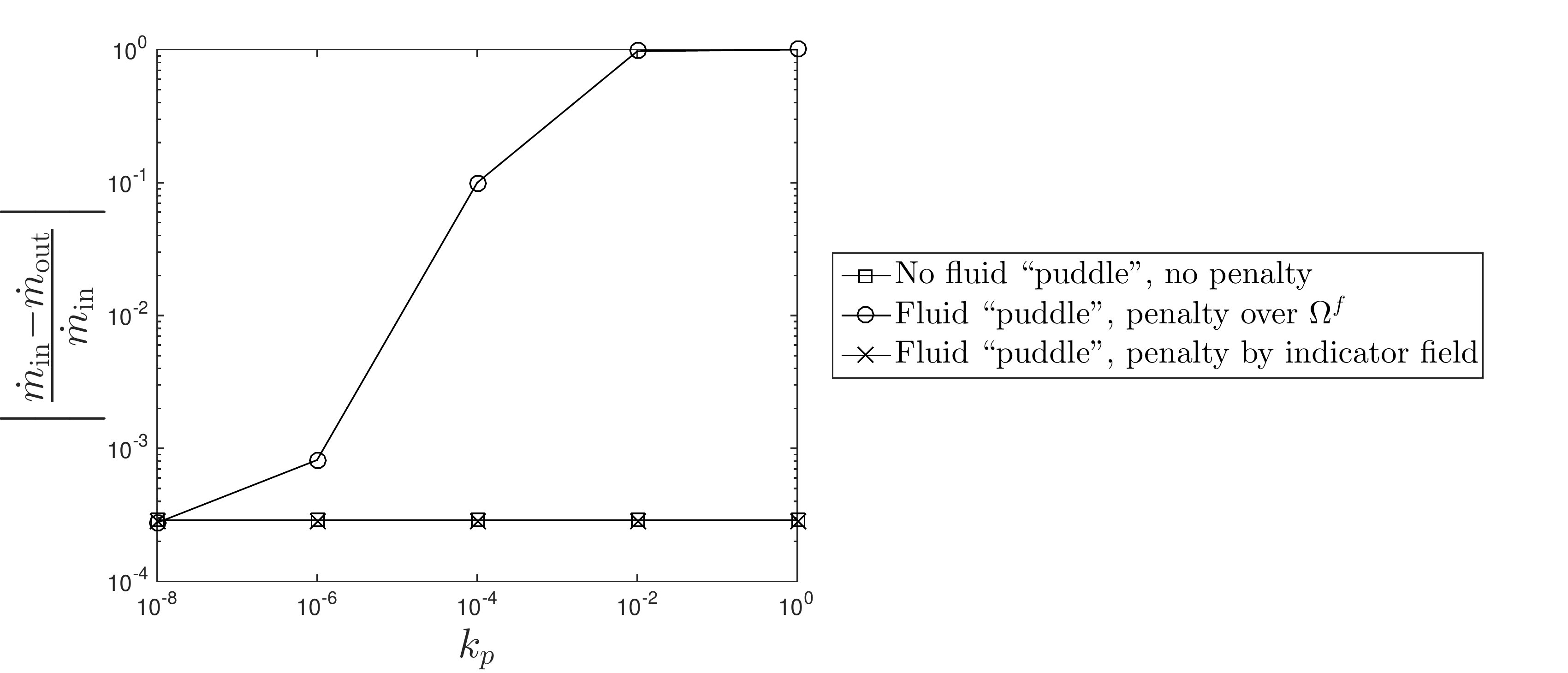}
	\caption{Comparison of the influence of the average pressure constraint on the relative mass flow rate difference between the inlet and the outlet for the bent pipe example.}
	\label{fig:pipebend_springs_comp}
\end{figure}
The relative error between the inlet and outlet mass flow rates is shown in Figure \ref{fig:pipebend_springs_comp}. The flow without the sphere results in a relative mass flux error of $0.03\%$. We obtain the same error when the sphere is included and if we apply the pressure penalty formulation \eqref{eq:pressure_penalty_formulation} exclusively to the isolated volume of fluid through the use of the indicator field in \eqref{eq:auxiliary_species_field_residual}. The error is insensitive to the penalty parameter value, which is set to $k_{p} = 1$ for all numerical examples below. Applying the penalty formulation over the entire domain, similar to the approach used by \cite{VM:14} for linear elasticity problems, can cause a significant error in the mass conservation if a large value of $k_{p}$ is chosen; see Figure \ref{fig:pipebend_springs_comp}.


\subsection{Design of a Manifold with Multiple Outlets}
\label{sec:multiple_outlets}

\begin{figure}
	\centering
	\includegraphics[width=0.4\linewidth]{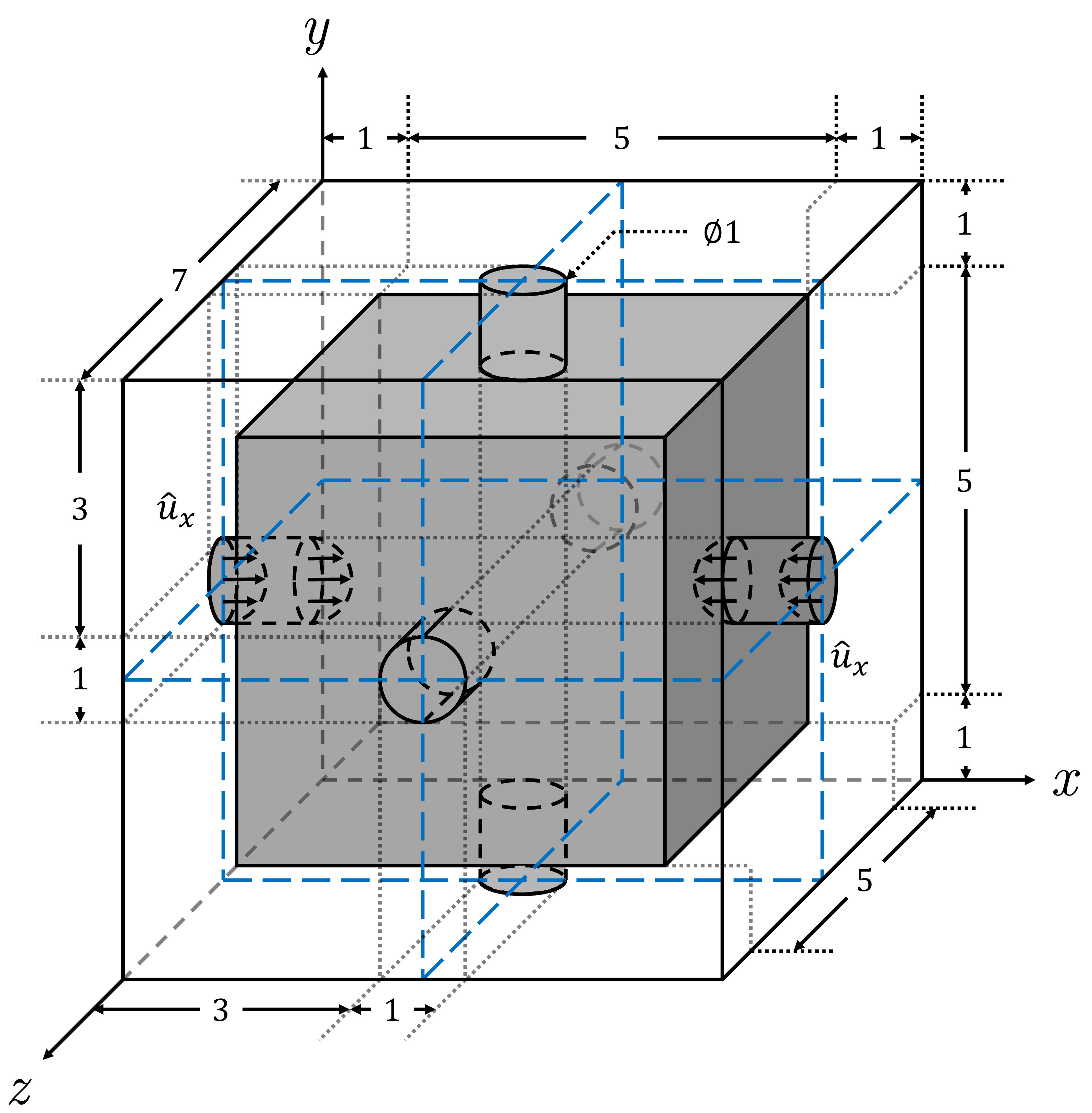}
	\caption{Problem setup for the multiple outlets example. Blue dashes lines denote the symmetry planes.}
	\label{fig:setup_skeleton}
\end{figure}

In this example, we apply the CutFEM framework to the design of a steady-state flow bend, with multiple inlets and outlets. We seek to minimize the total pressure drop between the inlets and the outlets at steady-state while controlling the mass flow rates through the outlet ports. The example is the 3D analog to the 2D problem found in \cite{PWE+:10}. The problem setup is shown in Figure \ref{fig:setup_skeleton}. The design domain has two inlets (on the left and right), and four outlets (on all other planes). The inflow condition is formulated in the same way as \eqref{eq:pipebend_inflow}, the characteristic velocity is $u_{c} = 200$, and the characteristic length is $L_{c}=1$. Traction-free boundary conditions are imposed on the outlets. No-slip boundary conditions are imposed on the fluid-solid interface. We only model an eighth of the domain, and symmetry boundary conditions are imposed on the  planes with $x=3.5$, $y=3.5$, and $z=3.5$.

In addition to minimizing the total pressure drop, we penalize the surface area of the fluid-solid interface. The use of the surface area as a contribution to the objective function has been applied previously to species transport topology optimization by \cite{MM:14} to improve the smoothness of the final design, and to regularize the optimization problem. The objective is defined as:
%
\begin{equation}
	\label{eq:multiple_outlets_objective}
	\mathcal{Z} = \frac{\sum\limits^{2}_{i = 1}\mathcal{T}_{\mathrm{in},i} - \sum\limits^{4}_{i = 1}\mathcal{T}_{\mathrm{out},i}}{\left\Vert\sum\limits^{2}_{i = 1}\mathcal{T}^{0}_{\mathrm{in},i} - \sum\limits^{4}_{i = 1}\mathcal{T}^{0}_{\mathrm{out},i}\right\Vert} + w_{\mathcal{S}} \frac{\mathcal{S}}{\left\Vert\mathcal{S}^{0}\right\Vert} \ ,
\end{equation}
%
where the superscript ``0'' denotes the values of the initial design, the subscript $i$ denotes the $i$-th inlet or outlet, and $w_{\mathcal{S}}$ is a constant scaling factor. The design is subject to a $5\%$ volume constraint of the fluid domain to suppress trivial solutions, and to promote the formation of distinct fluid channels:
%
\begin{equation}
	\label{eq:multiple_outlets_vol_const}
	g_{1} = \frac{\mathcal{V}^{f}}{0.05\left(\mathcal{V}^{f} + \mathcal{V}^{s}\right)} - 1\ .
\end{equation}
%
Further, we wish to impose a constraint such that the amount of mass flow exiting through each outlet is the same. Given that the GCMMA algorithm does not allow equality constraints, we impose inequality constraints with lower and upper limits on the mass flow rates. The upper and lower bounds are set to $25\% \pm 1.25\%$, respectively, where the tolerance value of $\pm 1.25\%$ was chosen in order to not overconstrain the optimization problem. The constraints are defined as follows:
%
\begin{align}
	\label{eq:multiple_outlets_mass_const}
	g_{i + 1} &= 1 - \frac{\dot{m}_{\mathrm{out}, i}}{\left(23.75\%\right) \left(\dot{m}_{\mathrm{in},1} + \dot{m}_{\mathrm{in},2}\right)}\ , \quad & i = 1 \dots N_{\mathrm{out}}\ , \\
	g_{i + 5} &= \frac{\dot{m}_{\mathrm{out}, i}}{\left(26.25\%\right) \left(\dot{m}_{\mathrm{in},1} + \dot{m}_{\mathrm{in},2}\right)} - 1\ , \quad & i = 1 \dots N_{\mathrm{out}}\ ,
\end{align}
%
where $\dot{m}_{\mathrm{in},i}$ and $\dot{m}_{\mathrm{out},i}$ are the mass flow rates at the $i$-th inlet and $i$-th outlet, respectively, and $N_{\mathrm{out}}$ is the number of outlets. 

\begin{table}
	\centering
	\begin{tabular*}{0.75\linewidth}{@{\extracolsep{\fill}} l l}
	\hline
	                                    & Value \\
  \hline
  Mesh size                           & $56 \times 56 \times 56$ (eighth of the domain) \\
	Element size                        & $h = 0.0625$ \\
	Characteristic velocity             & $u_{c} = 200$ \\
	Characteristic length               & $L_{c} = 1$ \\
	Dynamic viscosity                   & $\mu = 1$ \\
	Density                             & $\rho = 1$ \\
  Nitsche velocity penalty            & $\alpha_{\mathrm{N}, \bm{u}} = 100$ \\
  Nitsche indicator field penalty     & $\alpha_{\mathrm{N}, \psi} = 1$ \\
  Viscous ghost-penalty               & $\alpha_{\mathrm{GP}, \mu} = 0.5$ \\
  Pressure ghost-penalty              & $\alpha_{\mathrm{GP}, p} = 0.05$ \\
  Convective ghost-penalty            & $\alpha_{\mathrm{GP}, \bm{u}} = 0.5$ \\
	Pressure constraint parameter       & $k_{p} = 1$ \\
	Surface area scaling weight         & $w_{\mathcal{S}} = 0.01$ \\
	Volume constraint                   & $5\%$ \\
	Number of design variables          & $68,921$ (eighth of the domain) \\
	Design variables bounds             & $s^{L}_{i} = -0.03125$, $s^{U}_{i} = +0.03125$ \\
	Smoothing filter radius             & $r_{\phi} = 2.4h$ \\
  \hline
	\end{tabular*}
	\caption{Problem parameters for the multiple outlets example.}
	\label{tab:multiple_outlets_params}
\end{table}

\begin{figure}
	\centering
	\subfloat[Fluid domain, $\mathrm{\Omega}^{f}$]{
		\includegraphics[width=0.3\linewidth]{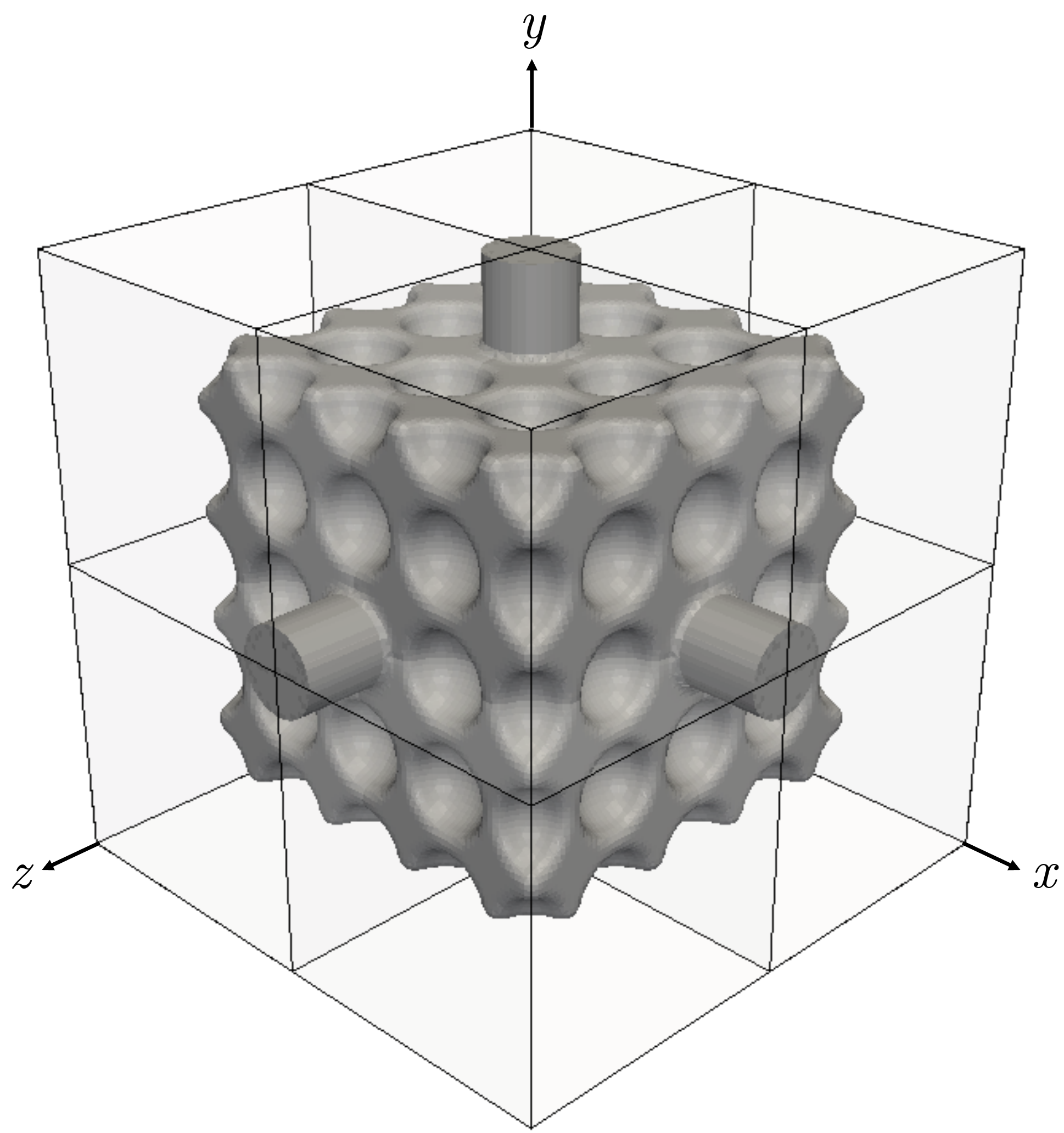}
		\label{fig:multiple_outlets_initial_-_design}
	}
	\subfloat[Half of the solid domain, $\mathrm{\Omega}^{s}$]{
		\includegraphics[width=0.3\linewidth]{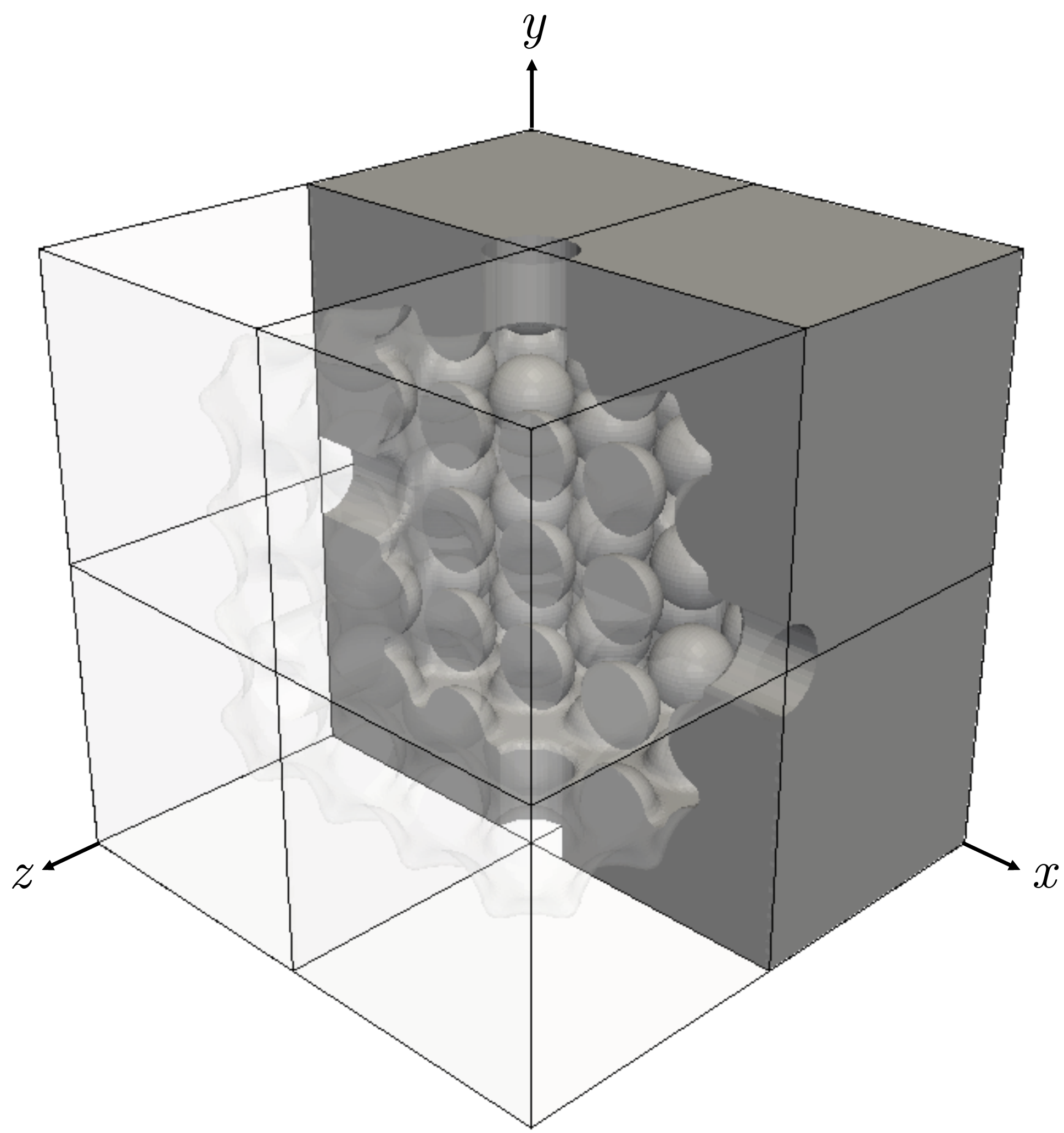}
		\label{fig:multiple_outlets_initial_+_design}
	}
	\caption{Initial design for the multiple outlets example.}
	\label{fig:multiple_outlets_initial_design}
\end{figure}

The remaining parameters are given in Table \ref{tab:multiple_outlets_params}. The design domain is initialized with $5 \times 5 \times 5$ spherical solid inclusions of radii $0.5$, as shown in Figure \ref{fig:multiple_outlets_initial_design}. In our experience, the flow topology optimization problems studied here are rather insensitive to the initial design as long as the number of inclusions is sufficiently large. Additional mechanisms for seeding solid inclusions could be added to the proposed CutFEM framework, such as topological derivatives \cite{SMN+:16}; however, they are outside the scope of this study.

\begin{figure}
	\centering
	\includegraphics[width=0.35\linewidth]{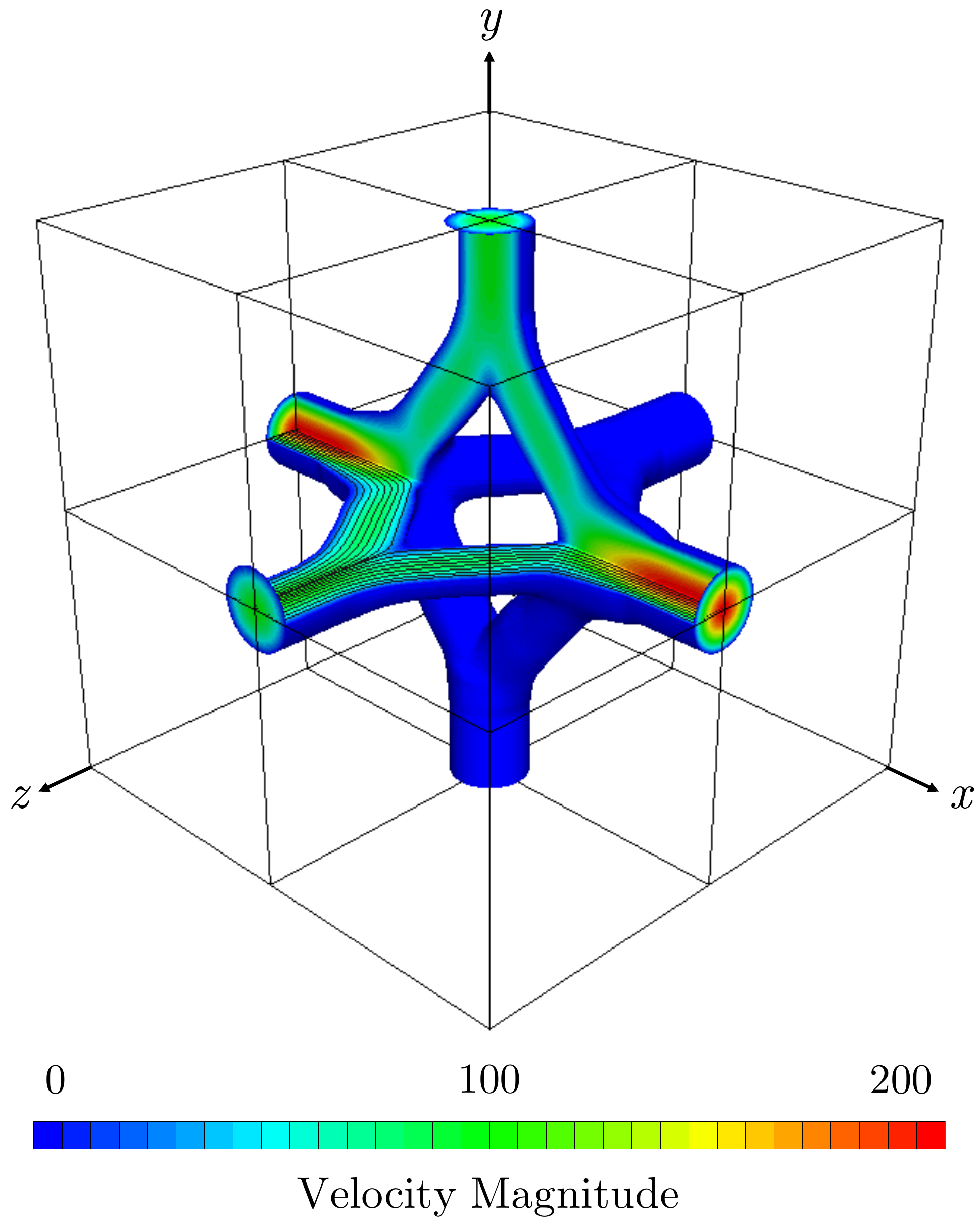}
	\caption{Velocity magnitude, with streamlines, of the optimized material layout for the multiple outlets example. A section of the design was removed for visualization purposes.}
	\label{fig:multiple_outlets_final_design}
\end{figure}

\begin{figure}
	\centering
	\subfloat[Objective and volume constraint.]{
		\includegraphics[width=0.32\linewidth]{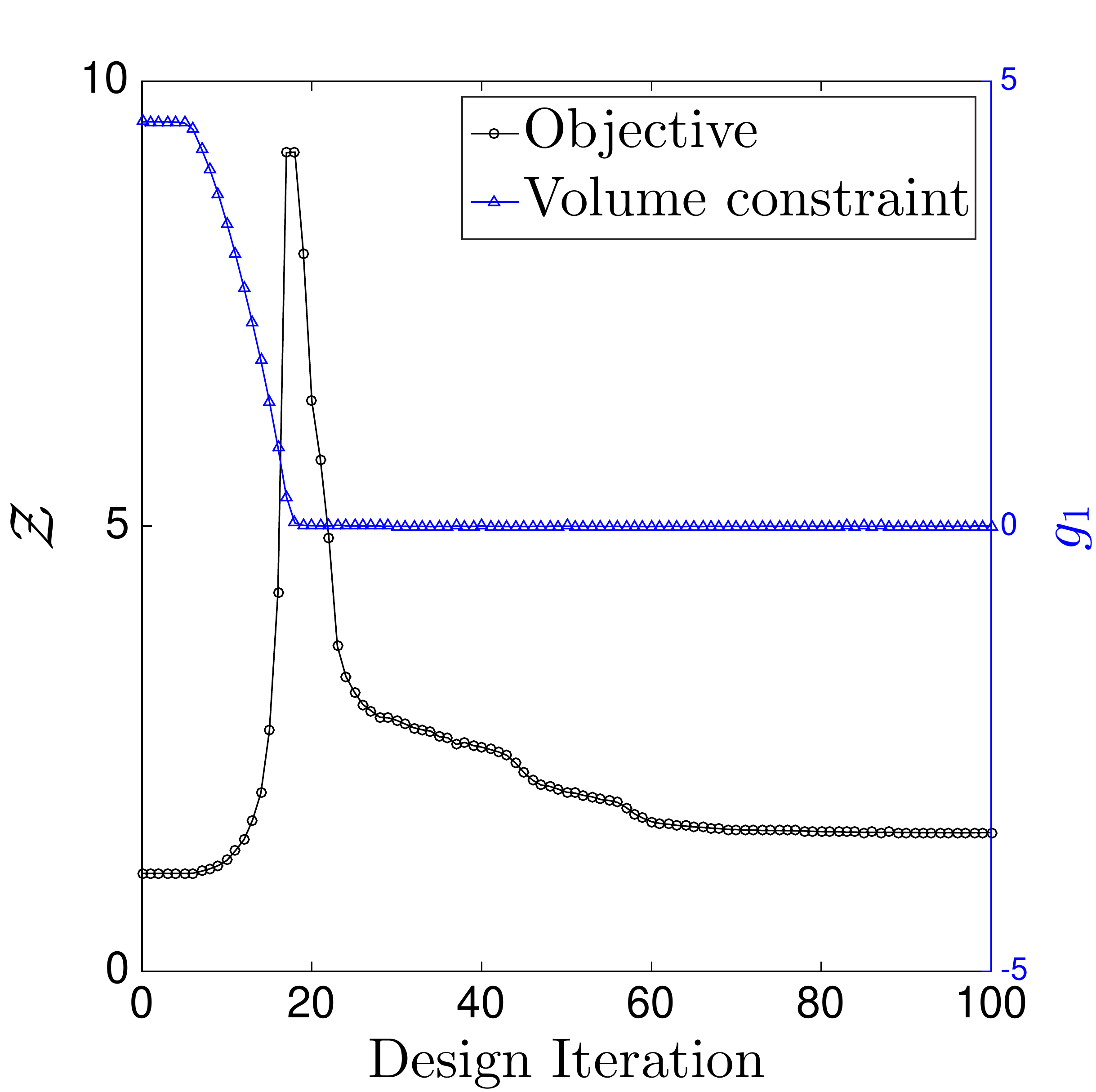}
		\label{fig:multiple_outlets_obj_and_const}
	}
	\subfloat[Lower bounds of mass inequality constraints.]{
		\includegraphics[width=0.32\linewidth]{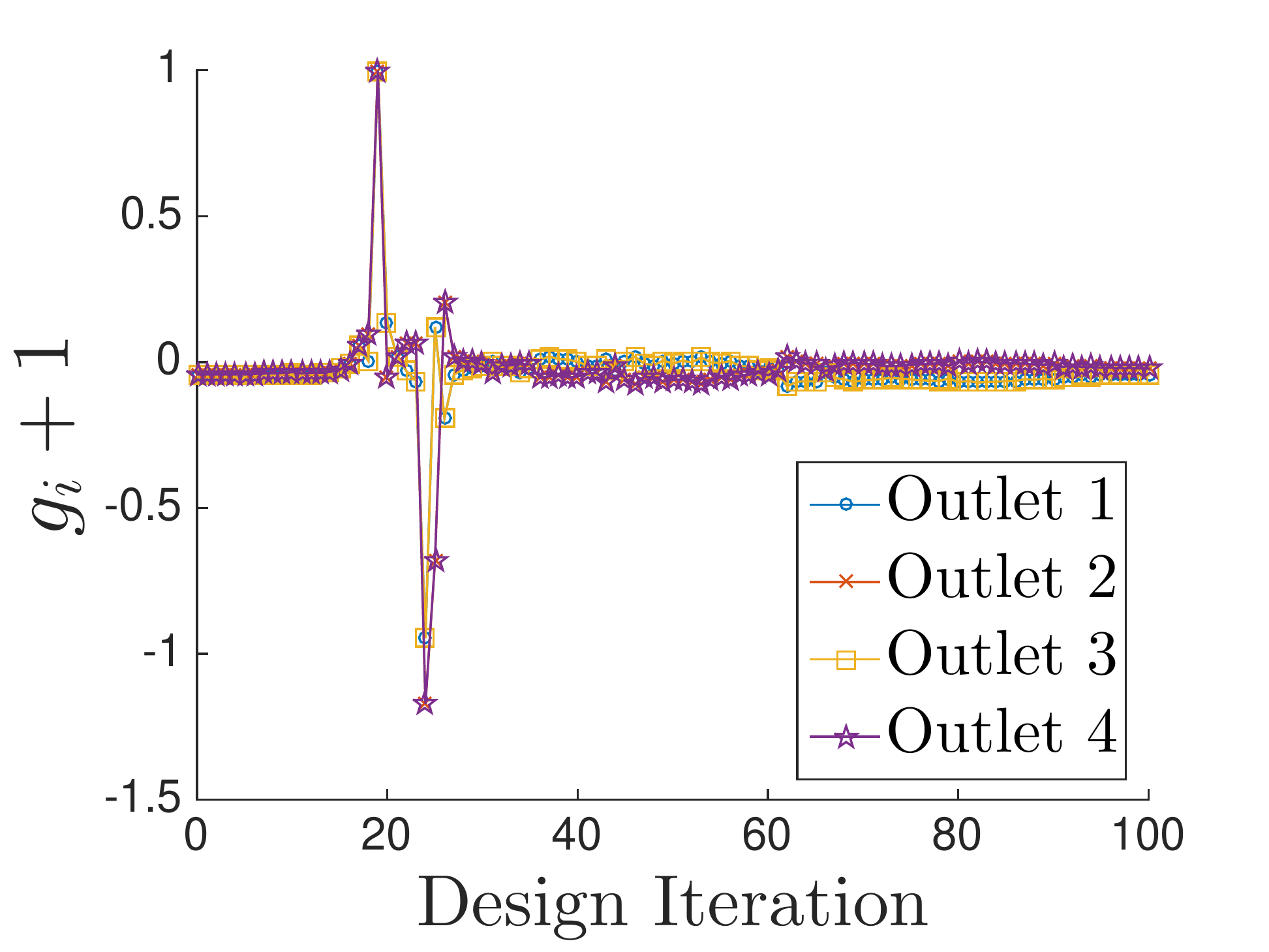}
		\label{fig:multiple_outlets_mass_lower}
	}
	\subfloat[Upper bounds of mass inequality constraints.]{
		\includegraphics[width=0.32\linewidth]{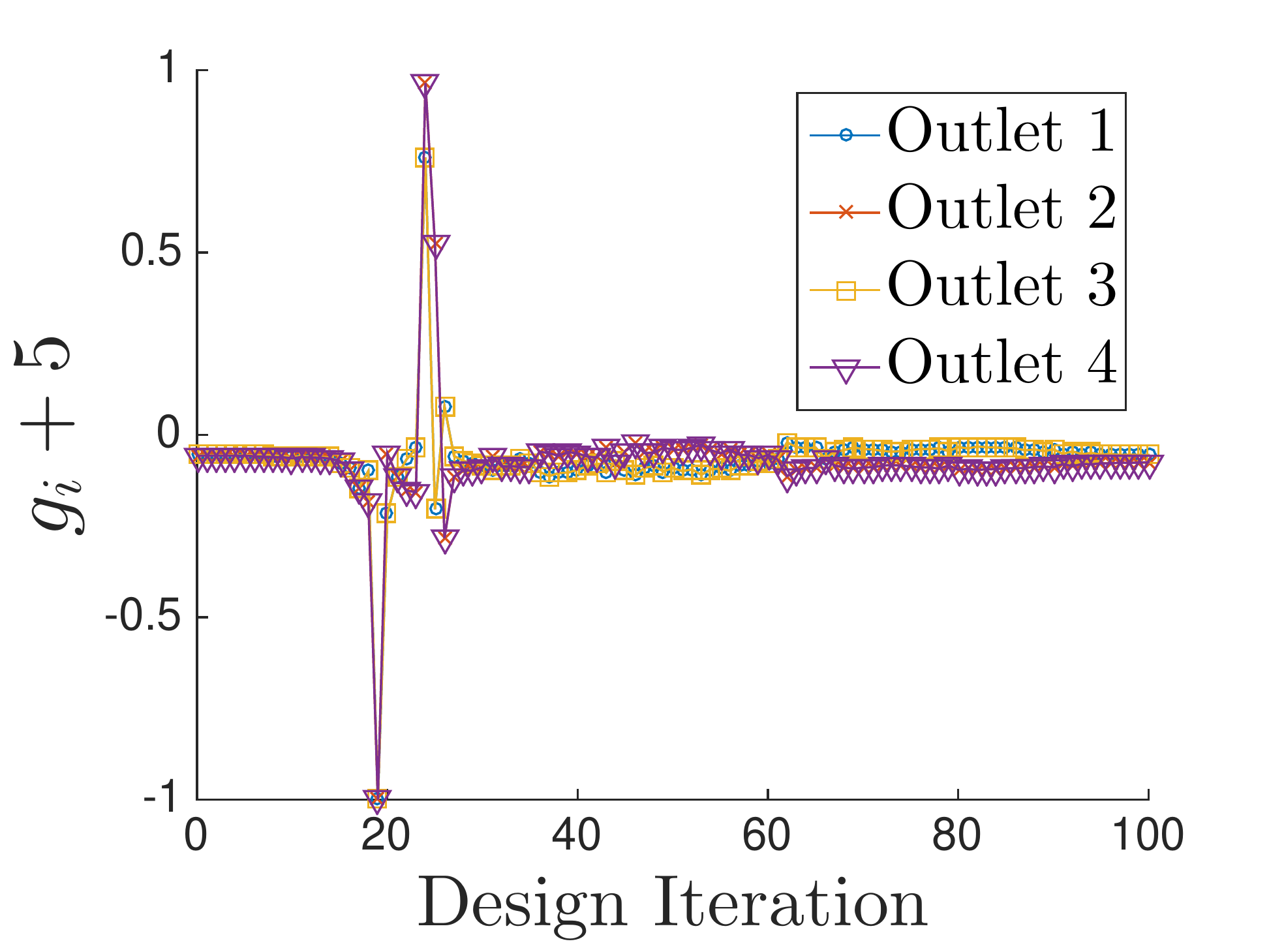}
		\label{fig:multiple_outlets_mass_upper}
	}
	\caption{Convergence plots of the objective and constraints for the multiple outlets example.}
	\label{fig:multiple_outlets_plots}
\end{figure}

The converged design after $99$ iterations is shown in Figure \ref{fig:multiple_outlets_final_design}, which resembles conceptually the 2D results from \cite{PWE+:10}. Figure \ref{fig:multiple_outlets_obj_and_const} shows the convergence plots of the objective and the volume constraint. The initialization process described above leads to an initial design that violates the volume constraint. Initially, the objective increases while the volume constraint value decreases. Once the constraint is satisfied, the objective decreases until a feasible minimum is found. The objective changes from a normalized value of $1.1$, where the mass constraint was violated, to $1.56$. Figures \ref{fig:multiple_outlets_mass_lower} and \ref{fig:multiple_outlets_mass_upper} show the convergence plots for the lower and upper bounds of the mass inequality constraints. We can observe that the mass flow rate constraints are also satisfied, and the amount of fluid flow exiting through each outlet is virtually the same. The mass flow rate at the inlets is $6.24$, while the rates at each of the outlets is $1.56$.

\subsection{Design of a Manifold with Variable Outlets}
\label{sec:variable_outlets}
\begin{figure*}
	\centering
	\includegraphics[width=0.5\linewidth]{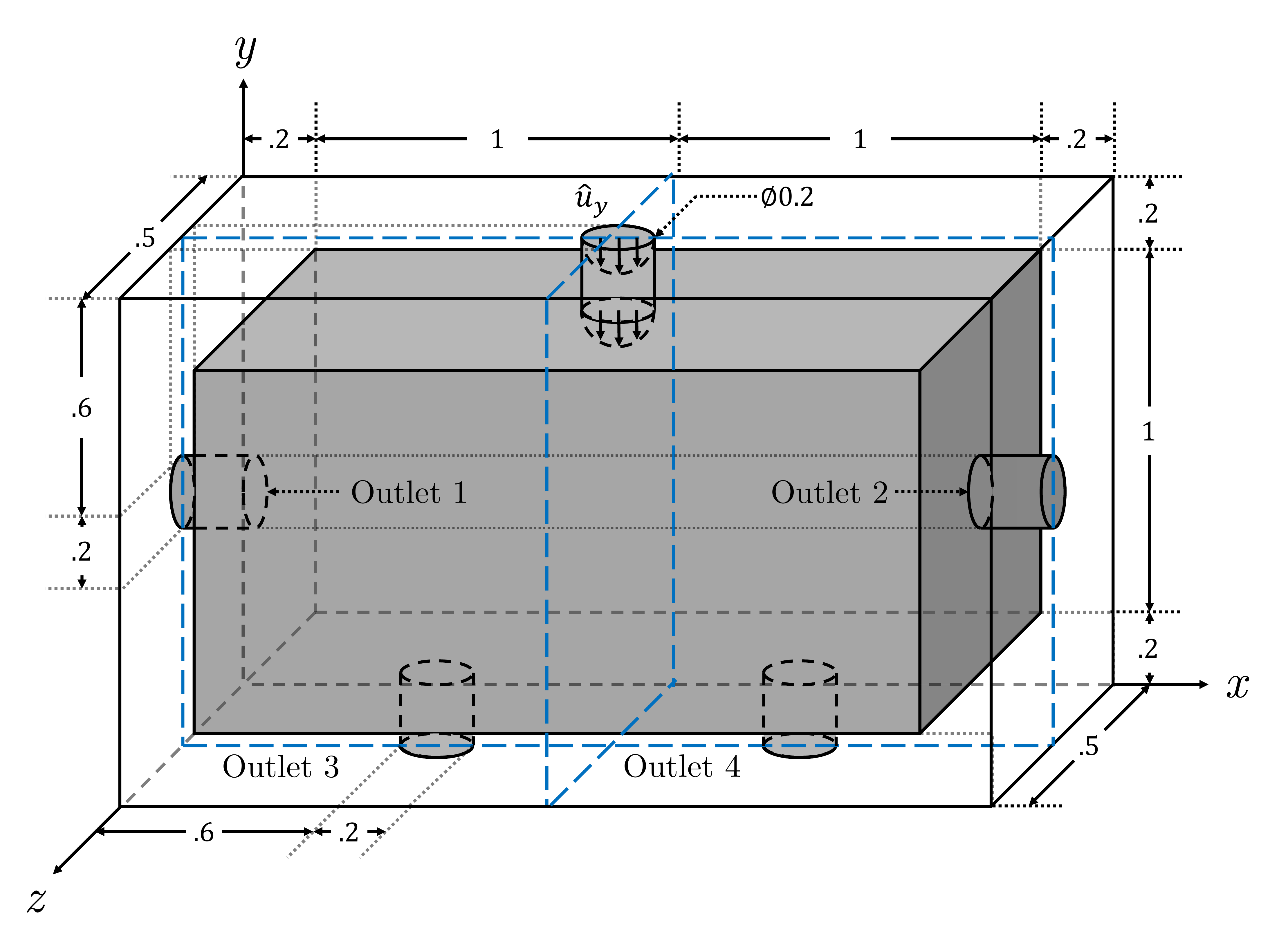}
	\caption{Problem setup for the variable outlets example. Blue dashes lines denote the symmetry planes.}
	\label{fig:setup_manifold}
\end{figure*}

In this example, we study a fluid flow problem with one inlet and multiple outlets, similar to the one from Example \ref{sec:multiple_outlets}, with the caveat that the outlets are allowed to vary in position and shape. The problem setup is shown in Figure \ref{fig:setup_manifold}. Fluid flows into the domain through the top inlet, and exits through the $2$ outlets on the left and right, and the $2$ outlets on the bottom. The outlets are described by cylinders using the level set parametrization of \eqref{eq:level_set_cylinder}. We introduce three additional design variables per outlet that define the in-plane coordinates at the center of the outlets and their radii; the lower and upper bounds are given in Table \ref{tab:multiple_outlets_bounds}. The interior of the domain is parametrized by \eqref{eq:smoothing_filter}, and the lower and upper bounds of the design variables are the same as in Table \ref{tab:multiple_outlets_params}. The inflow, outflow, and the interface boundary conditions are the same as in Example \ref{sec:multiple_outlets}. The problem parameters are given in Table \ref{tab:variable_outlets_params}. We model a quarter of the domain, and impose symmetry boundary conditions along the $x=1.2$ and the $z=0.25$ planes.

\begin{table}
	\centering
	\begin{tabular*}{0.75\linewidth}{@{\extracolsep{\fill}} l l}
	\hline
	                                    & Value \\
  \hline
  Mesh size                           & $56 \times 48 \times 10$ (quarter domain) \\
	Element size                        & $h = 0.025$ \\
	Characteristic velocity             & $u_{c} = 200$ \\
	Characteristic length               & $L_{c} = 0.2$ \\
	Dynamic viscosity                   & $\mu = 1$ \\
	Density                             & $\rho = 1$ \\
  Nitsche velocity penalty            & $\alpha_{\mathrm{N}, \bm{u}} = 100$ \\
  Nitsche indicator field penalty     & $\alpha_{\mathrm{N}, \psi} = 1$ \\
  Viscous ghost-penalty               & $\alpha_{\mathrm{GP}, \mu} = 0.5$ \\
  Pressure ghost-penalty              & $\alpha_{\mathrm{GP}, p} = 0.05$ \\
  Convective ghost-penalty            & $\alpha_{\mathrm{GP}, \bm{u}} = 0.5$ \\
	Pressure constraint parameter       & $k_{p} = 1$ \\
	Surface area scaling weight         & $w_{\mathcal{S}} = 0.01$ \\
	Volume constraint                   & $5\%$ \\
	Number of design variables          & $18,491$ (quarter domain) \\
	Design variables bounds             & $s^{L}_{i} = -0.0125$, $s^{U}_{i} = +0.0125$ \\
	Smoothing filter radius             & $r_{\phi} = 2.4h$ \\
  \hline
	\end{tabular*}
	\caption{Problem parameters for the multiple outlets example.}
	\label{tab:variable_outlets_params}
\end{table}

We seek to minimize the total pressure drop between the inlet and outlet ports and the surface area of the fluid-solid interface, and impose a $20\%$ fluid volume constraint \eqref{eq:multiple_outlets_vol_const}. Further, we formulate constraints on the mass flow rates such that the left and right outlets have $33.333\%$ each, and the bottom outlets have $16.667\%$ each, of the fluid flow entering through the inlet. Similar to \eqref{eq:multiple_outlets_mass_const}, we formulate these limits as inequality constraints and set the tolerances to $\Delta\dot{m}^{1,2}=\mp0.333\%$ for the left and right outlets and to $\Delta\dot{m}^{3,4}=\mp0.167\%$ for the bottom outlets. As we assume symmetry and model only one quarter of the design domain, constraints for only one outlet on a vertical face, i.e.~Outlet 1 or 2, and one outlet at the bottom face, i.e.~Outlet 3 or 4, need to be imposed; see Figure \ref{fig:setup_manifold}. This leads to the following four additional constraints:
%
\begin{equation}
	\label{eq:variable_outlets_mass_const}
	\begin{split}
	g_{2} = 1 - \frac{\dot{m}_{\mathrm{out}, 1,2}}{\left(33.333\% - \Delta\dot{m}^{1,2}\right)\dot{m}_{\mathrm{in}}}\ , \qquad 
	g_{3} = 1 - \frac{\dot{m}_{\mathrm{out}, 3,4}}{\left(16.667\% - \Delta\dot{m}^{3,4}\right)\dot{m}_{\mathrm{in}}}\ , \\
	g_{4} = \frac{\dot{m}_{\mathrm{out}, 1,2}}{\left(33.333\% + \Delta\dot{m}^{1,2}\right)\dot{m}_{\mathrm{in}}} - 1\ , \qquad 
	g_{5} = \frac{\dot{m}_{\mathrm{out}, 3,4}}{\left(16.667\% + \Delta\dot{m}^{3,4}\right)\dot{m}_{\mathrm{in}}} - 1\ .
	\end{split}
\end{equation}
%
We impose these constraints via a continuation method and gradually contract the tolerances from initially $\Delta\dot{m}^{1,2}=\mp1.667\%$ and $\Delta\dot{m}^{3,4}=\mp0.833\%$ to their target values in the course of the optimization process.

\begin{table}
	\centering
	\def\arraystretch{1.5}
	\begin{tabular*}{0.625\linewidth}{@{\extracolsep{\fill}} l l l l l l l l l}
	\hline
	         & $x^{L}_{c}$ & $x^{U}_{c}$ & $y^{L}_{c}$ & $y^{U}_{c}$ & $z^{L}_{c}$ & $z^{U}_{c}$ & $r^{L}_{c}$ & $r^{U}_{c}$ \\
  \hline
  Outlet 1 & $-$         & $-$         & $0.2$         & $1.2$        & $0.2$         & $1.2$        & $0$ & $0.5$ \\
	Outlet 2 & $-$         & $-$         & $0.2$         & $1.2$        & $0.2$         & $1.2$        & $0$ & $0.5$ \\
	Outlet 3 & $0.2$         & $2.2$       & $-$         & $-$         & $0.2$         & $1.2$        & $0$ & $0.5$ \\
	Outlet 4 & $0.2$         & $2.2$       & $-$         & $-$         & $0.2$         & $1.2$        & $0$ & $0.5$ \\
    \hline
	\end{tabular*}
	\caption{Lower and upper bounds of the design variables that control the position and size of the ports in the example with variables outlets.}
	\label{tab:multiple_outlets_bounds}
\end{table}

\begin{figure}
	\centering
	\includegraphics[width=0.45\linewidth]{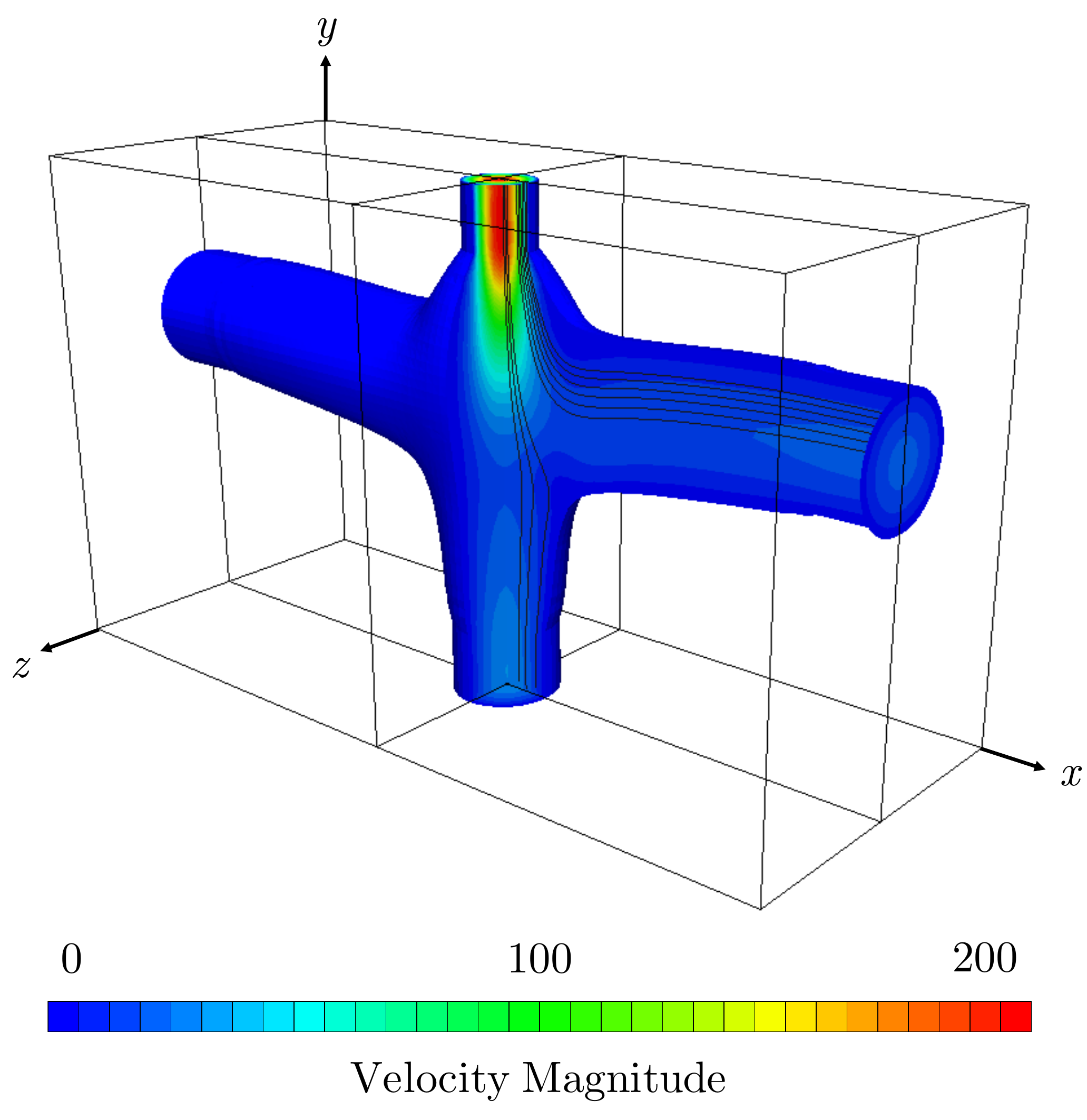}
	\caption{Velocity magnitude, with streamlines, of the optimized material layout for the variable outlets example. A section of the design was removed for visualization purposes.}
	\label{fig:manifold_final_design}
\end{figure}

The design domain is initialized with $3 \times 6 \times 1$ spherical solid inclusions of radii $0.125$. The optimized design after $350$ iterations is shown in Figure \ref{fig:manifold_final_design}. The mass flow through each outlet, and the radii and coordinates of the outlets are given in Table \ref{tab:manifold_outlets_solutions}. During the optimization process, the left and right outlets increase their radii and move up in the $y$-direction. The bottom outlets merge and form a single outlet. The numbers in Table \ref{tab:manifold_outlets_solutions} reveal that these outlets do not fully overlap, and that they have moved away from the $x = 1.2$ and $z=0.25$ planes of symmetry. The relative difference between their final positions and the axes of symmetry is small nonetheless, in the order of $0.1\%$. This value got increasingly closer to $0$ as the mass flow rate tolerances of \eqref{eq:variable_outlets_mass_const} were lowered. The mass flow exiting through each of the left and right outlets is $33.2\%$ of the mass flow entering through the inlet, and the mass flow exiting through each of the bottom outlets is $16.8\%$; these values are within the lower and upper bounds of the mass flow limits.

\begin{table*}
	\centering
	\def\arraystretch{1.5}
	\begin{tabular*}{0.75\linewidth}{@{\extracolsep{\fill}} l l l l l l l}
	\hline
	         & $\dot{m}$ & Fraction & Radius   & $x$ & $y$ & $z$ \\
    \hline
	Inlet    & $3.135$   & $100\%$   & $0.1$    & $1.2$    & $1.4$   & $0.7$ \\
	Outlet 1 & $1.040$   & $33.2\%$  & $0.171$  & $2.4$    & $0.897$ & $0.7114$ \\
	Outlet 2 & $1.040$   & $33.2\%$  & $0.171$  & $0$      & $0.897$ & $0.7114$ \\
	Outlet 3 & $0.528$   & $16.8\%$  & $0.1482$ & $1.2048$ & $0$     & $0.7092$ \\
	Outlet 4 & $0.528$   & $16.8\%$  & $0.1482$ & $1.1952$ & $0$     & $0.7092$ \\
    \hline
	\end{tabular*}
	\caption{Mass flow rate, fraction of the mass flow entering through the inlet, radii, and coordinates of the outlets in the variable outlets example.}
	\label{tab:manifold_outlets_solutions}
\end{table*}


\subsection{Design of a Micromixer}
\label{sec:micromixer}

\begin{figure*}
	\centering
	\includegraphics[width=0.65\linewidth]{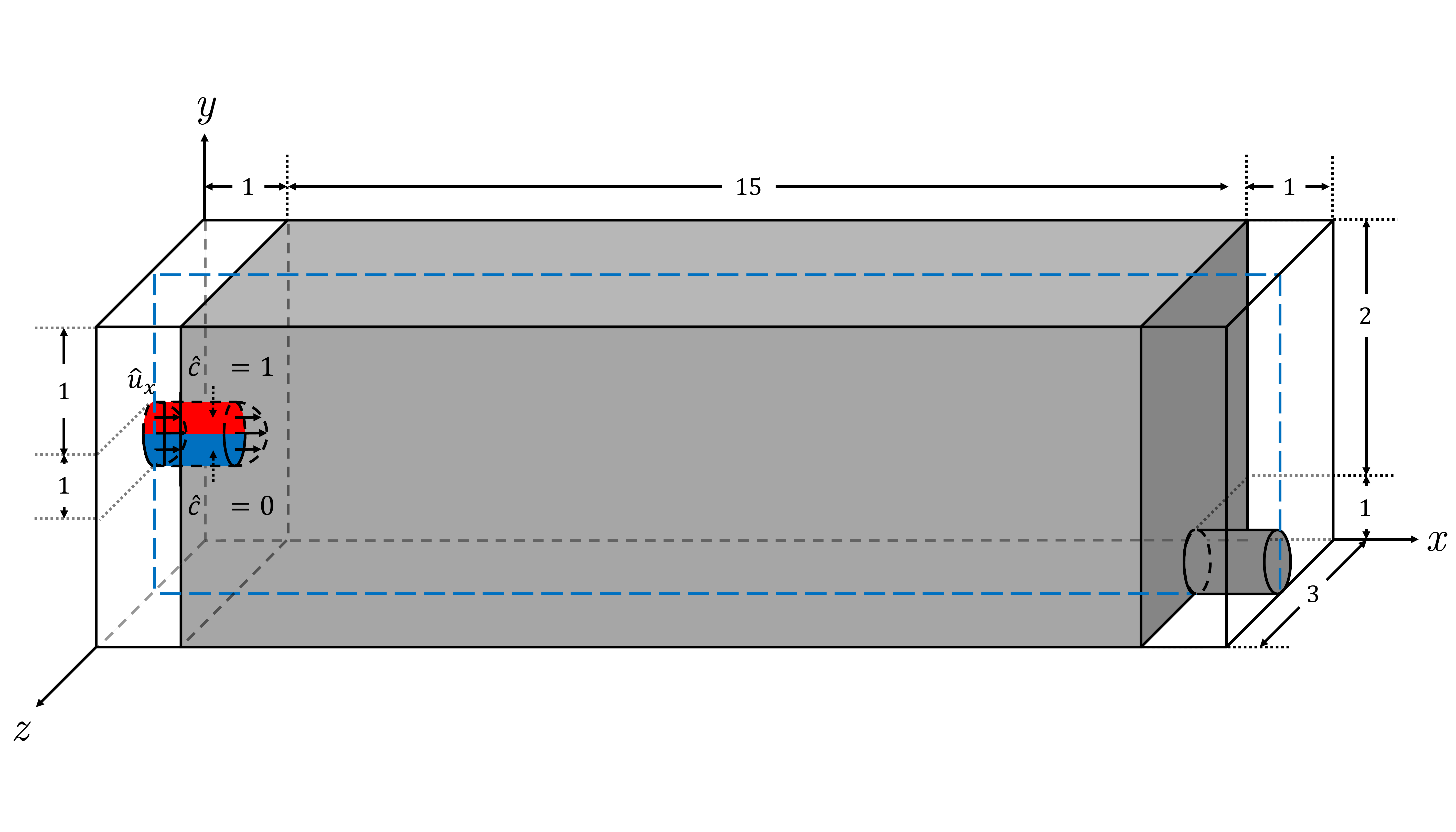}
	\caption{Problem setup for the micromixer example. Blue dashes lines denote the symmetry plane.}
	\label{fig:setup_micromixer}
\end{figure*}

In this example, we apply our CutFEM framework to the modeling and optimizing of a micromixer at a low Reynolds number and steady-state conditions. The example is the 3D analog to the 2D problem found in \cite{MM:14}, and is similar to the micromixer studies with flow topology optimization found in \cite{AGS:09}, \cite{LDZ+:13} and \cite{LR:16}. The problem setup is shown in Figure \ref{fig:setup_micromixer}. A ``red'' fluid and a ``blue'' fluid enter the design domain through the left inlet and exit it through the lower right side. We assume that the fluids are ideally miscible and have identical flow properties. The ``red'' fluid is represented by a species concentration value of $\hat{c} = 1$, and the ``blue'' fluid by $\hat{c} = 0$. We do not consider diffusion of the species field through the solid phase. The inflow condition is formulated using the same approach as in Example \ref{sec:pipebend}. No-slip boundary conditions are applied at the fluid-solid interface, and a traction-free boundary condition is imposed on the outlet. We only model one half of the domain, and symmetry boundary conditions are imposed on the $z=1.5$ plane. An adiabatic condition is imposed on the fluid-solid interface for the species field.

The objective uses the target scalar value formulation from \eqref{eq:target_scalar_value_criteria} at steady-state, and is defined as:
%
\begin{equation}
	\label{eq:micromixer_objective}
	\mathcal{Z} = \frac{\mathcal{K}_{\mathrm{out}}}{\left\Vert\mathcal{K}^{0}_{\mathrm{out}}\right\Vert} + w_{\mathcal{S}} \frac{\mathcal{S}}{\left\Vert\mathcal{S}^{0}\right\Vert} \ ,
\end{equation}
%
with $\beta_{KS} = 400$, and $c_{\mathrm{ref}} = 0.5$. A small surface area penalty with $w_{\mathcal{S}} = 0.001$  is applied to regularize the optimization problem. The design is subject to a volume constraint \eqref{eq:multiple_outlets_vol_const} of $35\%$ to suppress trivial solutions, and to promote the formation of distinct fluid channels. Similar to \cite{MM:14}, a constraint is imposed on the maximum pressure drop to prevent the formation of small geometric features:
%
\begin{equation}
	\label{eq:micromixer_pressure_constraint}
	g_{2} = \frac{\mathcal{T}_{\mathrm{in}} - \mathcal{T}_{\mathrm{out}}}{\mathrm{\Delta} p_{\mathrm{ref}}} - 1 \ .
\end{equation}
%
The problem is initialized with $15 \times 5 \times 5$ spherical solid inclusions of radii $0.1$, similar to the previous examples. The remaining parameters are given in Table \ref{tab:micromixer_params}.

\begin{table}
	\centering
	\begin{tabular*}{0.75\linewidth}{@{\extracolsep{\fill}} l l}
	\hline
	                                    & Value \\
  \hline
  Mesh size                           & $226 \times 40 \times 20$ (half domain) \\
	Element size                        & $h = 0.075$ \\
	Characteristic velocity             & $u_{c} = 1$ \\
	Characteristic length               & $L_{c} = 1$ \\
	Dynamic viscosity                   & $\mu = 1$ \\
	Density                             & $\rho = 1$ \\
	Specific heat capacity              & $c_{p} = 1$ \\
	Diffusivity                         & $k = 0.001$ \\
  Nitsche velocity penalty            & $\alpha_{\mathrm{N}, \bm{u}} = 100$ \\
  Nitsche species penalty             & $\alpha_{\mathrm{N}, c} = 1$ \\
  Nitsche indicator field penalty     & $\alpha_{\mathrm{N}, \psi} = 1$ \\
  Viscous ghost-penalty               & $\alpha_{\mathrm{GP}, \mu} = 0.05$ \\
  Pressure ghost-penalty              & $\alpha_{\mathrm{GP}, p} = 0.005$ \\
  Convective ghost-penalty            & $\alpha_{\mathrm{GP}, \bm{u}} = 0.05$ \\
  Species ghost-penalty               & $\alpha_{\mathrm{GP}, c} = 0.05$ \\
	Pressure constraint parameter       & $k_{p} = 1$ \\
	Pressure constraint reference       & $\mathrm{\Delta} p_{\mathrm{ref}} = 30$ \\
	Surface area objective weight       & $w_{\mathcal{S}} = 0.001$ \\
	Volume constraint                   & $35\%$ \\
	Number of design variables          & $173,061$ (half domain) \\
	Design variables bounds             & $s^{L}_{i} = -0.0375L$, $s^{U}_{i} = +0.0375L$ \\
	Smoothing filter radius             & $2.4h$ \\
  \hline
	\end{tabular*}
	\caption{Problem parameters for the micromixer example.}
	\label{tab:micromixer_params}
\end{table}

\begin{figure*}
	\centering
	\subfloat[Species field.]{
		\includegraphics[width=0.4\linewidth]{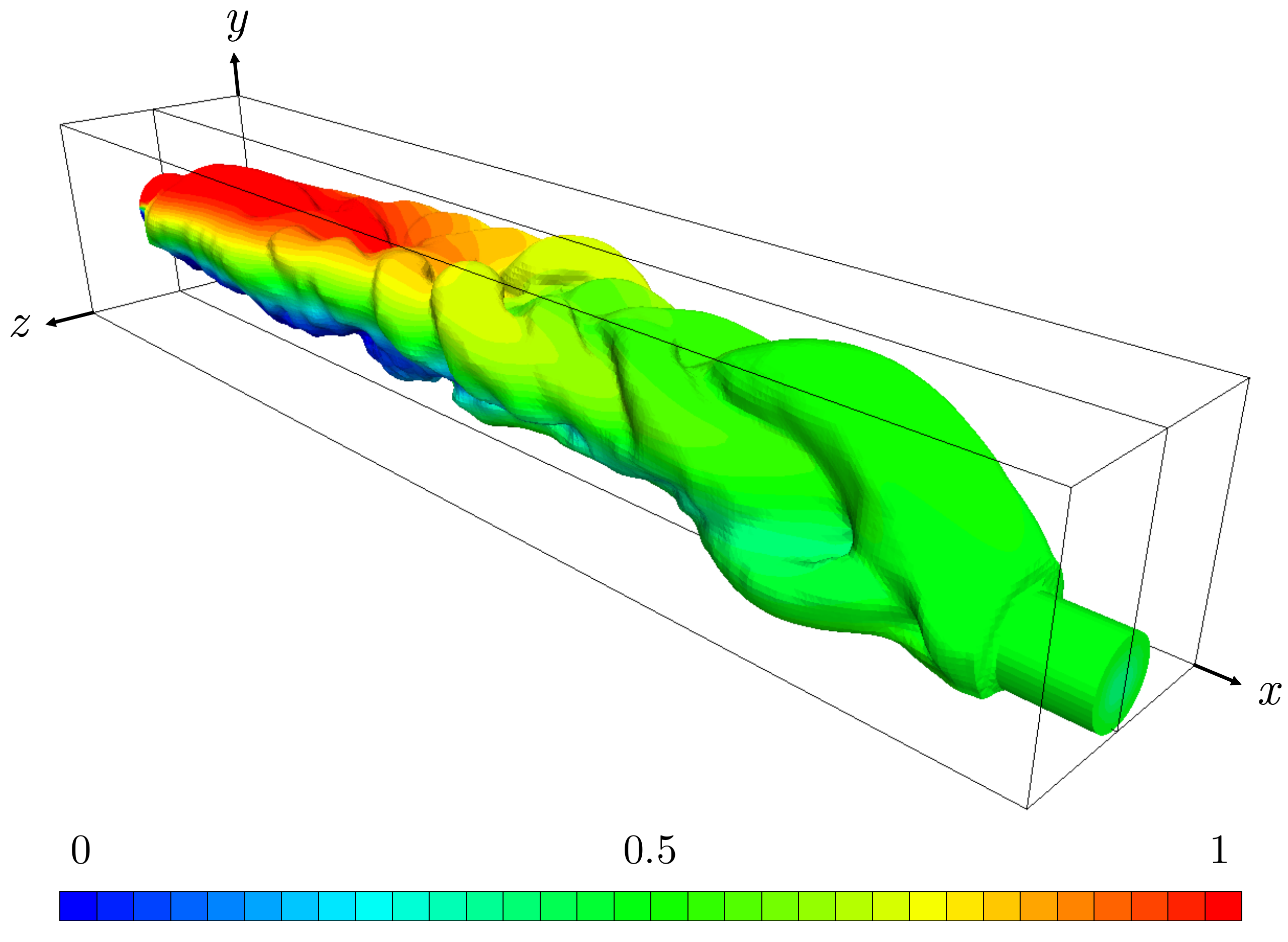}
		\label{fig:micromixer_final_design_2}
	} \qquad
	\subfloat[Velocity magnitude, with streamlines.]{
		\includegraphics[width=0.4\linewidth]{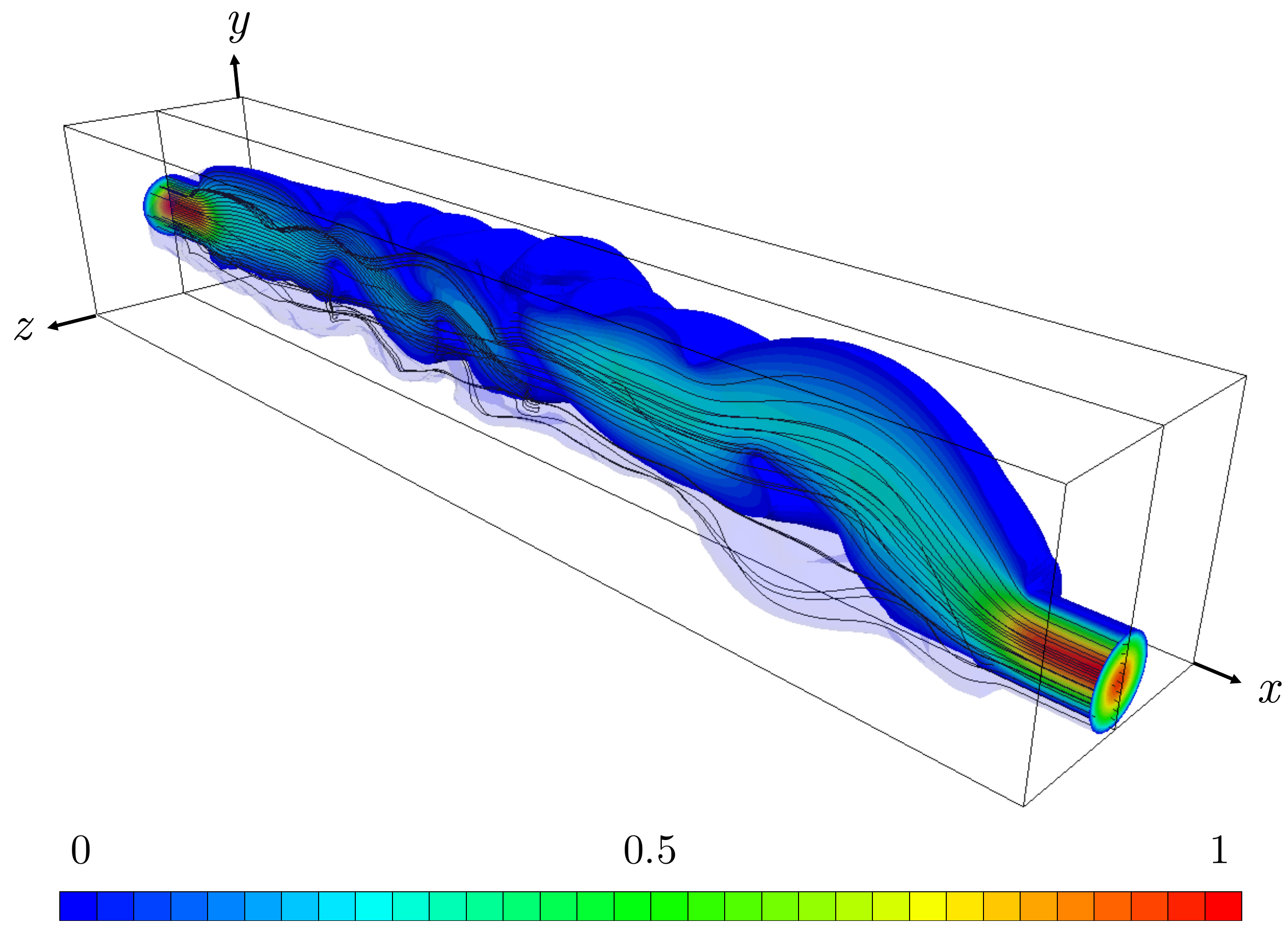}
		\label{fig:micromixer_final_design_1}
	} \\
	\subfloat[Species field, top view.]{
		\includegraphics[width=0.5\linewidth]{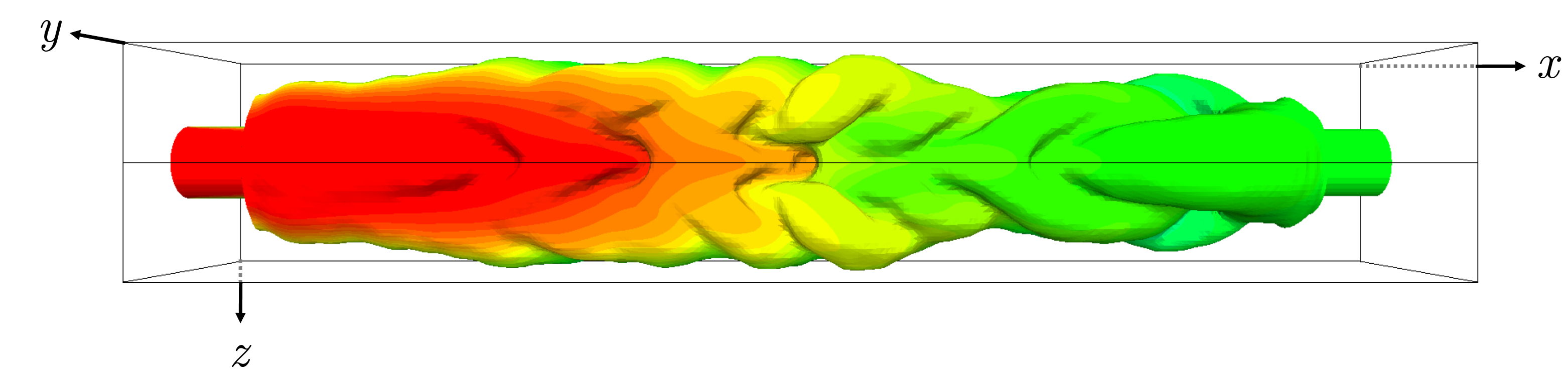}
		\label{fig:micromixer_final_design_3}
	}
	\subfloat[Species field, bottom view.]{
		\includegraphics[width=0.5\linewidth]{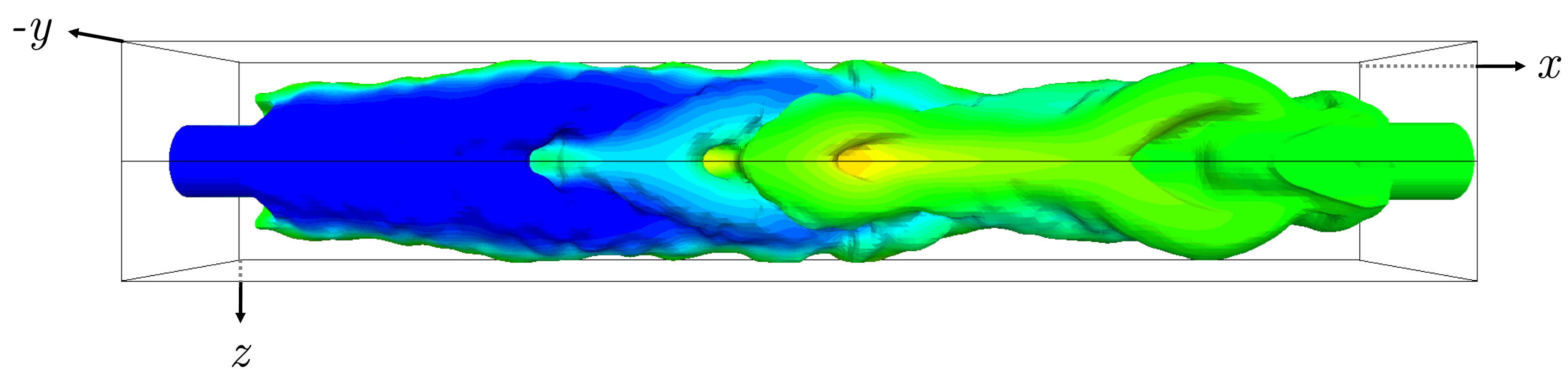}
		\label{fig:micromixer_final_design_4}
	}
	\caption{Optimized material layout for the micromixer example.}
	\label{fig:micromixer_final_design}
\end{figure*}

The optimized material layout, along with the species field and the velocity magnitude, is shown in Figure \ref{fig:micromixer_final_design}. The average species concentration at the outlet is $0.46$. Analogous to the 2D results from \cite{MM:14}, the length of the channel layout increases by creating an intricate wavy design, thereby increasing the path traveled by the fluids. This lengthening mechanism is the key to enhance the mixing of the fluids in 3D laminar flows. The channel does not contain internal isolated solid particles. The number of iterations is rather large: $1,500$ iterations are required to form the channel layout and fine-tune its shape. This behavior has also been seen by \cite{MM:14}, who attributes this to the interplay of localized sensitivities along the fluid-solid interface and the volume constraint on the fluid phase. In contrast to the density approach of \cite{MPY+:12}, we do not obtain numerical artifacts in the optimized design.


\subsection{Design of an Oscillating Pump}
\label{sec:pump}

\begin{figure}
	\centering
	\includegraphics[width=0.4\linewidth]{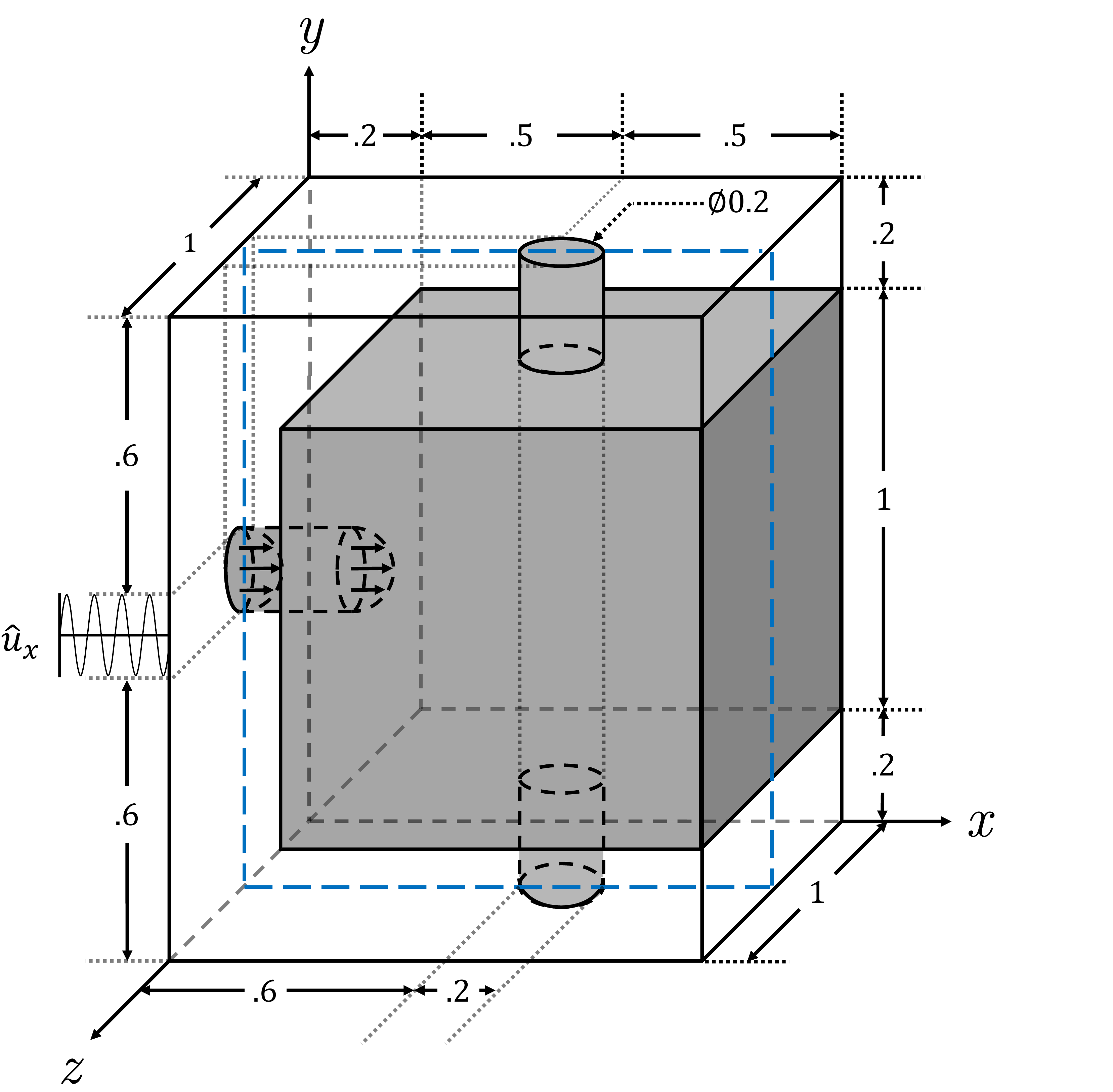}
	\caption{Problem setup for the transient pump example. Blue dashes lines denote the symmetry plane.}
	\label{fig:setup_pump}
\end{figure}

In the final example, we study the applicability of the CutFEM framework to problems with transient behavior. We consider the design of a simplified fluid pump. The example is the 3D analog to the 2D problem found in \cite{NSL:16}. The problem setup is shown in Figure \ref{fig:setup_pump}. The basic idea is to prescribe a harmonically oscillating velocity at the left port, and optimize the design domain to maximize the amount of fluid that is transported through the port at the top. The port at the bottom represents a reservoir from which additional fluid can enter the domain.

The objective is composed of the negative of the mass flow through the top outlet and the surface area of the fluid-solid interface:
%
\begin{equation}
	\label{eq:pump_objective}
	\mathcal{Z} = -\frac{\frac{1}{N_{t}}\sum\limits^{N_{t}}_{n=0}\dot{m}_{\mathrm{out}, 1}\left(t_{n}\right)}{\left\Vert\frac{1}{N_{t}}\sum\limits^{N_{t}}_{n=0}\dot{m}^{0}_{\mathrm{out}, 1}\left(t_{n}\right)\right\Vert} + w_{\mathcal{S}} \frac{\mathcal{S}}{\left\Vert\mathcal{S}^{0}\right\Vert}\ ,
\end{equation}
%
where $N_{t}$ is the total number of time iterations. The design domain is subject to a volume constraint \eqref{eq:multiple_outlets_vol_const} of $15\%$ to suppress trivial solutions and to promote the formation of smooth fluid channels.

We apply traction-free boundary conditions at the top and bottom ports. The flow velocity at the left port is prescribed as:
%
\begin{equation}
	\label{eq:pump_inflow}
	\hat{u}_{x}\left(0, y, z, t\right) = u_{c} \cdot \sin\left({v t}\right) \cdot \left(
	\left(-\frac{4}{L_{c}^{2}}\right) \cdot
	\left({\left(y - y_{c}\right)}^{2} + {\left(z - z_{c}\right)}^{2}\right)
	 + 1\right)\ , \qquad \hat{u}_{y} = \hat{u}_{z} = 0\ ,
\end{equation}
%
where the inflow frequency is set to $v=100\pi$. No-slip boundary conditions are imposed on the fluid-solid interface. We only model one half of the domain, and symmetry boundary conditions are imposed on the plane $z=0.5$. We model a single pumping cycle with $40$ time steps per cycle. Numerical experiments showed that this number of time steps is sufficient to model the harmonic behavior of the flow. The design is initialized with the same layout as in Example \ref{sec:multiple_outlets}. The remaining problem parameters are given in Table \ref{tab:pump_params}. 

\begin{table}
	\centering
	\begin{tabular*}{0.75\linewidth}{@{\extracolsep{\fill}} l l}
	\hline
	                                    & Value \\
  \hline
  Mesh size                           & $72 \times 84 \times 30$ (half domain) \\
	Element size                        & $h = 0.0167$ \\
	Characteristic velocity             & $u_{c} = 200$ \\
	Characteristic length               & $L_{c} = 0.2$ \\
	Dynamic viscosity                   & $\mu = 1$ \\
	Density                             & $\rho = 1$ \\
	Inlet velocity frequency            & $v = 100\pi$ \\
	Time step                           & $\mathrm{\Delta} t = 5 \times 10^{-4}$ \\
	Number of time iterations           & $N_{t} = 2 \times 10^{-2}$ \\
  Nitsche velocity penalty            & $\alpha_{\mathrm{N}, \bm{u}} = 100$ \\
  Nitsche indicator field penalty     & $\alpha_{\mathrm{N}, \psi} = 1$ \\
  Viscous ghost-penalty               & $\alpha_{\mathrm{GP}, \mu} = 0.05$ \\
  Pressure ghost-penalty              & $\alpha_{\mathrm{GP}, p} = 0.005$ \\
  Convective ghost-penalty            & $\alpha_{\mathrm{GP}, \bm{u}} = 0.05$ \\
	Pressure constraint parameter       & $k_{p} = 1$ \\
	Surface area objective weight       & $w_{\mathcal{S}} = 0.1$ \\
	Volume constraint                   & $15\%$ \\
	Number of design variables          & $132,928$ (half domain) \\
	Design variables bounds             & $s^{L}_{i} = -0.0083$, $s^{U}_{i} = +0.0083$ \\
	Smoothing filter radius             & $2.4h$ \\
  \hline
	\end{tabular*}
	\caption{Problem parameters for the transient pump example.}
	\label{tab:pump_params}
\end{table}

\begin{figure*}
	\centering
	\subfloat[Optimized layout.]{
		\includegraphics[width=0.3\linewidth]{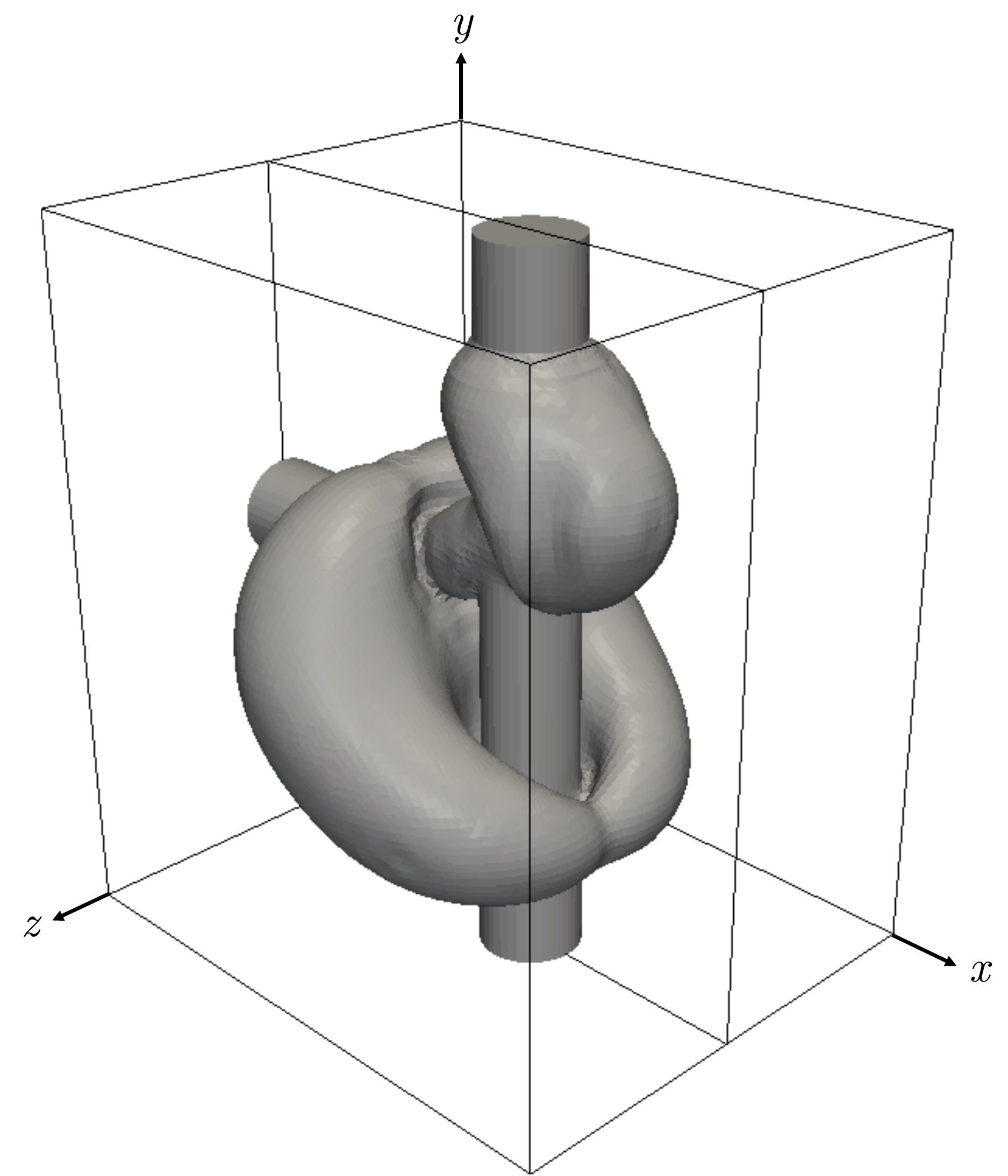}
		\label{fig:pump_final_design_1}
	}
	\subfloat[$xy$ plane view.]{
		\includegraphics[width=0.3\linewidth]{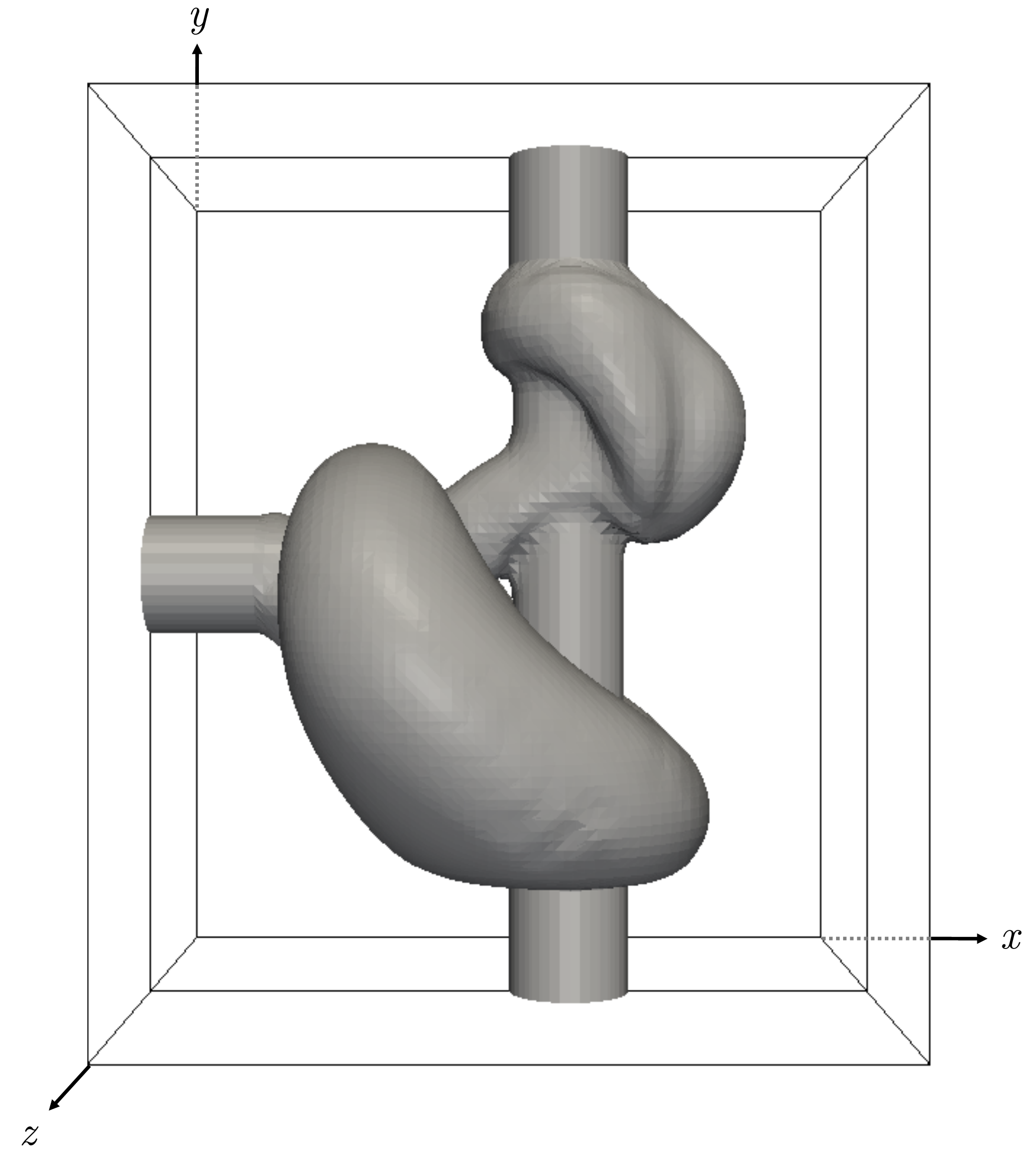}
		\label{fig:pump_final_design_2}
	}
	\subfloat[$yz$ plane view.]{
		\includegraphics[width=0.3\linewidth]{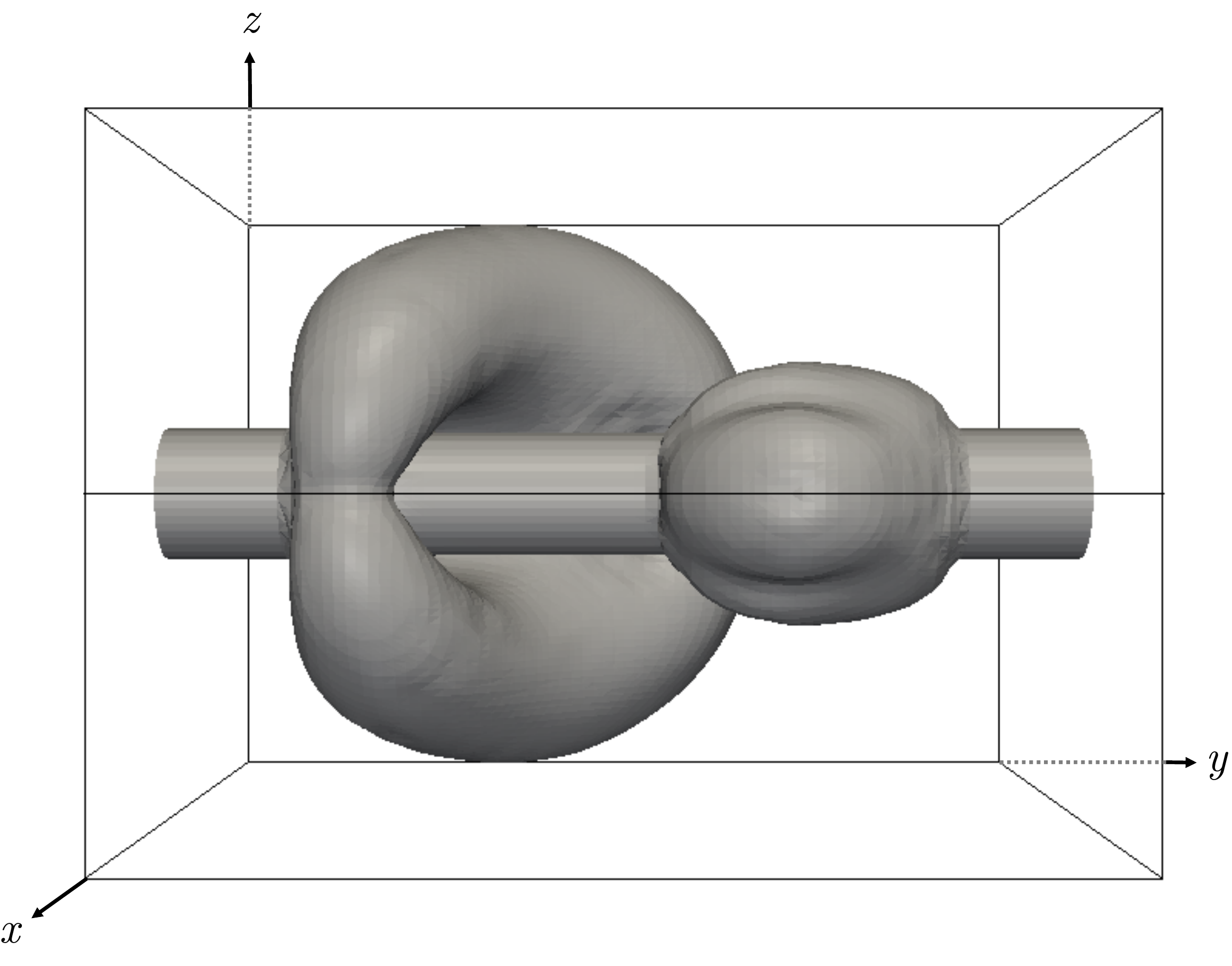}
		\label{fig:pump_final_design_3}
	}
	\caption{Optimized material layout for the transient pump example.}
	\label{fig:pump_final_design}
\end{figure*}

\begin{figure*}
	\centering
	\subfloat[Inflow cycle.]{
		\includegraphics[width=0.35\linewidth]{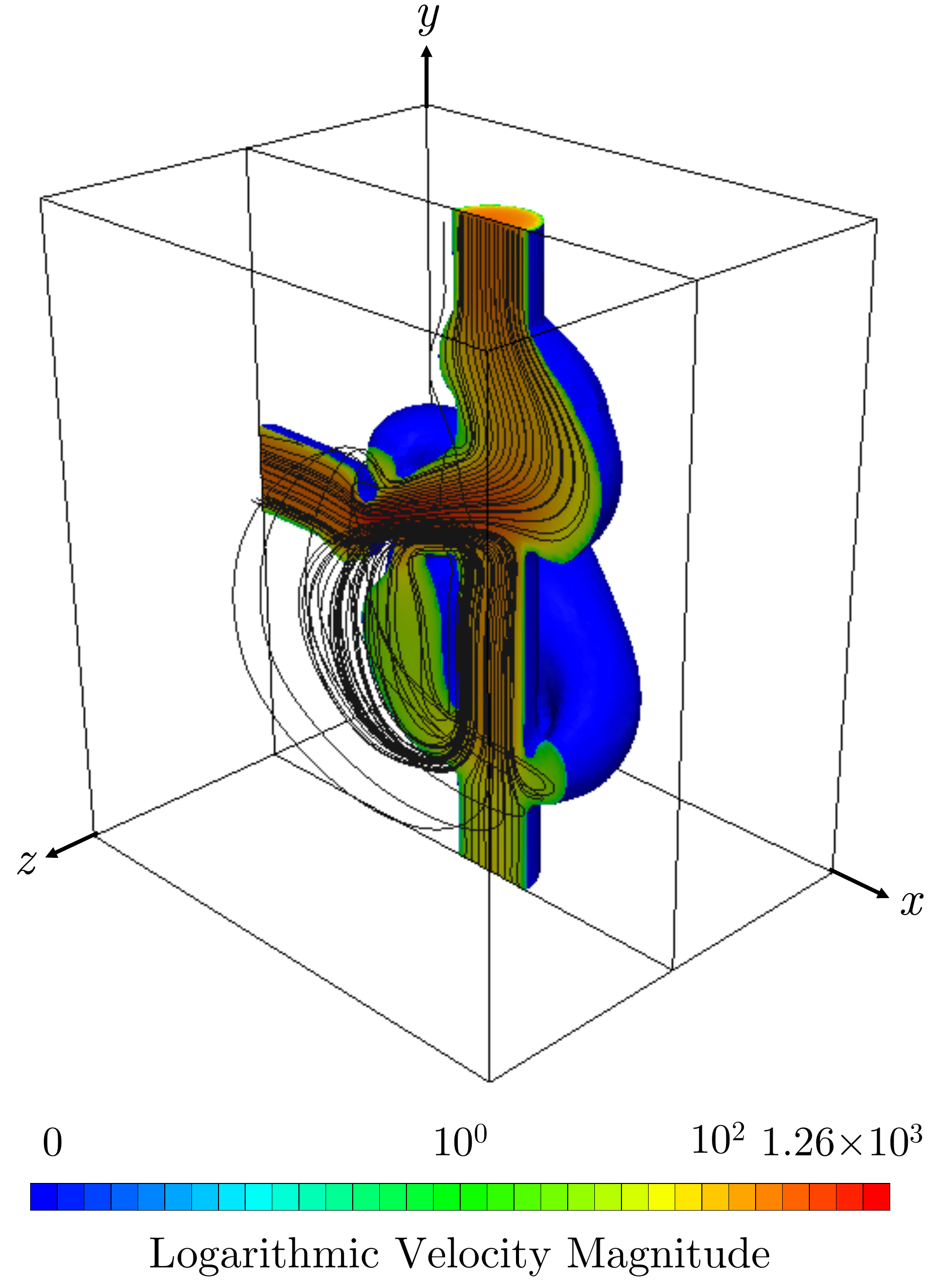}
		\label{fig:pump_final_design_stream_in}
	} \qquad
	\subfloat[Outflow cycle.]{
		\includegraphics[width=0.35\linewidth]{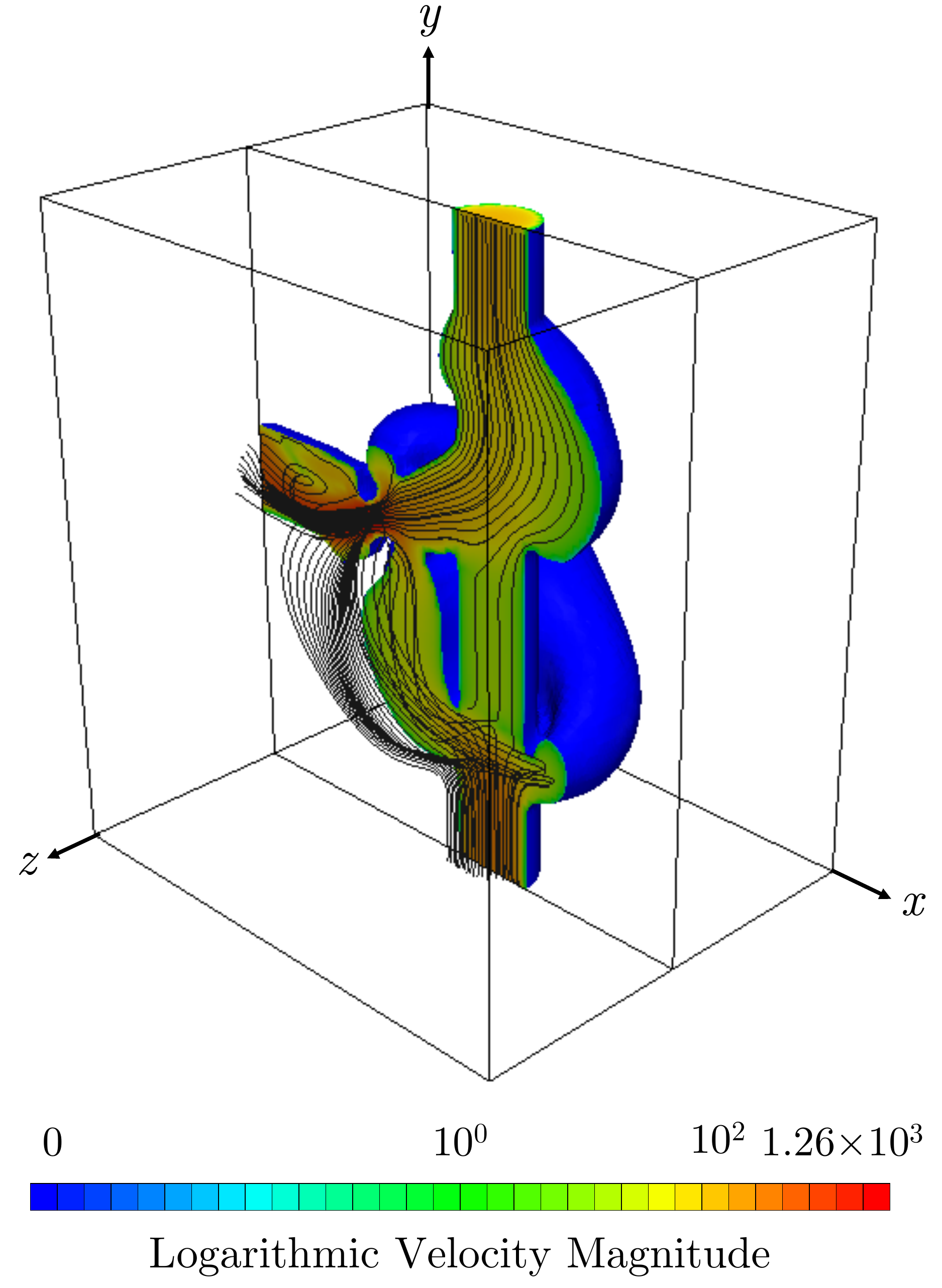}
		\label{fig:pump_final_design_stream_out}
	}
	\caption{Velocity streamlines for half the domain of the transient pump example.}
	\label{fig:pump_final_design_stream}
\end{figure*}

\begin{figure*}
	\centering
	\includegraphics[width=0.45\linewidth]{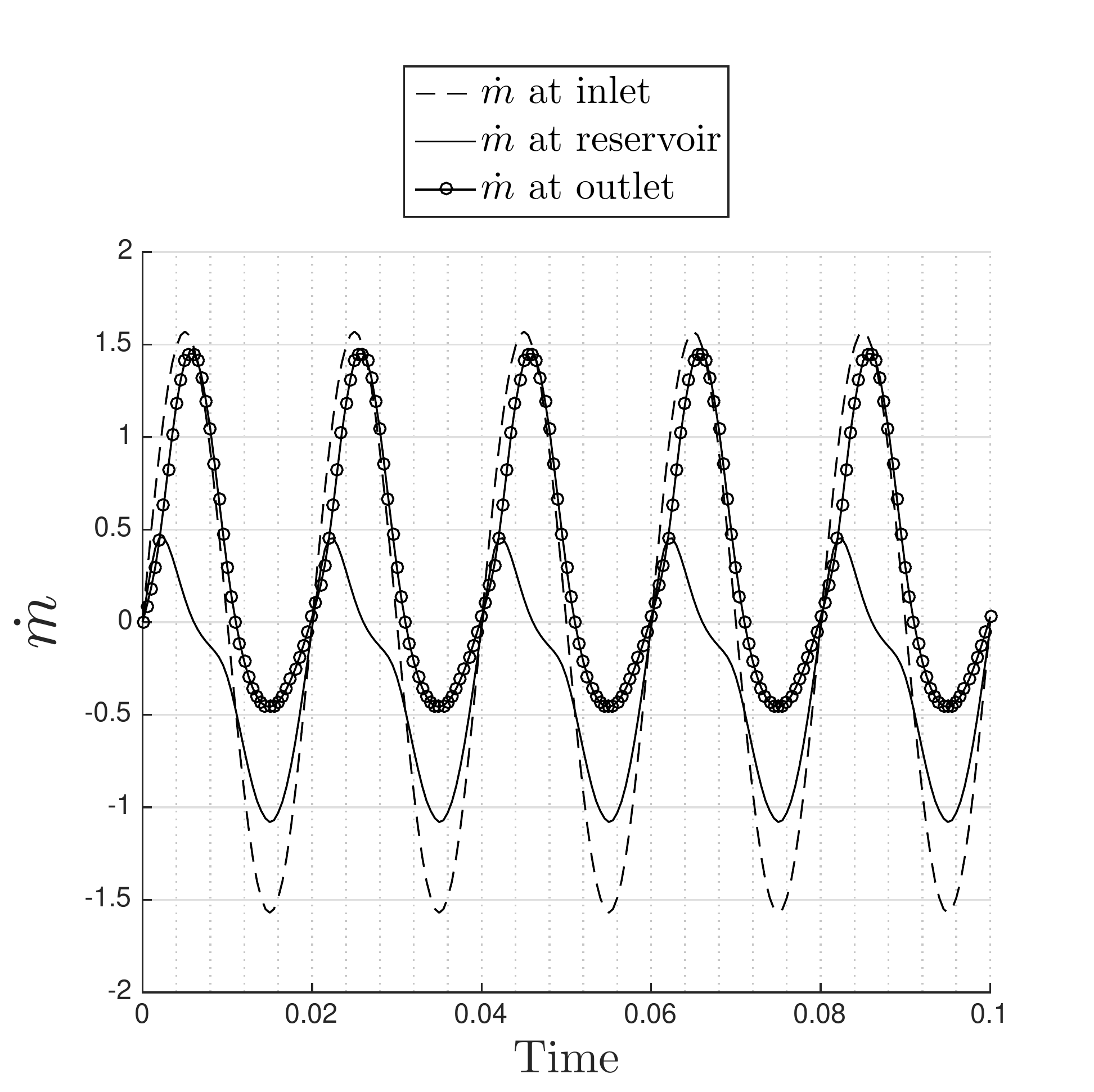}
	\caption{Mass flow through the inlet, reservoir, and outlet ports of the oscillating pump example. A positive mass flow value indicates flow entering the pump.}
	\label{fig:pump_outflow}
\end{figure*}

The design converged in $225$ iterations and the optimized layout is shown in Figure \ref{fig:pump_final_design}. Similar to the results from \cite{NSL:16}, the design exhibits a narrowing of the inflow channel, and the formation of a central reservoir from which fluid flows towards the pumping outlet during the outflow cycle. The pipe joining the bottom port to the top outlet forms naturally during the optimization process. During the inflow cycle, fluid flows from the left port towards the outlet port, while additional fluid is pulled from the bottom port towards the central reservoir. An area of recirculation then forms around the central reservoir. During the outflow phase, fluid flows from the central reservoir towards the outlet port, while fluid from the bottom port flows in a vortex like path around the left inlet; this fluid is then transported towards the outlet port during the next inflow phase. This process is shown by streamlines in Figure \ref{fig:pump_final_design_stream}. To verify the harmonic behavior of the pump for the optimized design, we modeled 5 pumping cycles of 40 time steps each; the corresponding mass flow rates for all ports are shown in Figure \ref{fig:pump_outflow}. We can observe that the optimized design maintains its oscillatory behavior. The average mass flow rate through the top port of the optimized design is $-0.328$, with the negative sign indicating that the flow leaves the domain. In contrast to the results shown by \cite{NSL:16}, we did not observe numerical artifacts in the optimized material layout of our pump design, such as isolated regions of fluid flow.

\section{Conclusions}
\label{sec:conclusions}

We presented and studied a design optimization framework based on the CutFEM method \cite{BCH+:14}, which combines the LSM, the XFEM, and face-oriented ghost-penalty methods. This framework builds upon previous optimization studies that utilized the LSM and the XFEM. The approaches presented in this paper add flexibility, accuracy, and robustness through weakly enforcing boundary conditions via a Nitsche formulation and by applying face-oriented ghost-penalty methods to the velocity and pressure fields, discretized by a generalized Heaviside-step enrichment strategy. An auxiliary indicator field was introduced to identify isolated volumes of fluid, and to enforce a constraint on the average pressure on these ``puddles'' to prevent a singular analysis problem.

The geometry of the fluid-solid interface was described by an explicit LSM where the parameters of the parametrized LSF are defined as functions of the optimization variables. The resulting optimization problem was solved in the reduced space by a nonlinear programming method. The flexibility of this scheme allowed us to define optimization variables in addition to the parameters resulting from a finite element discretization of the LSF. We demonstrated this feature by defining the location and shape of the outlet ports as optimization variables.

The accuracy of our CutFEM approach was verified quantitatively through comparison against benchmark studies, and qualitatively through application to topology optimization problems. The design optimization problems converged to intuitive designs and resembled well the results from previous 2D studies. A rigorous comparison of the proposed approach against topology optimization methods using Brinkman penalization was beyond the scope of this paper. However, such a comparison should be performed in future studies to determine which of the two approaches enjoys greater accuracy and computational efficiency.

As demonstrated by the numerical studies presented in this paper, our CutFEM approach is accurate, robust, and applicable to a broad range of design problems for low Reynolds number flows. However, the computational costs can be significant, in particular if accurate flow solutions are required. This is in parts due to the enlarged bandwidth of the linear systems caused by the face-oriented ghost-penalty formulation, but more importantly due to the lack of an adaptive mesh refinement strategy. Future research may focus on integrating mesh refinement strategies into our CutFEM framework such the boundary layer phenomena are captured accurately and efficiently. 


\section*{Acknowledgment}

The authors acknowledge the support of the National Science Foundation under grant EFRI-ODISSEI 1240374 and CBET 1246854. The second author also acknowledges support through Sandia National Laboratories under Contract Agreement 1396470. The opinions and conclusions presented in this paper are those of the authors and do not necessarily reflect the views of the sponsoring organization.



\end{document}